\documentclass[11pt]{amsart}
\usepackage{amsmath}
\usepackage{amsxtra}
\usepackage{amscd}
\usepackage{amsthm}
\usepackage{amsfonts}
\usepackage{amssymb}
\usepackage{eucal}

\input diagrams

\newtheorem{thm}{Theorem}[section]

\newtheorem{cor}[thm]{Corollary}
\newtheorem{lem}[thm]{Lemma}
\newtheorem{prop}[thm]{Proposition}

\theoremstyle{remark}
\newtheorem{remark}[thm]{Remark}

\theoremstyle{definition}

\numberwithin{equation}{section}

\newcommand{\bean}{\begin{eqnarray}}
\newcommand{\eean}{\end{eqnarray}}
\newcommand{\be}{\begin{displaymath}}
\newcommand{\ee}{\end{displaymath}}
\newcommand{\bea}{\begin{eqnarray*}}   
\newcommand{\eea}{\end{eqnarray*}}

\newcommand{\thmref}[1]{Theorem~\ref{#1}}
\newcommand{\secref}[1]{Section~\ref{#1}}
\newcommand{\lemref}[1]{Lemma~\ref{#1}}
\newcommand{\propref}[1]{Proposition~\ref{#1}}
\newcommand{\corref}[1]{Corollary~\ref{#1}}
\newcommand{\remref}[1]{Remark~\ref{#1}}

\newcommand{\nc}{\newcommand}
\nc{\on}{\operatorname}
\nc{\ch}{\mbox{ch}}
\nc{\Z}{{\mathbb Z}}
\nc{\C}{{\mathbb C}}
\nc{\pone}{{\mathbb C}{\mathbb P}^1}
\nc{\pa}{\partial}
\nc{\F}{{\mathcal F}}
\nc{\arr}{\rightarrow}
\nc{\larr}{\longrightarrow}
\nc{\al}{\alpha}
\nc{\ri}{\rangle}
\nc{\lef}{\langle}
\nc{\W}{{\mathcal W}}
\nc{\la}{\lambda}
\nc{\ep}{\epsilon}
\nc{\su}{\widehat{{\mathfrak s}{\mathfrak l}}_2}
\nc{\sw}{{\mathfrak s}{\mathfrak l}}
\nc{\g}{{\mathfrak g}}
\nc{\h}{{\mathfrak h}}
\nc{\n}{{\mathfrak n}}
\nc{\N}{\widehat{\n}}
\nc{\G}{\widehat{\g}}
\nc{\De}{\Delta}
\nc{\gt}{\widetilde{\g}}
\nc{\Ga}{\Gamma}
\nc{\one}{{\mathbf 1}}
\nc{\z}{{\mathfrak Z}}
\nc{\La}{\Lambda}
\nc{\wt}{\widetilde}
\nc{\wh}{\widehat}
\nc{\cri}{_{\kappa_c}}
\nc{\kk}{h^\vee}
\nc{\sun}{\widehat{\sw}_N}
\nc{\si}{\sigma}
\nc{\el}{\ell}
\nc{\bi}{\bibitem}
\nc{\om}{\omega}
\nc{\ol}{\overline}
\nc{\ds}{\displaystyle}
\nc{\dzz}{\frac{dz}{z}}
\nc{\Res}{\on{Res}}
\nc{\mc}{\mathcal}
\nc{\Cal}{\mathcal}
\nc{\bb}{{\mathfrak b}}
\nc{\ot}{\otimes}
\nc{\R}{{\mc R}}
\nc{\yy}{{\mc Y}}
\nc{\ga}{\gamma}

\nc{\us}{\underset}
\nc{\opl}{\oplus}
\nc{\beq}{\begin{equation}}
\nc{\Fq}{{\mathcal F}}
\nc{\Mq}{{\mathcal M}}
\nc{\Rep}{\on{Rep}}
\nc{\sssec}{\subsubsection}
\nc{\ssec}{\subsection}
\nc{\lan}{\langle}
\nc{\ran}{\rangle}

\nc{\D}{\mathcal D}
\nc{\Vect}{\on{Vect}}
\nc{\ghat}{\G}
\nc{\T}{\mc T}
\nc{\Tloc}{\T^\g_{\on{loc}}}
\nc{\vac}{|0\ran}
\nc{\Wick}{{\mb :}}
\nc{\mb}{\mathbf}
\nc{\delz}{\partial_z}
\nc{\K}{{\cali K}}
\nc{\cali}{\mathcal}
\nc{\li}{\mathfrak l}
\nc{\lt}{\widetilde{\li}}
\nc{\astar}{a^*}
\nc{\cA}{{\mc A}}
\nc{\ka}{\kappa}

\nc{\OO}{{\mc O}}
\nc{\AutO}{\on{Aut}\OO}
\nc{\DerO}{\on{Der}\OO}
\nc{\DerpO}{\on{Der}_+\OO}
\nc{\Au}{{\mc A}ut}
\nc{\mf}{\mathfrak}
\nc{\V}{{\mc V}}
\nc{\hh}{\wh{\h}}

\nc{\pp}{{\mathfrak p}}
\nc{\mm}{{\mathfrak m}}
\nc{\rr}{{\mathfrak r}}
\nc{\ket}{\rangle}
\nc{\zz}{{\mathfrak z}}
\nc{\gr}{\on{gr}}
\nc{\Spe}{\on{Spec}}
\nc{\rv}{\rho^\vee}
\nc{\can}{\on{can}}
\nc{\CC}{\on{Op}_G(D))}
\nc{\Op}{\on{Op}_G(D)}
\nc{\MOp}{\on{MOp}_G(D)}
\nc{\Db}{{\mathbb D}}
\nc{\ww}{w}

\begin{document}

\title[Wakimoto modules, opers and the center]{Lectures on Wakimoto
modules, opers and the center at the critical level}

\author{Edward Frenkel}\thanks{Partially supported by grants from the
  Packard Foundation and the NSF}

\address{Department of Mathematics, University of California,
  Berkeley, CA 94720, USA}

\date{September 2002}

\maketitle

\tableofcontents

\section*{Introduction}

Wakimoto modules are representations of affine Kac-Moody algebras in
Fock modules over infinite-dimensional Heisenberg algebras. They were
introduced in 1986 by M. Wakimoto \cite{Wak} in the case of $\su$ and
in 1988 by B. Feigin and the author \cite{FF:usp} in the general case.
Wakimoto modules have useful applications in representation theory and
conformal field theory. In particular, they have been used to
construct chiral correlation functions of the WZW models
\cite{ATY,FFR} (reproducing the Schechtman--Varchenko solutions of the
KZ equations), in the study of the Drinfeld-Sokolov reduction and
$\W$--algebras, and in the description of the center of the completed
enveloping algebra of an affine Lie algebra \cite{FF:gd,F:thesis}. The
Wakimoto modules also provide a bridge between representation theory
of affine algebras and the geometry of the semi-infinite flag manifold
\cite{FF:si}.

In these lectures we present the construction of the Wakimoto modules
from the point of view of the vertex algebra theory (see
\cite{K,vertex} for an introduction to vertex algebras). In the
original construction of \cite{FF:usp}--\cite{FF:si} the connection
with vertex algebras was not explicitly discussed. In these notes we
make this connection explicit. Let $\g$ be a simple Lie algebra and
$\ghat$ the corresponding (untwisted) affine Kac-Moody algebra. One
can associate to $\ghat$ and an invariant inner product $\ka$ on $\g$
(equivalently, a level) a vertex algebra, denoted by $V_\ka(\g)$. The
construction of Wakimoto modules amounts to constructing a
homomorphism from $V_\ka(\g)$ to a vertex algebra $M_\g \otimes
\pi^{\ka-\ka_c}_0$ associated to an infinite-dimensional Heisenberg
Lie algebra. In the case of $\g=\sw_2$ the construction is spelled out
in detail in \cite{vertex}, Ch. 10--11, where the reader is referred
for additional motivation and background. Here we generalize it to the
case of an arbitrary $\g$. We follow the original approach of
\cite{FF:si} and prove the existence of the homomorphism $V_\ka(\g)
\to M_\g \otimes \pi^{\ka-\ka_c}_0$ by cohomological methods. We
describe the cohomology class responsible for the obstruction and show
that it is equal to the class defining the affine Kac-Moody
algebra. We note that explicit formulas for the Wakimoto realization
have been given in \cite{FF:usp} for $\g=\sw_n$ and in \cite{dBF} for
general $\g$. Another proof of the existence of this realization has
been presented in \cite{FF:wak}. The construction has been extended to
twisted affine algebras in \cite{Sz}.

The Wakimoto modules may be viewed as representations of $\ghat$ which
are ``semi-infinitely induced'' from representations of its Heisenberg
subalgebra $\wh{\h}$ (extended by zero to $\wh{\bb}_-$). In contrast
with the usual induction, however, the level of the
$\wh{\h}$--module gets shifted by the so-called critical value
$\ka_c$. In particular, if we start with an $\wh{\h}$--module of level
zero (e.g., a one-dimensional module corresponding to a character of
the abelian Lie algebra $L\h$), then the resulting $\ghat$--module
will be at the critical level. Using these modules, we prove the
Kac-Kazhdan conjecture on the characters of irreducible
$\ghat$--modules of critical level (this proof appeared first in
\cite{F:thesis}).

We then extend the construction of the Wakimoto modules to a more
general context in which the Lie subalgebra $\hh$ is replaced by a
central extension of the loop algebra of the Levi subalgebra of an
arbitrary parabolic Lie subalgebra of $\g$ following the ideas of
\cite{FF:si}. Thus, we establish a ``semi-infinite parabolic
induction'' pattern for representations of affine Kac-Moody algebras,
similar to the parabolic induction for reductive groups over local
non-archimedian fields.

Next, we give a uniform construction of intertwining operators between
Wakimoto modules, the so-called screening operators (of two kinds)
following \cite{FF:kn,FF:lmp,F:thesis,FFR}.

Using the screening operators, we describe the center of the vertex
algebra $V_{\ka_c}(\g)$ at the critical level. We show that it is
canonically isomorphic to the classical $\W$--algebra associated to
the Langlands dual Lie algebra $^L\g$. The proof presented in these
lectures is different from the original proof from
\cite{FF:gd,F:thesis} in two respects. First of all, we use here the
screening operators of the second kind rather than the first kind;
their $\ka \to \ka_c$ limits are easier to study. Second, we use the
isomorphism between the Verma module of critical level with highest
weight $0$ and a certain Wakimoto module and the computation of the
associated graded of the spaces of singular vectors to estimate the
character of the center.

Finally, we show that the above classical $\W$--algebra is isomorphic
to the algebra of functions on the space of $^L G$--opers on the
formal disc. The notion of $^L G$--oper (on an arbitrary smooth curve)
was introduced by A. Beilinson and V. Drinfeld
\cite{BD},\cite{BD:opers} following the work of V. Drinfeld and
V. Sokolov \cite{DS} on the generalized KdV hierarchies. The algebra
$\on{Fun} \on{Op}_{^L G}(\Db^\times)$ of funtions on the space of $^L
G$--opers on the punctured disc carries a canonical Poisson structure
defined in \cite{DS} by means of hamiltonian reduction. At the same
time, the algebra $\on{Fun} \on{Op}_{^L G}(\Db)$ of functions on the
space of $^L G$--opers on the (unpunctured) disc carries a vertex
Poisson structure. We show that the above isomorphism between the
center of $V_{\ka_c}(\g)$ (with its canonical vertex Poisson structure
coming from the deformation of the level) and $\on{Fun} \on{Op}_{^L
G}(\Db)$ preserves vertex Poisson structures. Using this isomorphism,
we show that the center of the completed enveloping algebra of $\ghat$
at the critical level is isomorphic to $\on{Fun} \on{Op}_{^L
G}(\Db^\times)$ as a Poisson algebra. This isomorphism was
conjectured by V. Drinfeld.

These lectures are organized as follows. We begin in Sect. 1 with the
geometric construction of representations of a simple
finite-dimensional Lie algebra $\g$ using an embedding of $\g$ into a
Weyl algebra which is obtained from the infinitesimal action of $\g$
on the flag manifold. This will serve as a prototype for the
construction of Wakimoto modules presented in the subsequent
sections. In Sect. 2 we introduce the main ingredients needed for the
constructions of Wakimoto modules: the infinite-dimensional Weyl
algebra ${\mc A}^\g$, the corresponding vertex algebra $M_\g$ and the
infinitesimal action of the loop algebra $L\g$ on the space $LU$ of
formal loops to the big cell of the flag manifold of $\g$. We also
introduce the local Lie algebra ${\mc A}^\g_{\leq 1,\g}$ of
differential operators on this space of order less than or equal to
one. We show that it is a non-trivial extension of the Lie algebra of
local vector fields on $LU$ by local functionals on $LU$ and compute
the corresponding two-cocycle (most of this material has already been
presented in \cite{vertex}, Ch. 11). In Sect. 3 we prove that the
embedding of the loop algebra $L\g$ into the Lie algebra of local
vector fields on $LU$ may be lifted to an embedding of the central
extension $\ghat$ to ${\mc A}^\g_{\leq 1,\g}$. In order to do that we
need to show that the restriction of the above two-cocycle to $L\g$ is
cohomologically equivalent to the two-cocycle corresponding to its
Kac-Moody central extension (of level $\ka_c$). This is achieved by
replacing the standard cohomological Chevalley complex by a much
smaller local subcomplex (where both cocycles belong) and computing
the cohomology of the latter.

Having proved the existence of a homomorphism $\ghat_{\ka_c} \to {\mc
A}^\g_{\leq 1,\g}$, we derive in Sect. 4 the existence of a vertex
algebra homomorphism $V_{\ka_c}(\g) \to M_\g \otimes \pi_0$, where
$\pi_0$ is the commutative vertex algebra associated to $L\h$. Using
this homomorphism, we construct a $\ghat$--module structure of critical
level on the tensor product $M_\g \otimes N$, where $N$ is an
arbitrary (smooth) $L\h$--module. These are the Wakimoto modules of
critical level. We show that for a generic one-dimensional module $N$
the corresponding Wakimoto module is irreducible. This enables us to
prove the Kac--Kazhdan conjecture on the characters of irreducible
representations of critical level.

In Sect. 5 we deform the above construction of Wakimoto modules to
other levels. Thus, we construct a homomorphism of vertex algebras
$\ww_\ka: V_\ka(\g) \to M_g \otimes \pi_0^{\ka-\ka_c}$, where
$\pi_0^{\ka-\ka_c}$ is the vertex algebra associated to the Heinseberg
Lie algebra $\hh_{\ka-\ka_c}$. We describe the (quasi)conformal
structure on $M_g \otimes \pi_0^{\ka-\ka_c}$ corresponding to the
Segal-Sugawara (quasi)conformal structure on $V_\ka(\g)$. We then use
these structures to describe the Wakimoto modules in a
coordinate-independent fashion.

In Sect. 6 we generalize the construction to the case of an arbitrary
parabolic subalgebra $\pp$ of $\g$. We define the functors of
``semi-infinite parabolic induction'' from the category of smooth
representations of a central extension of $L\mm$, where $\mm$ is the
Levi subalgebra of $\pp$, to the category of $\ghat$--modules.

Sect. 7 is devoted to a detailed study of the Wakimoto modules over
$\su$. We write down explicit formulas for the homomorphism
$V_\ka(\sw_2) \to M_{\sw_2} \otimes \pi^{\ka-\ka_c}_0$ and use these
formulas to construct intertwining operators between the Wakimoto
modules. They are called screening operators, of the first and second
kind. We use these operators and the functor of parabolic induction to
construct intertwining operators between Wakimoto modules over an
arbitrary affine Kac-Moody algebra in Sect. 8.

In Sect. 9 we identify the center of the vertex algebra
$V_{\ka_c}(\g)$ with the intersection of the kernels of certain
operators, which we identify in turn with the classical $\W$--algebra
associated to the Langlands dual Lie algebra $^L\g$. We also show that
the corresponding vertex Poisson algebra structures coincide. In
Sect. 10 we define opers and Miura opers. We then show in Sect. 11
that this classical $\W$--algebra (resp., the commutative vertex
algebra $\pi_0$) is nothing but the algebra of functions on the space
of $^L G$--opers on the disc (resp., Miura $^L G$--opers on the
disc). Furthermore, the embedding of the classical $\W$--algebra into
$\pi_0$ coincides with the Miura map between the two algebras of
functions. Finally, in Sect. 12 we consider the center of the
completed universal enveloping algebra of $\ghat$ at the critical
level. We identify it with the algebra of functions on the space of
$^L G$--opers on the punctured disc. We then prove that this
identification satisfies various compatibilities. In particular, we
show that an affine analogue of the Harish-Chandra homomorphism
obtained by evaluating central elements on the Wakimoto modules is
nothing but the Miura transformation from Miura $^L G$--opers to $^L
G$--opers on the punctured disc.

\medskip

\noindent {\bf Acknowledgments.} Most of the results presented in this
paper are results of joint works with B. Feigin. I am grateful to him
for his collaboration. I thank D. Gaitsgory for valuable discussions,
in particular, for his suggestion to use the associated graded of the
spaces of singular vectors in the proof of \thmref{isom with w}.

I thank M. Szczesny for letting me use his notes of my lectures in
preparation of Sections 1--4. I am grateful to D. Kazhdan for his
valuable comments on the draft of this paper.

\section{Finite-dimensional case}

In this section we recall the realization of $\g$--modules in the
space of functions on the big cell of the flag manifold. This will
serve as a blueprint for the construction of the Wakimoto modules over
affine algebras in the following sections.

\subsection{Flag manifold}    \label{flag manifold}

Let $\g$ be a simple Lie algebra of rank $\ell$. As a vector space, it
has a Cartan decomposition
\begin{equation}    \label{cardec}
\g={\mathfrak n}_+ \oplus \h \oplus {\mathfrak n}_-,
\end{equation}
where $\h$ is the Cartan subalgebra and $\n_\pm$ are the upper and
lower nilpotent subalgebras. Let
$$
{\mathfrak b}_\pm = \h\oplus {\mathfrak n}_\pm
$$
be the upper and lower Borel subalgebras. Let $h_i = \al_i^\vee$ be
the $i$th coroot of $\g$. The set $\{ h_i \}_{i=1\ldots,\ell}$ is a
basis of $\h$. We choose generators $e_i, i=1,\ldots,\ell$, and $f_i,
i=1,\ldots,\ell$, of $\n_+$ and $\n_-$, respectively, corresponding to
the simple roots. We also choose a root basis of $\n_+$,
$\{ e^\al \}_{\al \in \De_+}$, where $\De_+$ is the set of positive
roots of $\g$, so that $[h,e^\al] = \al(h) e^\al$ for all $h \in
\h$. In particular, $e_i = e_{\al_i}, i=1,\ldots,\ell$.

Let $G$ be the connected simply-connected Lie group corresponding to
$\g$, and $N_\pm$ (respectively, $B_\pm$) the upper and lower
unipotent subgroups (respectively, Borel subgroups) of $G$
corresponding to ${\mathfrak n}_\pm$ (respectively, ${\mathfrak
  b}_\pm$).

Consider the flag manifold $G/B_-$. It has a unique open $N_+$--orbit,
the so-called big cell $U=N_+\cdot [1] \subset G/B_-$, which is
isomorphic to $N_+$. Since $N_+$ is a unipotent Lie group, the
exponential map $\n_+ \arr N_+$ is an isomorphism. Therefore $N_+$ is
isomorphic to the vector space $\n_+$. Thus, $N_+$ is isomorphic to
the affine space ${\mathbb A}^{|\Delta_+|}$, and the space $\C[N_+]$
of regular functions on $N_+$ is isomorphic to a free polynomial
algebra. We will call a system of coordinates $\{ y_\al \}_{\al \in
\De_+}$ on $N_+$ {\em homogeneous} if $h \cdot y_\al = - \al(h) y_\al$
for all $h \in \h$. In what follows we will consider only homogeneous
coordinate systems on $N_+$.

\medskip

\begin{remark}
  Note that in order to define $U$ it is sufficient to choose only a
  Borel subgroup $B_+$ of $G$. Then $N_+ = [B_+,B_+]$ and $U$ is the
  open $N_+$--orbit in the flag manifold defined as the variety of all
  Borel subgroups of $G$ (so $U$ is an $N_+$--torsor). All
  constructions in this paper make sense with the choice of $B_+$
  only, i.e., without making the choice of $H$ and $B_-$. However, to
  simplify the exposition we will fix a Cartan subgroup $H \subset
  B_+$ as well.\qed
\end{remark}

\medskip

The action of $G$ on $G/B_-$ gives us a map from $\g$ to the
Lie algebra of vector fields on $G/B_-$, and hence on its open subset
$U \simeq N_+$. Thus, we obtain a Lie algebra homomorphism $\g \to
\on{Vect} N_+$.

This homomorphism may be described explicitly as follows. Let
$G^\circ$ denote the dense open submanifold of $G$ consisting of
elements of the form $g_{+}g_{-}, g_{+} \in N_+, g_{-} \in B_{-}$
(note that such an expression is necessarily unique since $B_{-} \cap
N_+ = 1)$. In other words, $G^\circ = p^{-1}(N_+)$, where $p$ is the
projection $G \rightarrow G/B_-$. Given $a \in \g$, consider the
one-parameter subgroup $\gamma(\epsilon)=\exp(\epsilon a)$ in $G$.
Since $G^\circ$ is open and dense in $G$, $\gamma(\epsilon)x \in
G^\circ$ for $\epsilon$ in the formal neighborhood of $0$, so we can
write
$$\gamma(\epsilon)=Z_{+}(\epsilon)Z_{-}(\epsilon), \qquad
Z_{+}(\epsilon) \in N_+, \quad Z_{-}(\epsilon) \in B_{-}.$$
The factor
$Z_{+}(\epsilon)$ just expresses the projection of the subgroup
$\gamma(\epsilon)$ onto $N_+ \simeq U \subset G/B_-$ under the map
$p$. Then the vector field $\xi_a$ (equivalently, a derivation of
$\C[N_+]$) corresponding to $a$ is given by the formula
\begin{equation}    \label{def of xia}
(\xi_a f)(x) = \left. \left(\frac{d}{d\ep} f(Z_+(\ep))\right)
\right|_{\ep = 0}.
\end{equation}

To write a formula for $\xi_a$ in more concrete terms, we choose a
faithful representation $V$ of $\g$ (say, the adjoint representation).
Since we only need the $\ep$--linear term in our calculation, we can
and will assume that $\ep^2=0$. Considering $x \in N_+$ as a matrix
whose entries are polynomials in the coordinates $y_{\alpha}, \alpha
\in \Delta_{+}$, which expresses a generic element of $N_+$ in
$\on{End} V$, we have
\begin{equation}    \label{split a matrix}
        (1+\epsilon a) x = Z_{+}(\epsilon)Z_{-}(\epsilon).
\end{equation}
We find from this formula that $Z_{+}(\epsilon) = x + \ep Z_+^{(1)}$,
where $Z_+^{(1)} \in \n_+$, and $Z_{-}(\epsilon) = 1 + \ep Z_-^{(1)}$,
where $Z_-^{(1)} \in \bb_{-}$. Therefore we obtain from formulas
\eqref{def of xia} and \eqref{split a matrix} that
\begin{equation}    \label{expl xia}
\xi_a \cdot x = x(x^{-1}ax)_+,
\end{equation}
where $z_+$ denotes the projection of an element $z \in \g$ onto
$\n_+$ along $\bb_-$.

\subsection{The algebra of differential operators}

The algebra $\D(U)$ of differential operators on $U$ is isomorphic to
the Weyl algebra with generators $\{ y_\al, \pa/\pa y_\al \}_{\al \in
\De_+}$, and the standard relations
$$
\left[\frac{\pa}{\pa y_\al},y_\beta\right] = \delta_{\al,\beta}, \qquad
\left[\frac{\pa}{\pa y_\al},\frac{\pa}{\pa y_\beta}\right] =
[y_\al,y_\beta] = 0.
$$

The algebra $\D(U)$ has a natural filtration $\{ \D_{\leq i}(U) \}$ by
the order of the differential operator. In particular, we have an
exact sequence
\begin{equation}\label{Func and Vect}
0 \to \on{Fun} U\to \D_{\leq 1}(U) \to \on{Vect}U\to 0,
\end{equation}
where $\on{Fun} U \simeq \on{Fun} N_+$ denotes the ring of regular
functions on $U$, and $\on{Vect} U$ denotes the Lie algebra of vector
fields on $U$. This sequence has a canonical splitting: namely, we
lift $\xi\in\on{Vect} U$ to the unique first order differential
operator $D_\xi$ whose symbol equals $\xi$ and which kills the
constant functions, i.e., such that $D_\xi \cdot 1 = 0$. Using this
splitting, we obtain an embedding $\g \arr \D_{\leq 1}(N_+)$, and
hence the structure of a $\g$--module on the space of functions
$\on{Fun} N_+ = \C[y_{\al}]_{\al \in \De_+}$.

\subsection{Verma modules and contragredient Verma modules}
\label{contra}

By construction, the action of $\n_+$ on $\on{Fun} N_+$ satisfies $e^\al
\cdot y_\beta = \delta_{\al,\beta}$. Therefore, for any $A \in
\on{Fun} N_+$ there exists $P \in U(\n_+)$ such that $P \cdot A = 1$.
Consider the pairing $U(\n_+) \times \on{Fun} N_+ \to \C$ which maps
$(P,A)$ to the value of the function $P \cdot A$ at the identity
element of $N_+$. This pairing is $\n_+$--invariant. The
Poincar\'e-Birkhoff-Witt basis $\{ e^{\al(1)} \ldots e^{\al(k)} \}$
(with respect to some ordering on $\Delta_+$) of $U(\n_+)$ and the
monomial basis $\{ y_{\al(1)} \ldots y_{\al(k)} \}$ of $\on{Fun} N_+$ are
dual to each other with respect to this pairing. Moreover, both
$U(\n_+)$ and $\on{Fun} N_+$ are graded by the root lattice of $\g$ (under
the action of the Cartan subalgebra $\h$), and this pairing identifies
each graded component of $\on{Fun} N_+$ with the dual of the corresponding
graded component of $U(\n_+)$. Therefore the $\n_+$--module $\on{Fun} N_+$
is isomorphic to the restricted dual $U(\n_+)^\vee$ of $U(\n_+)$
(i.e., $\on{Fun} N_+$ is the direct sum of the dual spaces to the
homogeneous components of $U(\n_+)$).

Recall the definition of the Verma module and the contragredient Verma
modules.  For each $\chi \in \h^*$, consider the one-dimensional
representation $\C_\chi$ of $\bb_\pm$ on which $\h$ acts according to
$\chi$, and $\n_\pm$ acts by $0$. The {\em Verma module} $M_\chi$
\index{Verma module} with highest weight $\chi \in \h^*$ is the
induced module
\begin{equation*}    \label{defverma}
M_\chi \overset{\on{def}}{=} U(\g) \otimes_{U(\bb_+)} \C_\chi.
\end{equation*}
The Cartan decomposition \eqref{cardec} gives us the isomorphisms
$$
U(\g) \simeq U(\n_-) \otimes U(\bb_+), \quad \quad U(\g) \simeq
U(\n_+) \otimes U(\bb_-).
$$
Therefore, as an $\n_-$--module, $M_\chi \simeq U(\n_-)$.

The {\em contragredient Verma module} $M^*_\chi$ with highest weight
$\chi \in \h^*$ is defined as the coinduced module
\index{representation, coinduced}
$$
M^*_\chi = \on{Hom}^{\on{res}}_{U(\bb_-)}(U(\g),\C_\chi).
$$
Here $U(\g)$ is considered as a $U(\bb_-)$--module with respect to
the right action, and we consider only those homomorphisms $U(\g) \to
\C_\chi$ which belong to $U(\bb_-)^* \otimes U(\n_+)^\vee \subset
U(\g)^*$. Thus, as an $\n_+$--module, $M^*_\chi \simeq U(\n_+)^\vee$.

\subsection{Identification of $\on{Fun} N_+$ with $M^*_\chi$}

The module $M^*_0$ is isomorphic to $\on{Fun} N_+$ with its $\g$--module
structure defined above. Indeed, the vector $1 \in \on{Fun} N_+$ is
annihilated by $\n_+$ and has weight $0$ with respect to $h$. Hence
there is a non-zero homomorphism $\on{Fun} N_+ \to M^*_0$ sending $1 \in
\on{Fun} N_+$ to a non-zero vector $v^*_0 \in M^*_\chi$ of weight $0$.
Since both $\on{Fun} N_+$ and $M^*_\chi$ are isomorphic to $U(\n_+)^\vee$
as $\n_+$--modules, this homomorphism is an isomorphism.

Now we identify the module $M^*_\chi$ with an arbitrary weight $\chi$
with $\on{Fun} N_+$, where the latter is equipped with a modified action of
$\g$.

Recall that we have a canonical lifting of $\g$ to $\D_{\leq 1}(N_+)$,
$a \mapsto \xi_a$. But this lifting is not unique. We can modify it by
adding to each $\xi_a$ a function $\phi_a \in \on{Fun} N_+$ so that
$\phi_{a+b} = \phi_a + \phi_b$. One readily checks that the modified
differential operators $\xi_a + \phi_a$ satisfy the commutation
relations of $\g$ if and only if the linear map $\g \arr \on{Fun} N_+$
given by $a \mapsto \phi_a$ is a one-cocycle of $\g$ with coefficients
in $\on{Fun} N_+$.

Each such lifting gives $\on{Fun} N_+$ the structure of a
$\g$--module. Let us impose the extra condition that the modified
action of $\h$ on $V$ remains diagonalizable. This means that $\phi_h$
is a constant function on $N_+$ for each $h \in \h$, and therefore our
one-cocycle should be $\h$--invariant: $\phi_{[h,a]} = \xi_h \cdot
\phi_a$, for all $h \in \h, a \in \g$. It is easy to see that the
space of $\h$--invariant one-cocycles of $\g$ with coefficients in
$\on{Fun} N_+$ is canonically isomorphic to the first cohomology of
$\g$ with coefficients in $\on{Fun} N_+$, i.e., $H^1(\g,\on{Fun} N_+)$
(see \cite{vertex}, \S~10.2.6). By Shapiro's lemma (see \cite{Fuchs},
\S~5.4)
\begin{eqnarray*}
H^1(\g,\on{Fun} N_+) = H^1(\g,M_0^*) \simeq H^1({\mathfrak b}_-,\C_0)
=({\mathfrak b}_-/[{\mathfrak b}_-,{\mathfrak b}_-])^* \simeq \h^*.
\end{eqnarray*}
Thus, for each $\chi \in \h^*$ we obtain an embedding $\rho_\chi: \g
\hookrightarrow \D_{\leq 1}(N_+)$ and hence the structure of an
$\h^*$--graded $\g$--module on $\on{Fun} N_+$. Let us analyze this
$\g$--module in more detail.

We have $\xi_h \cdot y_\al = -\al(h) y_\al, \al \in \De_+$, so the
weight of any monomial in $\on{Fun} N_+$ is equal to a sum of negative
roots. Since our one-cocycle is $\h$--invariant, we obtain $\xi_h
\cdot \phi_{e^{\al}} = \phi_{[h,e^\al]} = \al(h) \phi_{e^{\al}}$, $\al
\in \De_+$, so the weight of $\phi_{e^{\al}}$ has to be equal to the
positive root $\al$. Therefore $\phi_{e^\al} = 0$ for all $\al \in
\De_+$. Thus, the action of $\n_+$ on $\on{Fun} N_+$ is not modified. On
the other hand, by construction, the action of $h \in \h$ is modified
by $h \mapsto h + \chi(h)$. Therefore the vector $1 \in \on{Fun} N_+$ is
still annihilated by $\n_+$, but now it has weight $\chi$ with respect
to $h$. Hence there is a non-zero homomorphism $\on{Fun} N_+ \to M^*_\chi$
sending $1 \in \on{Fun} N_+$ to a non-zero vector $v^*_\chi \in M^*_\chi$
of weight $\chi$. Since both $\on{Fun} N_+$ and $M^*_\chi$ are isomorphic
to $U(\n_+)^\vee$ as $\n_+$--modules, this homomorphism is an
isomorphism. Thus, under the modified action obtained via the lifting
$\rho_\chi$, the $\g$--module $\on{Fun} N_+$ is isomorphic to the
contragredient Verma module $M^*_\chi$.

\subsection{Explicit formulas}    \label{sec on formulas}

Choose a basis $\{ J^a \}_{a=1,\ldots,\dim \g}$ of $\g$. Under the
homomorphism $\rho_\chi$, we have
\begin{equation}    \label{Pa}
J^a \mapsto P_a\left( y_\al,\frac{\pa}{\pa y_\al} \right) +
f_a(y_\al),
\end{equation}
where $P_a$ is a polynomial in the $y_\al$'s and $\pa/\pa y_\al$'s of
degree one in the $\pa/\pa y_\al$'s, which is independent of $\chi$,
and $f_a$ is a polynomial in the $y_\al$'s only which depends on
$\chi$.

Let $e_i, h_i, f_i, i=1\ldots,\ell$, be the generators of $\g$. Using
formula \eqref{expl xia} we find the following explicit formulas:
\begin{align}    \label{formulas1}
\rho_\chi(e_i) &= \frac{\pa}{\pa y_{\al_i}} + \sum_{\beta \in \Delta_+}
P^i_\beta(y_\al) \frac{\pa}{\pa y_\beta}, \\ \label{formulas2}
\rho_\chi(h_i) &= - \sum_{\beta \in \Delta_+} \beta(h_i) y_\beta
\frac{\pa}{\pa y_\beta} + \chi(h_i), \\ \label{formulas3}
\rho_\chi(f_i) &= \sum_{\beta \in \Delta_+}
Q^i_\beta(y_\al) \frac{\pa}{\pa y_\beta} +  \chi(h_i) y_{\al_i},
\end{align}
for some polynomials $P^i_\beta, Q^i_\beta$ in $y_\al, \al \in \De_+$.

In addition, we have a Lie algebra anti-homomorphism $\rho^R: \n_+ \to
\D_{\leq 1}(N_+)$ which corresponds to the {\em right} action of $\n_+$
on $N_+$. The differential operators $\rho^R(x), x \in \n_+$, commute
with the differential operators $\rho_\chi(x'), x' \in \n_+$ (though
their commutation relations with $\rho_\chi(x'), x' \not\in \n_+$, are
complicated in general). We have
$$
\rho^R(e_i) = \frac{\pa}{\pa y_{\al_i}} + \sum_{\beta \in \Delta_+}
P^{R,i}_\beta(y_\al) \frac{\pa}{\pa y_\beta}
$$
for some polynomials $P^{R,i}_\beta, Q^i_\beta$ in $y_\al, \al \in \De_+$.

\section{The case of affine algebras}

\subsection{The infinite-dimensional Weyl algebra}    \label{infdim
Weyl}

Our goal is to generalize the above construction to the case of affine
Kac-Moody algebras.  Let again $U$ be the open $N_+$--orbit of the
flag manifold of $G$, which we identify with the group $N_+$ and hence
with the Lie algebra $\n_+$.  Consider the formal loop space
$LU=U((t))$ as a complete topological vector space with the basis of
open neighborhoods of zero formed by the subspaces $U \otimes t^N
\C[[t]] \subset LU, N \in \Z$. Thus, $LU$ is an affine ind-scheme
$$
LU = \underset{\longrightarrow}\lim \; U \otimes t^N \C[[t]], \qquad
N<0.
$$
Using the coordinates $y_\al, \al \in \De_+$, on $U$, we can write
$$U \otimes t^N \C[[t]] \simeq \left\{ \sum_{n\geq N} y_{\al,n} t^n
\right\}_{\al \in \De_+} = \on{Spec} \C[y_{\al,n}]_{n\geq N}.$$
Therefore we obtain that the ring of functions on $LU$, denoted by
$\on{Fun} LU$, is the inverse limit of the rings $\C[y_{\al,n}]_{\al
  \in \De_+,n\geq N}, N<0$, with respect to the natural surjective
homomorphisms $$s_{N,M}: \C[y_{\al,n}]_{\al \in \De_+,n\geq N} \to
\C[y_{\al,n}]_{\al \in \De_+,n\geq M}, \qquad N<M,$$ such that
$y_{\al,n} \mapsto 0$ for $N \leq n < M$ and $y_{\al,n} \mapsto
y_{\al,n}, n\geq M$. This is a complete topological ring, with the
basis of open neighborhoods of $0$ given by the ideals generated by
$y_{\al,n}, n<N$, i.e., the kernels of the homomorphisms
$s_{\infty,M}: \on{Fun} LU \to \on{Fun} U \otimes t^N \C[[t]]$.

A vector field on $LU$ is by definition a continuous linear
endomorphism $\xi$ of $\on{Fun} LU$ which satisfies the Leibniz rule:
$\xi(fg)= \xi(f)g + f\xi(g)$.  In other words, a vector field is a
linear endomorphism $\xi$ of $\on{Fun} LU$ such that for any $M < 0$
there exist $N\leq M$ and a derivation $$\xi_{N,M}: \C[y_{\al,n}]_{\al
  \in \De_+,n\geq N} \to \C[y_{\al,n}]_{\al \in \De_+,n\geq M}$$ which
satisfies
$$s_{\infty,M}(\xi \cdot f) = \xi_{N,M} \cdot s_{\infty,N}(f)$$
for
all $f \in \on{Fun} LU$. The space of vector fields is naturally a
topological Lie algebra, which we denote by $\on{Vect} LU$.

\begin{remark}
More concretely, an element of $\on{Fun} LU$ may be represented
as a (possibly infinite) series
$$
P_0 + \sum_{n<0} P_n y_{\al,n},
$$
where $P_0 \in \C[y_{\al,n}]_{\al \in \De_+,n\geq 0}$, and the $P_n$'s
are arbitrary (finite) polynomials in $y_{\al,m}, m \in \Z$.

The Lie algebra $\on{Vect} LU$ may also be described as
follows. Identify the tangent space $T_0 LU$ to the origin in $LU$
with $LU$, equipped with the structure of a complete topological
vector space. Then $\on{Vect} LU$ is isomorphic to the completed
tensor product of $\on{Fun} LU$ and $LU$. This means that vector
fields on $LU$ can be described more concretely as series
$$\sum_{n \in \Z} P_n \frac{\pa}{\pa y_{\al,n}},$$ where $P_n
\in \on{Fun} LU$ satisfies the following property: for each $M\geq 0$,
there exists $N\leq M$ such that each $P_n, n\leq N$, lies in the
ideal generated by the $y_{\al,m}, \al \in \De_+, m \leq M$. The
commutator between two such series is computed in the standard way
(term by term), and it is again a series of the above form.\qed
\end{remark}

\subsection{Action of $L\g$ on $LU$}    \label{action of Lg}

We have a natural Lie algebra homomorphism $$\wh{\rho}: L\g \to
\on{Vect} LU,$$
which may be described explicitly by the formulas that
we obtained in the finite-dimensional case, in which we replace the
ordinary variables $y_{\alpha}$ with the ``loop variables''
$y_{\alpha,n}$. More precisely, we have the following analogue of
formula \eqref{split a matrix},
\begin{equation*}
        (1+\epsilon A \otimes t^m) x(t) = Z_{+}(\epsilon)
        Z_{-}(\epsilon)
\end{equation*}
where $x \in N_+((t)), Z_{+}(\epsilon) = x(t) + \ep Z_+^{(1)}$,
$Z_+^{(1)} \in \n_+((t))$, and $Z_{-}(\epsilon) = 1 + \ep Z_-^{(1)}$,
$Z_-^{(1)} \in \bb_{-}((t))$. As before, we choose a faithful
finite-dimensional representation $V$ of $\g$ and consider $x(t)$ as
a matrix whose entries are Laurent power series in $t$ with
coefficients in the ring of polynomials in the coordinates
$y_{\alpha,n}, \alpha \in \Delta_{+}, n \in \Z$, expressing a generic
element of $N_+((t))$ in $\on{End} V((t))$. We define $\wh{\rho}$ by
the formula
$$
\wh\rho(a \otimes t^m) \cdot x(t) = Z^{(1)}_+.
$$
Then we have the following analogue of formula \eqref{expl xia},
\begin{equation}    \label{expl xia1}
\wh\rho(a \otimes t^m) \cdot x(t) = x(t)\left( x(t)^{-1}(a \otimes t^m)
x(t) \right)_+,
\end{equation}
where $z_+$ denotes the projection of an element $z \in \g((t))$ onto
$\n_+((t))$ along $\bb_-((t))$.

This formula implies that for any $a \in \g$ the series
$$
\wh\rho(a(z)) = \sum_{n \in \Z} \wh\rho(a \otimes t^n) z^{-n-1}
$$
may be obtained from the formula for $\rho_0(a) = \xi_a$ by the
substitution
\begin{align*}
y_\al & \mapsto \sum_{n \in \Z} y_{\al,n} z^{n}, \\ \frac{\pa}{\pa
y_\al} & \mapsto \sum_{n \in \Z} \frac{\pa}{y_{\al,n}} z^{-n-1}.
\end{align*}

\subsection{The Weyl algebra}    \label{local filtr} \label{hom of VAs}

Let $\cA^{\g}$ be the Weyl algebra with generators $$a_{\al,n} =
\frac{\partial}{\partial y_{\al,n}}, \quad a^*_{\al,n} = y_{\al,-n},
\qquad \al \in \De_+, n \in \Z,$$ and relations
\begin{equation}    \label{commina}
[a_{\al,n},a_{\beta,m}^*] = \delta_{\al,\beta}
\delta_{n,-m}, \hskip.3in
[a_{\al,n},a_{\beta,m}] = [a_{\al,n}^*,a_{\beta,m}^*] = 0.
\end{equation}
Introduce the generating functions
\begin{eqnarray}    \label{az}
a_\al(z)&=&\sum_{n\in\Z} a_{\al,n} z^{-n-1}, \\  \label{a*z}
a_\al^*(z)&=&\sum_{n\in\Z} a_{\al,n}^* z^{-n}.
\end{eqnarray}

Consider a topology on $\cA^\g$ in which the basis of open
neighborhoods of $0$ is formed by the left ideals $I_{N,M}, N,M \in
\Z$, generated by $a_{\al,n}, \al \in \Delta_+, n\geq N$, and
$a^*_{\al,m}, \al \in \De_+, m \geq M$. The {\em completed Weyl
  algebra} $\wt{\cA}^\g$ is by definition the completion of $\cA^\g$
with respect to this topology. The algebra $\wt{\cA}^\g$ should be
thought of as the algebra of differential operators on $LU$ (see
\cite{vertex}, Ch. 11, for more details).

In concrete terms, elements of $\wt{\cA}^\g$ may be viewed as
arbitrary series of the form
\begin{equation}    \label{explat}
\sum_{n>0} P_n a_n + \sum_{m>0} Q_m a^*_m + R, \quad \quad P_n, Q_m,
R \in \cA^\g.
\end{equation}

Let $\cA_0^\g$ be the (commutative) subalgebra of $\cA^\g$ generated
by $a^*_{\al,n}, \al \in \De_+, n \in \Z$, and $\wt{\cA}_0$ its
completion in $\wt{\cA}$. Next, let $\cA^\g_{\leq 1}$ be the subspace
of $\cA^\g$ spanned by the products of elements of $\cA^\g_0$ and the
generators $a_{\al,n}$. Denote by $\wt{\cA}^\g_{\leq 1}$ its
completion in $\wt{\cA}^\g$. Thus, $\wt{\cA}^\g_{\leq 1}$ consists of
all elements $P$ of $\wt{\cA}^\g$ with the property that $$P \;
\on{mod} I_{N,M} \in \cA_{\leq 1}^\g \on{mod} I_{N,M}, \qquad \forall
N,M \in \Z.$$

Here is a more concrete description of $\wt{\cA}^\g_0$ and
$\wt{\cA}^\g_{\leq 1}$ using the realization of $\wt{\cA}^\g$ by the
series of the form \eqref{explat}. The space $\wt{\cA}^\g_0$ consists
of series of the form \eqref{explat}, where all $P_n = 0$, and $Q_m, R
\in \cA^\g_0$. In other words, these are the series which do not
contain $a_n, n \in \Z$. The space $\wt{\cA}^\g_{\leq 1}$ consists of
elements of the form \eqref{explat}, where $P_n \in \cA^\g_0$, and
$Q_m, R \in \cA^\g_{\leq 1}$. The following assertion is proved in
Proposition 11.1.6 of \cite{vertex}.

\begin{prop}    \label{three parts}

\hfill

\begin{enumerate}
\item[(1)] $\wt{\cA}^\g_{\leq 1}$ is a Lie algebra, and $\wt{\cA}^\g_0$ is
its ideal.

\item[(2)] $\wt{\cA}^\g_0 \simeq \on{Fun} LU$.

\item[(3)] $\wt{\cA}^\g_{\leq 1}/\wt{\cA}^\g_0 \simeq \on{Vect} LU$ as Lie
algebras.
\end{enumerate}
\end{prop}

According to \propref{three parts}, $\wt{\cA}^\g_{\leq 1}$ is an
extension of the Lie algebra $\on{Vect} LU$ by its module
$\wt{\cA}^\g_0$,
\begin{equation}\label{local sequence}
0\to \on{Fun} LU \to \wt{\cA}^\g_{\leq 1} \to \on{Vect} LU
\to 0.
\end{equation}
This extension is however different from the standard (split)
extension defining the Lie algebra of the usual differential operators
on $LU$ of order less than or equal to $1$. In particular, our
extension \eqref{local sequence} is non-split (see Lemma 11.1.9 of
\cite{vertex}). Therefore we cannot expect to lift the homomorphism
$\wh\rho: L\g \rightarrow \Vect LU$ to a homomorphism $L\g \rightarrow
\wt{\cA}^\g_{\leq 1}$, as in the finite-dimensional case.
Nevertheless, we will show that the homomorphism $\wh\rho$ may be
lifted to a homomorphism $\ghat_\ka \rightarrow \wt{\cA}^\g_{\leq 1}$,
where $\ghat_\ka$ is the universal one-dimensional central extension
of $L\g$ defined as follows.

Fix an invariant inner product $\kappa$ on $\g$ and let $\ghat_\ka$
denote the central extension of $L\g = \g \otimes \C((t))$ with the
commutation relations
\begin{equation}    \label{KM rel}
[A \otimes f(t),B \otimes g(t)] = [A,B] \otimes f(t) g(t) - (\kappa(A,B)
\on{Res} f dg) K,
\end{equation}
where $K$ is the central element. The Lie algebra $\ghat_\ka$ is the
{\em affine Kac-Moody algebra} associated to $\ka$.

We will begin by showing that the image of the embedding $L\g
\rightarrow \Vect LU$ belongs to a Lie subalgebra of ``local'' vector
fields $\Tloc \subset \Vect LU$. This observation will allow us to
replace the extension \eqref{local sequence} by its ``local'' part.

\subsection{Heisenberg vertex algebra}    \label{hva}

Let $M_{\g}$ be the Fock representation of $\cA^{\g}$ generated by a
vector $\vac$ such that
$$
a_{\al,n} \vac = 0, \quad n\geq 0; \qquad a^*_{\al,n} \vac = 0, \quad
n>0.
$$
It is clear that the action of $\cA^{\g}$ on $M_\g$ extends to a
continuous action of the topological algebra $\wt{\cA}^\g$ (here we
equip $M_\g$ with the discrete topology). Moreover, $M_{\g}$ carries
a vertex algebra structure defined as follows (see Definition 1.3.1 of
\cite{vertex} for the definition of vertex algebra).

\begin{enumerate}
\item[$\bullet$] Gradation: $\deg a_{\al,n} = \deg a^*_{\al,n} = -n,
\deg \vac = 0$;

\item[$\bullet$] Vacuum vector: $\vac$;

\item[$\bullet$] Translation operator: $T \vac = 0, [T,a_{\al,n}] = -
n a_{\al,n-1}, [T,a^*_{\al,n}] = - (n-1) a^*_{\al,n-1}$.

\item[$\bullet$] Vertex operators:
$$
Y(a_{\al,-1}\vac,z) = a_\al(z), \quad \quad
Y(a_{\al,0}^*\vac,z)=a_\al^*(z),
$$
\begin{multline*}
Y(a_{\al_1,n_1} \ldots a_{\al_k,n_k} a^*_{\beta_1,m_1} \ldots
a^*_{\beta_l,m_l} \vac,z) = \prod_{i=1}^k \frac{1}{(-n_i-1)!}
\prod_{j=1}^l \frac{1}{(-m_j)!} \cdot \\ \Wick \delz^{-n_1-1}
a_{\al_1}(z) \dots \delz^{-n_k-1} a_{\al_k}(z) \delz^{-m_1}
a^*_{\beta_1}(z) \dots \delz^{-m_l} a^*_{\beta_l}(z) \Wick \, .
\end{multline*}
\end{enumerate}
(see Lemma 10.3.8 of \cite{vertex}).

Set
$$
U(M_\g) \overset{\on{def}}{=} (M_\g \otimes \C((t)))/\on{Im}(T \otimes
1 + \on{Id} \otimes \pa_t).
$$
We have a linear map $U(M_\g) \to \wt{\cA}^\g$ sending $A \otimes
f(t)$ to $\on{Res} Y(A,z) f(z) dz$. According to \cite{vertex},
\S~3.5, this is a Lie algebra, which may be viewed as the
completion of the span of all Fourier coefficients of vertex operators
from $M_\g$. Moreover, the above map is a homomorphism of Lie algebras
(actually, an embedding in this case).

Note that $U(M_\g)$ is not an algebra. For instance, it contains the
generators $a_{\al,n}, a^*_{\al,n}$ of the Heisenberg algebra, but
does not contain monomials in these generators of degree greater than
one. However, we will only need the Lie algebra structure on $U(M_\g)$.

The elements of $\wt{\cA}^\g_0 = \on{Fun} LU$ which lie in the image
of $U(M_\g)$ are called ``local functionals'' on $LU$. The elements of
$\wt{\cA}^\g$ which belong to $U(M_\g)$ are given by (possibly
infinite) linear combinations of Fourier coefficients of normally
ordered polynomials in $a_\al(z), a^*_\al(z)$ and their derivatives.
We refer to them as ``local'' elements of $\wt{\cA}^\g$.

\subsection{Coordinate independent definition of $M_\g$}

The above definition of the vertex algebra $M_\g$ referred to a
particular system of coordinates $y_\al, \al \in \De_+$, on the group
$N_+$. If we choose a different coordinate system $y'_\al, \al \in
\De_+$, on $N_+$, we obtain another Heisenberg algebra with generators
$a'_{\al,n}$ and $a^{*}_{\al,n}{}'$ and a vertex algebra $M'_\g$.
However, the vertex algebras $M'_\g$ and $M_\g$ are canonically
isomorphic to each other. In particular, it is easy to express the
vertex operators $a'_\al(z)$ and $a^{*}_\al{}'(z)$ in terms of
$a_\al(z)$ and $a^*_\al(z)$. Namely, if $y'_\al = F_\al(y_\beta)$,
then
\begin{align*}
a'_\al(z) & \mapsto \sum_{\gamma \in \De_+} \Wick \pa_{y_\gamma}
F_\al(a^*_\beta(z)) \; a_\gamma(z) \Wick, \\
a^{*}_\al{}'(z) & \mapsto F_\al(a^*_\beta(z)).
\end{align*}

It is also possible to define $M_\g$ without any reference to a
coordinate system on $N_+$. Namely, we set
$$
M_\g = \on{Ind}_{\n_+[[t]]}^{\n_+((t))} \on{Fun}(N_+[[t]]) \simeq
U(\n_+ \otimes t^{-1}\C[t^{-1}]) \otimes \on{Fun}(N_+[[t]]),
$$
where $\on{Fun}(N_+[[t]])$ is the ring of regular functions on the
proalgebraic group $N_+[[t]]$, considered as an $\n_+[[t]]$--module.
If we choose a coordinate system $\{ y_\al \}_{\al \in \De_+}$ on
$N_+$, then we obtain a coordinate system $\{ y_{\al,n} \}_{\al \in
  \De_+,n\geq 0}$ on $N_+[[t]]$. Then $u \otimes P(y_{\al,n}) \in
M_\g$, where $u \in U(\n_+ \otimes t^{-1}\C[t^{-1}])$ and
$P(y_{\al,n}) \in \on{Fun}(N_+[[t]]) = \C[y_{\al,n}]_{\al \in
  \De_+,n\geq 0}$, corresponds to $u \cdot P(a^*_{\al,n})$ in our
previous description of $M_\g$.

It is straightforward to define a vertex algebra structure on $M_\g$
(see \cite{FF:wak}, \S~2). Namely, the vacuum vector of $M_\g$ is
the vector $1 \otimes 1 \in M_\g$. The translation operator $T$ is
defined as the operator $-\pa_t$, which naturally acts on
$\on{Fun}(N_+[[t]])$ as well as on $\n_+((t))$ preserving $\n_+[[t]]$.
Next, we define the vertex operators corresponding to the elements
of $M_\g$ of the form $x_{-1} \vac$, where $x \in \n_+$, by the formula
$$
Y(x_{-1}\vac,z) = \sum_{n \in \Z} x_n z^{-n-1},
$$
where $x_{n} = x \otimes t^{n}$, and we consider its action on
$M_\g$ viewed as the induced representation of $\n_+((t))$. We also
need to define the vertex operators
$$
Y(P\vac,z) = \sum_{n \in \Z} P_{(n)} z^{-n-1}
$$
for $P \in \on{Fun}(N_+[[t]])$. The corresponding linear operators
$P_{(n)}$ are completely determined by their action on $\vac$:
\begin{align*}
P_{(n)} \vac &= 0, \qquad n\geq 0, \\
P_{(n)} \vac &= \frac{1}{(-n-1)!} T^{-n-1} P\vac,
\end{align*}
their mutual commutativity and the following commutation relations
with $\n_+((t))$:
$$
[x_m,P_{(k)}] = \sum_{n\geq 0} \left(
                        \begin{array}{c} m \\ n \end{array} \right)
                        (x_{n} \cdot P)_{(m+k-n)}.
$$

Using the Reconstruction Theorem 3.6.1 of \cite{vertex}, it is easy to
prove that these formulas define a vertex algebra structure on $M_\g$.

In fact, the same definition works if we replace $N_+$ by any
algebraic group $G$. In the general case, it is natural to consider
the central extension $\ghat_\ka$ of the loop algebra $\g((t))$
corresponding to an invariant inner product $\ka$ on $\g$ defined as
in \secref{hom of VAs}. Then we have the induced module
$$
\on{Ind}_{\g[[t]] \oplus \C K}^{\ghat_\ka} \on{Fun}(G[[t]]),
$$
where the central element $K$ acts on $\on{Fun}(G[[t]])$ as the
identity. The corresponding vertex algebra is the algebra of chiral
differential operators on $G$, considered in \cite{GMS,AG}.  As shown
in \cite{GMS,AG}, in addition to the natural (left) action of
$\ghat_\ka$ on this vertex algebra, there is another (right) action of
$\ghat_{-\ka-\ka_{_K}}$, which commutes with the left action. Here
$\ka_{_K}$ is the Killing form on $\g$, defined by the formula
$\ka_{_K}(x,y) =\on{Tr}_\g (\on{ad} x \on{ad} y)$. In the case when
$\g=\n_+$, there are no non-zero invariant inner products (in
particular, $\ka_{_K} = 0$), and so we obtain a commuting right action
of $\n_+((t))$ on $M_\g$. We will use this right action below (see
\remref{right action}).

\begin{remark}    \label{indep of point}
The above formulas in fact define a canonical vertex algebra structure
on
$$
\on{Ind}_{\n_+[[t]]}^{\n_+((t))} \on{Fun}(U[[t]]),
$$
which is independent of the choice of identification $N_+ \simeq
U$. Recall that in order to define $U$ we only need to fix a Borel
subgroup $B_+$ of $G$. Then $U$ is defined as the open $B_+$--orbit of
the flag manifold and so it is naturally an $N_+$--torsor. In order to
identify $U$ with $N_+$ we need to choose a point in $U$, i.e., an
opposite Borel subgroup $B_-$, or equivalently, a Cartan subgroup $H =
B_+ \cap B_-$ of $B_+$. But in the above formulas we never used an
identification of $N_+$ and $U$, only the action of $L\n_+$ on $LU$,
which is determined by the canonical $\n_+$--action on $\C[U]$.

If we do not fix an identification $N_+ \simeq U$, then the right
action of $\n_+((t))$ discussed above becomes an action of the
``twisted'' Lie algebra $\n_{+,U}((t))$, where $\n_{+,U} = U
\underset{N_+}\times \n_+$. It is interesting to observe that unlike
$\n_+$, this twisted Lie algebra $\n_{+,U}$ has a canonical
decomposition into the one-dimensional subspaces $\n_{\al,U}$
corresponding to the roots $\De_+$. Therefore there exist canonical
(up to a scalar) generators $e^R_i$ of the right action of $\n_+$ on
$U$ (this observation is due to D.~Gaitsgory). We use these operators
below to define the screening operators, and their independence of the
choice of the Cartan subgroup in $B_+$ implies the independence of the
kernel of the screening operator from any additional choices.\qed
\end{remark}

\subsection{Local extension}    \label{local extension}

For our purposes we may replace $\wt{\cA}^\g$, which is a very large
topological algebra, by a relatively small ``local part'' $U(M_\g)$.
Accordingly, we replace $\wt{\cA}^\g_0$ and $\wt{\cA}^\g_{\leq 1}$ by
their local versions $\cA^\g_{0,\on{loc}} = \wt{\cA}^\g_0 \cap
U(M_\g)$ and $\cA^\g_{\leq 1,\on{loc}} = \wt{\cA}^\g_{\leq 1} \cap
U(M_\g)$.

Let us describe $\cA^\g_{0,\on{loc}}$ and $\cA^\g_{\leq 1,\on{loc}}$ more
explicitly. The space $\cA^\g_{0,\on{loc}}$ is spanned (topologically)
by the Fourier coefficients of all polynomials in the $\pa_z^n
a^*_\al(z), n \geq 0$. Note that because the $a^*_{\al,n}$'s commute
among themselves, these polynomials are automatically normally
ordered. The space $\cA^\g_{\leq 1,\on{loc}}$ is spanned by the Fourier
coefficients of the fields of the form $\Wick P(a^*_\al(z),\pa_z
a^*_\al(z),\ldots) a_\beta(z)\Wick $ (the normally ordered product of
$P(a^*_\al(z),\pa_z a^*_\al(z),\ldots)$ and $a_\beta(z)$). Observe
that the Fourier coefficients of all fields of the form
$$
\Wick P(a^*_\al(z),\pa_z a^*_\al(z),\ldots) \pa_z^m a_\beta(z)\Wick \;
, \quad \quad m>0,
$$
may be expressed as linear combinations of the fields of the form
\begin{equation}    \label{xi}
\Wick P(a^*_\al(z),\pa_z a^*_\al(z),\ldots) a_\beta(z)\Wick  \; .
\end{equation}

Further, we define a local version $\T^\g_{\on{loc}}$ of $\on{Vect} LU$
as the subspace, which consists of finite linear combinations of
Fourier coefficients of the formal power series
\begin{equation}    \label{xi1}
P(a^*_\al(z), \pa_z a^*_\al(z),\ldots) a_\beta(z),
\end{equation}
where $a_\al(z)$ and $a^*_\al(z)$ are given by formulas
\eqref{az}, \eqref{a*z}.

Since $\cA^\g_{\leq 1,\on{loc}}$ is the intersection of Lie subalgebras
of $\wt{\cA}^\g$, it is also a Lie subalgebra of $\wt{\cA}$.  By
construction, its image in $\on{Vect} LU$ under the homomorphism
$\wt{\cA}^\g_{\leq 1} \to \on{Vect} LU$ equals
$\T^\g_{\on{loc}}$. Finally, the kernel of the resulting surjective
Lie algebra homomorphism $\cA^\g_{\leq 1,\on{loc}} \to \T^\g_{\on{loc}}$
equals $\cA^\g_{0,\on{loc}}$. Hence we obtain that the extension
\eqref{local sequence} restricts to the ``local'' extension
\begin{equation}    \label{local ext}
0 \to \cA^\g_{0,\on{loc}} \to \cA^\g_{\leq 1,\on{loc}} \to
\T^\g_{\on{loc}} \to 0.
\end{equation}
This sequence is non-split (see Lemma 11.1.9 of \cite{vertex}). The
corresponding two-cocycle will be computed explicitly in
\lemref{WickOPE} using the Wick formula (it comes from the ``double
contractions'' of the corresponding vertex operators).

According to \secref{action of Lg}, the image of $L\g$ in $\on{Vect}
LU$ belongs to $\T^\g_{\on{loc}}$. We will show that the homomorphism
$L\g \to \T^\g_{\on{loc}}$ may be lifted to a homomorphism
$\ghat_{\ka} \to \cA^\g_{\leq 1,\on{loc}}$, where $\ghat_\ka$ is the
central extension of $L\g$ defined in \secref{hom of VAs}.

\subsection{Computation of the two-cocycle of the extension
  \eqref{local ext}}    \label{two cocycle}

Recall that an exact sequence of Lie algebras
$$
0 \to \h \to \wt\g \to \g \to 0,
$$
where $\h$ is an abelian ideal, with prescribed $\g$--module
structure, gives rise to a two-cocycle of $\g$ with coefficients in
$\h$. It is constructed as follows. Choose a splitting $\imath: \g \to
\wt\g$ of this sequence (considered as a vector space), and define
$\sigma: \bigwedge^2 \g \to \h$ by the formula
$$
\sigma(a,b) = \imath([a,b]) - [\imath(a),\imath(b)].
$$
One checks that $\sigma$ is a two-cocycle in the Chevalley complex
of $\g$ with coefficients in $\h$, and that changing the splitting
$\imath$ amounts to changing $\sigma$ by a coboundary.

Conversely, suppose we are given a linear functional $\sigma:
\bigwedge^2 \g \to \h$. Then we associate to it a Lie algebra
structure on the direct sum $\g \oplus \h$. Namely, the Lie bracket of
any two elements of $\h$ is equal to $0$, $[X,h] = X \cdot h$ for all
$X \in \g, h \in \h$, and
$$
[X,Y] = [X,Y]_\g + \omega(X,Y), \qquad X,Y \in \g.
$$
These formulas define a Lie algebra structure on $\wt\g$ if and only if
$\sigma$ is a two-cocycle in the standard Chevalley complex of $\g$
with coefficients in $\h$. Therefore we obtain a bijection between the
set of isomorphism classes of extensions of $\g$ by $\h$ and the
cohomology group $H^2(\g,\h)$.

Consider the extension \eqref{local ext}. The operation of normal
ordering gives us a splitting $\imath$ of this extension as vector
space. Namely, $\imath$ maps the $n$th Fourier coefficient of the
series \eqref{xi1} to the $n$th Fourier coefficient of the series
\eqref{xi}. To compute the corresponding two-cocycle we have to learn
how to compute commutators of Fourier coefficients of generating
functions of the form \eqref{xi} and \eqref{xi1}. Those may be
computed from the operator product expansion (OPE) of the
corresponding vertex operators. We now explain how to compute the OPEs
of vertex operators using the Wick formula following \cite{vertex},
\S\S~11.2.4--11.2.10.

Let us describe the operation of normal ordering more explicitly. In
order to simplify notation, we will assume that $\g=\sw_2$ and
suppress the index $\al$ in $a_{\al,n}$ and $a^*_{\al,n}$. The general
case is treated in the same way. We call the generators $a_n, n\geq
0$, and $a^*_m, m>0$, annihilation operators, and the generators $a_n,
n<0$, and $a^*_m, m\leq 0$, creation operators. \index{annihilation
operator} \index{creation operator} A monomial $P$ in $A$ is called
normally ordered if all factors of $P$ which are annihilation
operators stand to the right of all factors of $P$ which are creation
operators. Given any monomial $P$, we define the normally ordered
monomial $\Wick P\Wick $ as the monomial obtained by moving all
factors of $P$ which are annihilation operators to the right, and all
factors of $P$ which are creation operators to the left. Note that
since the annihilation operators commute among themselves, it does not
matter how we order them among themselves. The same is true for the
creation operators. This shows that $\Wick P\Wick $ is well-defined by
the above conditions.

Given two monomials $P$ and $Q$, their normally ordered product is by
definition the normally ordered monomial $\Wick PQ\Wick $. By
linearity, we define the normally ordered product of any two vertex
operators from the vertex algebra $M_\g$ by applying the above
definition to each of their Fourier coefficients.

\subsection{Contractions of fields}

In order to state the Wick formula, we have to introduce the notion of
contraction of two fields. \index{contraction}

{}From the commutation relations, we obtain the following OPEs:
\begin{align*}
a(z) a^*(w) &= \frac{1}{z-w} + \Wick a(z) a^*(w)\Wick \; , \\
a^*(z) a(w) &= - \frac{1}{z-w} + \Wick a^*(z) a(w)\Wick \; .
\end{align*}
We view them now as identities on formal power series, in which by
$1/(z-w)$ we understand its expansion in positive powers of
$w/z$. Differentiating several times, we obtain
\begin{align}    \label{panpam}
\pa_z^n a(z) \pa^m_w a^*(w) &= (-1)^n \frac{(n+m)!}{(z-w)^{n+m+1}} +
\Wick \pa_z^n a(z) \pa^m_w a^*(w)\Wick \; , \\
\pa^m_z a^*(z) \pa^n_w a(w) &= (-1)^{m+1} \frac{(n+m)!}{(z-w)^{n+m+1}} +
\Wick \pa^m_z a^*(z) \pa^n_w a(w)\Wick 
\end{align}
(here again by $1/(z-w)^n$ we understand its expansion in positive
powers of $w/z$).

Suppose that we are given two normally ordered monomials in $a(z),
a^*(z)$ and their derivatives. Denote them by $P(z)$ and $Q(z)$. A
single {\em pairing} between $P(z)$ and $Q(w)$ is by definition either
the pairing $(\pa^n_z a(z),\pa^m_w a^*(w))$ of $\pa^n_z a(z)$ occurring
in $P(z)$ and $\pa^m_w a^*(w)$ occurring in $Q(w)$, or $(\pa^m_z
a^*(z),\pa^n_w a(w))$ of $\pa^m_z a^*(z)$ occurring in $P(z)$ and
$\pa^n_w a(w)$ occurring in $Q(w)$. We attach to it the functions $$\ds
(-1)^n
\frac{(n+m)!}{(z-w)^{n+m+1}} \quad \quad \on{and} \quad \quad \ds
(-1)^{m+1} \frac{(n+m)!}{(z-w)^{n+m+1}},$$ respectively. A multiple
pairing $B$ is by definition a disjoint union of single pairings. We
attach to it the function $f_B(z,w)$, which is the product of the
functions corresponding to the single pairings in $B$.

Note that the monomials $P(z)$ and $Q(z)$ may well have multiple
pairings of the same type. For example, the monomials $\Wick a^*(z)^2
\pa_z a(z)\Wick $ and $\Wick a(w) \pa_z^2 a^*(w)\Wick $ have two
different pairings of type $(a^*(z),a(w))$; the corresponding function
is $-1/(z-w)$. In such a case we say that the multiplicity of this
pairing is $2$. Note that these two monomials also have a unique
pairing $(\pa_z a(z),\pa_w^2 a^*(w))$, and the corresponding function
is $-6/(z-w)^4$.

Given a multiple pairing $B$ between $P(z)$ and $Q(w)$, we define
$(P(z)Q(w))_B$ as the product of all factors of $P(z)$ and $Q(w)$
which do not belong to the pairing (if there are no factors left, we
set $(P(z)Q(w))_B = 1$). The {\em contraction} \index{contraction} of
$P(z)Q(w)$ with respect to the pairing $B$, denoted $\Wick
P(z)Q(w)\Wick _B$, is by definition the normally ordered formal power
series $\Wick (P(z)Q(w)_B)\Wick$ multiplied by the function
$f_B(z,w)$. We extend this definition to the case when $B$ is the
empty set by stipulating that
$$
\Wick P(z)Q(w)\Wick_{\emptyset} = \Wick P(z)Q(w)\Wick \; .
$$

Now we are in a position to state the {\em Wick formula}, which gives
the OPE of two arbitrary normally ordered monomial vertex
operators. The proof of this formula is straightforward and is left to
the reader.

\begin{lem}    \label{Wick form}
Let $P(z)$ and $Q(w)$ be two monomials as above. Then $P(z)Q(w)$
equals the sum of terms $\Wick P(z)Q(w)\Wick_B$ over all pairings $B$
between $P$ and $Q$ including the empty one, counted with
multiplicity.
\end{lem}

Now we can compute our two--cocycle. For that, we need to apply the
Wick formula to the fields of the form
$$\Wick R(a^*(z),\pa_z a^*(z),\ldots) a(z)\Wick \; ,
$$
whose Fourier coefficients span the pre-image of $\T_{\on{loc}}$ in
$A_{\leq 1,\on{loc}}$ under our splitting $\imath$. Two fields of this
form may have only single or double pairings, and therefore their OPE
can be written quite explicitly.

A field of the above form may be written as $Y(P(a^*_n) a_{-1},z)$ (or
$Y(P a_{-1},z)$ for short), where $P$ is a polynomial in the $a^*_n, n
\leq 0$ (recall that $a^*_n, n\leq 0$, corresponds to $\pa_z^{-n}
a^*(z)/(-n)!$). Applying the Wick formula, we obtain

\begin{lem}    \label{WickOPE}
\begin{align*}
Y(P a_{-1},z) Y(Q a_{-1},w)  &=
\Wick Y(Pa_{-1},z) Y(Qa_{-1},w)\Wick \; \\ &+ \sum_{n\geq 0}
\frac{1}{(z-w)^{n+1}}
\Wick Y(P,z) Y\left( \frac{\pa Q}{\pa a^*_{-n}} a_{-1},w \right)\Wick
\\ &- \sum_{n\geq 0} \frac{1}{(z-w)^{n+1}} \Wick Y\left( \frac{\pa
  P}{\pa a^*_{-n}} a_{-1},z \right) Y(Q,w)\Wick  \\ &- \sum_{n,m\geq 0}
\frac{1}{(z-w)^{n+m+2}} Y\left( \frac{\pa P}{\pa a^*_{-n}},z
\right) Y\left( \frac{\pa Q}{\pa a^*_{-m}},w \right)
\end{align*}
\end{lem}

Note that we do not need to put normal ordering in the last summand.

\subsection{Double contractions}    \label{double contractions}

Using this formula, we can now easily obtain the commutators of the
Fourier coefficients of the fields $Y(P a_{-1},z)$ and $Y(Q
a_{-1},w)$, using the residue calculus (see \cite{vertex}, \S~3.3).

The first two terms in the right hand side of the formula in
\lemref{WickOPE} correspond to single contractions between $Y(P
a_{-1},z)$ and $Y(Q a_{-1},w)$. The part in the commutator of the
Fourier coefficients induced by these terms will be exactly the same
as the commutator of the corresponding vector fields, computed in
$\T_{\on{loc}}$. Thus, we see that the discrepancy between the
commutators in $A_{\leq 1,\on{loc}}$ and in $\T_{\on{loc}}$ (as
measured by our two--cocycle) is due to the last term in the formula
from \lemref{WickOPE}, which comes from the {\em double contractions}
between $Y(P a_{-1},z)$ and $Y(Q a_{-1},w)$.

Explicitly, we obtain the following formula for our two-cocycle
\begin{equation}    \label{omega}
\omega((Pa_{-1})_{(k)},(Qa_{-1})_{(s)})
\end{equation}
$$
= - \sum_{n,m\geq 0} \int
\left. \frac{1}{(n+m+1)!} \pa_z^{n+m+1} Y\left( \frac{\pa P}{\pa
  a^*_{-n}},z \right) Y\left( \frac{\pa Q}{\pa a^*_{-m}},w \right) z^k
w^s \right|_{z=w} \; dw.
$$

\subsection{A reminder on cohomology}

Let again
\begin{equation*}
        0 \rightarrow \h \rightarrow \lt \rightarrow
        \li \rightarrow 0
\end{equation*}
be an extension of Lie algebras, where $\h$ is an abelian Lie
subalgebra and an ideal in $\lt$. Choosing a splitting $\imath$ of
this sequence considered as a vector space we define a two-cocycle of
$\li$ with coefficients in $\h$ as in \secref{two cocycle}. Suppose
that we are given a Lie algebra homomorphism $\al: \g \rightarrow \li$
for some Lie algebra $\g$. Pulling back our two-cocycle under $\al$ we
obtain a two-cocycle of $\g$ with coefficients in $\h$.

On the other hand, given a homomorphism $\h' \xrightarrow{i} \h$ of
$\g$--modules we obtain a map $i_*$ between the spaces of two-cocycles
of $\g$ with coefficients in $\h$ and $\h'$. The corresponding map of
the cohomology groups $H^2(\g,\h') \rightarrow H^2(\g,\h)$ will also
be denoted by $i_*$.

\begin{lem} \label{L:diag}
Suppose that we are given a two-cocycle $\sigma$ of $\g$ with
coefficients in $\h'$ such that the cohomology classes of
$i_{*}(\sigma)$ and $\omega$ are equal in $H^2(\g,\h')$. Denote by
$\widetilde{\g}$ the extension of $\g$ by $\h'$ corresponding to
$\sigma$ defined as above. Then the map $\g \rightarrow \li$ may be
augmented to a map of commutative diagrams
\begin{equation}    \label{nu vot}
\begin{CD}
0 @>>> \h @>>> \lt @>>> \li @>>> 0 \\
& & @AAA @AAA @AAA \\
0 @>>> \h' @>>> \wt{\g} @>>> \g @>>> 0
\end{CD}
\end{equation}
Moreover, the set of isomorphism classes of such diagrams is a
torsor over $H^1(\g,\h)$.
\end{lem}

\begin{proof}
If the cohomology classes of $i_{*}(\sigma)$ and $\omega$ coincide,
then $i_*(\sigma) + d\ga = \omega$, where $\ga$ is a one-cochain,
i.e., a linear functional $\g \to \h$. Define a linear map $\beta:
\wt{\g} \to \lt$ as follows. By definition, we have a splitting
$\wt{\g} = \g \oplus \h'$ as a vector space. We set $\beta(X) =
\imath(\al(X)) + \ga(X)$ for all $X \in \g$ and $\beta(h) = i(h)$ for
all $h \in \h$. Then the above equality of cocycles implies that
$\beta$ is a Lie algebra homomorphism which makes the diagram
\eqref{nu vot} commutative. However, the choice of $\ga$ is not unique
as we may modify it by adding to it an arbitrary one-cocycle
$\gamma'$. But the homomorphisms corresponding to $\ga$ and to $\ga +
\ga'$, where $\ga'$ is a coboundary, lead to isomorphic diagrams. This
implies that the set of isomorphism classes of such diagrams is a
torsor over $H^1(\g,\h)$.
\end{proof}

\subsection{Two cocycles}

Restricting the two-cocycle $\omega$ of $\T_{\on{loc}}$ with
coefficients in $\cA^\g_{0,\on{loc}}$ corresponding to the extension
\eqref{local ext} to $L\g \subset \T_{\on{loc}}$ we obtain a
two-cocycle of $L\g$ with coefficients in $\cA^\g_{0,\on{loc}}$. We
also denote it by $\omega$. The $L\g$--module $\cA^\g_{0,\on{loc}}$
contains the trivial subrepresentation ${\mathbb C}$, i.e., the span
of $\int Y(\vac,z) dz/z$ (which we view as the constant function on
$LU$), and the inclusion ${\mathbb C} \xrightarrow{i}
\cA^\g_{0,\on{loc}}$ induces a map $i_*$ of the corresponding spaces
of two-cocycles and the cohomology groups $H^2(L\g,\C) \rightarrow
H^2(L\g,\cA^\g_{0,\on{loc}})$.

It is well known that $H^2(L\g,\C)$ is one-dimensional and is
isomorphic to the space of invariant inner products on $\g$ (the
corresponding Kac-Moody two-cocyles are described in \secref{hom of
VAs}). We denote by $\sigma$ the class corresponding to the inner
product $\ka_c = - \frac{1}{2} \ka_{_K}$, where $\ka_{_K}$ denotes the
Killing form on $\g$. Thus, by definition, $$\ka_c(x,y) = -
\frac{1}{2} \on{Tr} (\on{ad} x \on{ad} y).$$ We will show that the
cohomology classes of $i_*(\sigma)$ and $\omega$ are
equal. \lemref{L:diag} will then imply that there exists a family of
Lie algebra homomorphisms $\ghat_{\ka_c} \rightarrow \cA^\g_{\leq
1,\on{loc}}$ such that $K \mapsto 1$.

Unfortunately, the Chevalley complex that calculates
$H^2(L\g,\cA^\g_{0,\on{loc}})$ is unmanageably large. So as the first
step we will show in the next section that $\omega$ and $\sigma$
actually both belong to a much smaller subcomplex, where they are more
easily compared.

\section{Comparison of cohomology classes}

\subsection{Clifford Algebras}

Choose a basis $\{ J^a \}_{a=1,\ldots,\dim \g}$ of $\g$, and set
$J^a_{n}=J^a \otimes t^{n}$. Introduce the Clifford algebra with
generators $\psi_{a,n}, \psi^{*}_{a,m}, a=1,\ldots,\dim \g; m,n \in
\Z$, with anti-commutation relations
$$
[\psi_{a,n},\psi_{b,n}]_{+}=[\psi^{*}_{a,n},\psi^{*}_{b,m}]_{+}=0,
\qquad [\psi_{a,n}, \psi^*_{b,m}]_{+}=\delta_{a,b} \delta_{n,-m}.
$$

Let $\bigwedge_{\g}$ be the module over this Clifford algebra
generated by a vector $\vac$ such that
$$
\psi_{a,n} \vac = 0, \quad n>0, \qquad
\psi^*_{a,n} \vac = 0, \quad n\geq 0.
$$
Then $\bigwedge_{\g}$ carries the following structure of a vertex
superalgebra (see \cite{vertex}, \S~14.1.1):

\begin{enumerate}

\item[$\bullet$] Gradation: $\on{deg} \psi_{a,n} = \deg \psi^*_{a,n}
= -n, \deg \vac = 0$;

\item[$\bullet$] Vacuum vector: $\vac$;

\item[$\bullet$] Translation operator: $T\vac=0, [T,\psi_{a,n}] =
-n\psi_{a,n-1}, [T,\psi^*_{a,n}] =$ \newline $-(n-1)
\psi^*_{a,n-1}$.

\item[$\bullet$] Vertex operators:
\begin{align*} 
        Y(\psi_{a,-1} \vac,z) &= \psi_a(z) = \sum_{n \in \Z}
        \psi_{a,n} z^{-n-1}, \\
        Y(\psi^*_{a,0} \vac,z) &= \psi^*_a(z) = \sum_{n \in \Z}
        \psi^*_{a,n} z^{-n},
\end{align*}
\begin{multline*}
Y(\psi_{a_1,n_1} \ldots \psi_{a_k,n_k} \psi^*_{b_1,m_1} \ldots
\psi^*_{b_l,m_l} \vac,z) = \prod_{i=1}^k \frac{1}{(-n_i-1)!}
\prod_{j=1}^l \frac{1}{(-m_j)!} \cdot \\ \Wick
\delz^{-n_1-1}\psi_{a_1}(z)\dots\delz^{-n_k-1}\psi_{a_k}(z)
\delz^{-m_1}\psi^*_{b_1}(z)\dots\delz^{-m_l}\psi^*_{b_l}(z) \Wick \; .
\end{multline*}

\end{enumerate}
The tensor product of two vertex superalgebras is naturally a vertex
superalgebra (see Lemma 1.3.6 of \cite{vertex}), and so $M_\g \otimes
\bigwedge_\g$ is a vertex superalgebra.

\subsection{The local Chevalley complex}

The ordinary Chevalley complex computing the Lie algebra cohomology
$H^\bullet(L\g,\cA^\g_{0,\on{loc}})$ is
$C^\bullet(L\g,\cA^\g_{0,\on{loc}}) = \bigoplus_{i\geq 0}
C^i(L\g,\cA^\g_{0,\on{loc}})$,
$$C^i(L\g,\cA^\g_{0,\on{loc}}) =
\on{Hom}_{\on{cont}}(\bigwedge{}^i(L\g),\cA^\g_{0,\on{loc}}),$$ where
$\bigwedge^i(L\g)$ stands for the natural completion of the ordinary
$i$th exterior power of the topological vector space $L\g$, and we
consider all continuous linear maps. For $f \in \cA^\g_{0,\on{loc}}$ we
denote the linear functional $\phi \in
\on{Hom}_{\on{cont}}(\bigwedge^i(L\g),\cA^\g_{0,\on{loc}})$ defined by
the formula
$$
\phi(J^{a_1}_{m_1} \wedge \ldots \wedge J^{a_i}_{m_i}) = \begin{cases}
f, & \on{if} \quad ((a_1,m_1),\ldots,(a_i,m_i)) =
((b_1,k_1),\ldots,(b_i,k_i)),\\ 0, & \on{if} \quad
((a_1,m_1),\ldots,(a_i,m_i)) \neq \tau((b_1,k_1),\ldots,(b_i,k_i)),
\end{cases}
$$
where $\tau$ runs over the symmetric group on $i$ letters, by
$\psi^*_{b_1,-k_1} \ldots \psi^*_{b_i,-k_i} f$. Then any element of
the space $C^i(L\g,\cA^\g_{0,\on{loc}})$ may be written as a (possibly
infinite) linear combination of terms of this form.

The differential $d:
C^i(L\g,\cA^\g_{0,\on{loc}}) \to C^{i+1}(L\g,\cA^\g_{0,\on{loc}})$ is
given by the formula
\begin{align*} 
        (d \phi)(X_{1},\ldots,X_{i+1}) &= \sum_{j=1}^{i+1}
(-1)^{j+1} X_{j} \phi(X_{1},\ldots,\wh{X_{j}},\ldots,X_{i+1}) \\ &+
\sum_{j<k} (-1)^{j+k+1} \phi([X_{j},X_{k}],X_{1},\ldots,
\wh{X_{j}},\ldots,\wh{X_{k}},\ldots,X_{i+1}).
\end{align*}

Let $C^i_{\on{loc}}(L\g,\cA^\g_{0,\on{loc}})$ be the subspace of
$\on{Hom}_{\on{cont}}(\bigwedge^i(L\g),\cA^\g_{0,\on{loc}})$, spanned
by all linear maps of the form
\begin{equation}    \label{A form}
        \int Y(\psi^{*}_{a_{1},n_{1}} \cdots \psi^{*}_{a_{i},n_{i}}
\astar_{\alpha_{1},m_{1}} \cdots \astar_{\alpha_{j},m_{j}} \vac, z) \;
dz, \qquad n_{p} \leq 0, m_{p} \leq 0.
\end{equation}

\begin{lem}    \label{preserved}
$C^\bullet_{\on{loc}}(L\g,\cA^\g_{0,\on{loc}})$ is preserved by the
differential $d$, and so it forms a subcomplex of
$C^\bullet(L\g,\cA^\g_{0,\on{loc}})$.
\end{lem}

\begin{proof}
The action of the differential $d$ on $\int Y(A,z) dz$, where $A$ is
of the form \eqref{A form}, may be written as the commutator
$$\left[\int Q(z) dz,\int Y(A,z) dz \right],$$ where $Q(z)$ is the
field
\begin{equation}    \label{Qz}
Q(z) = Y(Q,z) = \sum_{a} J^a(z) \psi^*_a(z) - \frac{1}{2} \sum_{a,b,c}
\mu^{ab}_c \Wick \psi^*_a(z) \psi^*_b(z) \psi_c(z)\Wick \, .
\end{equation}
Therefore it is equal to $Y\left(\int Q(z) dz \cdot A,z \right)$,
which is of the form \eqref{A form}.
\end{proof}

We call $C^i_{\on{loc}}(L\g,\cA^\g_{0,\on{loc}})$ the {\em local}
subcomplex of $C^\bullet(L\g,\cA^\g_{0,\on{loc}})$.

\begin{lem}    \label{repr in loc}
Both $\omega$ and $i_{*}(\sigma)$ belong to
$C^2_{\on{loc}}(L\g,\cA^\g_{0,\on{loc}})$.
\end{lem}

\begin{proof}
We begin with $\sigma$ given by the formula
\begin{equation}    \label{sigma}
        \sigma(J^a_{n},J^b_{m}) = n \delta_{n,-m} \ka_c(J^a,J^b)
\end{equation}
(see \secref{hom of VAs}). Therefore the cocycle $i_*(\sigma)$ is
equal to the following element of
$C^2_{\on{loc}}(L\g,\cA^\g_{0,\on{loc}})$,
\begin{align} \notag
i_*(\sigma) & = \sum_{a\leq b} \sum_{n \in \Z} \ka_c (J^a,J^b)
        n \psi^{*}_{a,-n} \psi^{*}_{b,n} \\ & = \sum_{a \leq b}
        \ka_c(J^a,J^b) \int Y(\psi^{*}_{a,-1} \psi^{*}_{b,0} \vac,z)
        \; dz.  \label{istar}
\end{align}

Next we consider $\omega$. Combining the discussion of \secref{action
of Lg} with formulas of \secref{sec on formulas}, we obtain that
$$
\imath(J^a_n) = \sum_{\beta \in \De_+} \int
Y(R_{a}^\beta(\astar_{\alpha,0})a_{\beta,-1} \vac,z) \; z^{n} dz,
$$
where $R_a^\beta$ is a polynomial. Then formula \eqref{omega} implies
that
\begin{multline*}
\omega(J^a_n,J^b_m) = \\ - \sum_{\al,\beta \in \De_+} \int \left(
  Y\left( T \frac{\pa R_a^\al}{\pa a^*_{\beta,0}} \cdot \frac{\pa
  R^\beta_b}{\pa a^*_{\al,0}},w \right) w^{n+m} + n Y\left( \frac{\pa
  R^\al_a}{\pa a^*_{\beta,0}} \frac{\pa R^\beta_b}{\pa a^*_{\al,0}},w
  \right) w^{n+m-1} \right) \; dw.
\end{multline*}
Therefore
\begin{multline}  \label{omega1}
\omega = \\ - \sum_{a\leq b; \al,\beta \in \De_+} \int \left(
        Y\left(\psi^{*}_{a,0} \psi^{*}_{b,0} T \frac{\pa R^\al_a}{\pa
        a^*_{\beta,0}} \cdot \frac{\pa R^\beta_b}{\pa a^*_{\al,0}},z
        \right) + Y\left(\psi^{*}_{a,-1} \psi^{*}_{b,0} \frac{\pa
        R^\al_a}{\pa a^*_{\beta,0}} \frac{\pa R^\beta_b}{\pa
        a^*_{\al,0}},z \right) \right) dz
\end{multline}
(in the last two formulas we have omitted $\vac$). Hence it belongs to
$C^2_{\on{loc}}(L\g,\cA^\g_{0,\on{loc}})$.
\end{proof}

We need to show that the cocycles $i_*(\sigma)$ and $\omega$ represent
the same cohomology class in the local complex
$C^\bullet_{\on{loc}}(L\g,\cA^\g_{0,\on{loc}})$.

\subsection{Another complex}    \label{M+}

Given a Lie algebra ${\mathfrak l}$, we denote the Lie subalgebra
${\mathfrak l}[[t]]$ of $L{\mathfrak l} = {\mathfrak l}((t))$ by
$L_+{\mathfrak l}$.

Observe that the Lie algebra $L_+\n_+$ acts naturally on the space
$$M_{\g,+} = \C[a^*_{\al,n}]_{\al \in \De_+,n\leq 0} \simeq
\C[y_{\al,n}]_{\al \in \De_+,n\geq 0},$$
which is identified with the
ring of functions on the space $U[[t]] \simeq N_+[[t]]$. We identify
the standard Chevalley complex $C^\bullet(L_+\g,M_{\g,+}) =
\on{Hom}_{\on{cont}}(\bigwedge^\bullet L_+\g,M_{\g,+})$ with the
tensor product $M_{\g,+} \otimes \bigwedge^\bullet_{\g,+}$, where
$$\bigwedge{}^\bullet_{\g,+} = \bigwedge(\psi^*_{a,n})_{n\leq 0}.$$
Introduce a superderivation $T$ on $C^\bullet(L_+\g,M_{\g,+})$ acting
by the formulas $$ T \cdot a^*_{a,n} = -(n-1) a^*_{a,n-1}, \qquad T
\cdot \psi^*_{a,n} = -(n-1) \psi^*_{a,n-1}.  $$ We have a linear map
$C^\bullet(L_+\g,M_{\g,+}) \to
C^\bullet_{\on{loc}}(L\g,\cA^\g_{0,\on{loc}})$ sending $A \in
C^\bullet(L_+\g,M_{\g,+})$ to $\int Y(A,z) dz$. It is clear that if $A
\in \on{Im} T$, then $\int Y(A,z) dz = 0$.

\begin{lem}
The map $\int$ defines an isomorphism
$$C^\bullet_{\on{loc}}(L\g,\cA^\g_{0,\on{loc}}) \simeq
C^\bullet(L_+\g,M_{\g,+})/(\on{Im} T + \C),$$ and $\on{Ker} T =
\C$. Moreover, the following diagram is commutative
$$
\begin{CD}
C^\bullet(L_+\g,M_{\g,+}) @>{T}>> C^\bullet(L_+\g,M_{\g,+}) @>{\int}>>
C^\bullet_{\on{loc}}(L\g,\cA^\g_{0,\on{loc}}) \\ @A{d}AA @A{d}AA
@A{d}AA \\ C^\bullet(L_+\g,M_{\g,+}) @>{T}>>
C^\bullet(L_+\g,M_{\g,+}) @>{T}>>
C^\bullet_{\on{loc}}(L\g,\cA^\g_{0,\on{loc}})
\end{CD}
$$
\end{lem}

\begin{proof}
It is easy to see that the differential of the Chevalley complex
$C^\bullet(L_+\g,M_{\g,+})$ acts by the formula $A \mapsto \int Q(z)
dz \cdot A$, where $Q(z)$ is given by formula \eqref{Qz}. The lemma
now follows from the argument used in the proof of \lemref{preserved}.
\end{proof}

Consider the double complex
$$
\begin{CD}
\C @>>> C^\bullet(L_+\g,M_{\g,+}) @>{T}>> C^\bullet(L_+\g,M_{\g,+})
@>>> \C \\
& & @AAA @AAA  \\
\C @>>> C^\bullet(L_+\g,M_{\g,+}) @>{T}>> C^\bullet(L_+\g,M_{\g,+})
@>>> \C
\end{CD}
$$
According to the lemma, the cohomology of the complex
$C^\bullet_{\on{loc}}(L\g,\cA^\g_{0,\on{loc}})$ is given by the second
term of the spectral sequence, in which the zeroth differential is
vertical.

We start by computing the first term of this spectral
sequence. Observe that the module $M_{\g,+}$ is isomorphic to the
coinduced module $\on{Coind}_{L_+\bb_-}^{L_+\g} \C$.  Define a map of
complexes $$\mu': C^\bullet(L_+\g,M_{\g,+}) \to
C^\bullet(L_+\bb_-,\C)$$ as follows. If $\ga$ is an $i$--cochain in the
complex $C^\bullet(L_+\g,M_{\g,+}) = \on{Hom}_{\on{cont}}(\bigwedge^i
L_+\g,M_{\g,+})$, then $\mu'(\ga)$ is by definition the restriction of
$\ga$ to $\bigwedge^i L_+\bb_-$ composed with the natural projection
$M_{\g,+} = \on{Coind}_{L_+\bb_-}^{L_+\g} \C \to \C$. It is clear that
$\mu'$ is a morphism of complexes. The following is an example of the
Shapiro lemma (see \cite{Fuchs}, \S~5.4, for the proof).

\begin{lem}    \label{shapiro lemma}
The map $\mu'$ induces an isomorphism at the level of
cohomologies, i.e.,
\begin{equation}    \label{shapiro}
H^\bullet(L_+\g,M_{\g,+}) \simeq H^\bullet(L_+\bb_-,\C).
\end{equation}
\end{lem}

Now we compute the right hand side of \eqref{shapiro}. Since $L_+\n_-
\in L_+\bb_-$ is an ideal, and $L_+\bb_-/L_+\n_- \simeq L_+\h$, we
obtain from the Serre-Hochshield spectral sequence (see \cite{Fuchs},
\S~5.1) that
\begin{equation*}
H^n(L_+\bb_-,\C) = \bigoplus_{p+q=n} H^p(L_+\h,H^q(L_+\n_-,\C)).
\end{equation*} 
But $\mathfrak{h}\otimes 1 \in L_+\h$ acts diagonally on
$H^\bullet(L_+\n_-,\C)$, inducing an inner grading. According to
\cite{Fuchs}, \S~5.2, $$H^p(L_+\h,H^q(L_+\n_-,\C)) =
H^p(L_+\h,H^q(L_+\n_-,\C)_0),$$ where $H^q(L_+\n_-,\C)_0$ is the
subspace where $\mathfrak{h}\otimes 1$ acts by $0$. Clearly,
$H^0(L_+\n_-,\C)_0$ is the one-dimensional subspace of the scalars $\C
\subset H^0(L_+\n_-,\C)$, and $H^q(L_+\n_-,\C)_0 = 0$ for $q\neq
0$. Thus, we find that
\begin{equation*} 
H^p(L_+\h,H^q(L_+\n_-,\C)) = H^p(L_+\h,\C). \\
\end{equation*}
Furthermore, we have the following result. Define a map of complexes
\begin{equation}    \label{def of mu}
\mu: C^\bullet(L_+\g,M_{\g,+}) \to C^\bullet(L_+\h,\C)
\end{equation}
as follows. If $\ga$ is an $i$--cochain in the complex
$C^\bullet(L_+\g,M_{\g,+}) = \on{Hom}_{\on{cont}}(\bigwedge^i
L_+\g,M_{\g,+})$, then $\mu(\ga)$ is by definition the restriction of
$\ga$ to $\bigwedge^i L_+\h$ composed with the natural projection
\begin{equation}    \label{projection}
p: M_{\g,+} = \on{Coind}_{L_+\bb_-}^{L_+\g} \C \to \C.
\end{equation}

\begin{lem}    \label{map mu}
The map $\mu$ induces an isomorphism at the level of cohomologies,
i.e.,
$$
H^\bullet(L_+,\g,M_{\g,+}) \simeq H^\bullet(L_+\h,\C).
$$
In particular, the cohomology class of a cocycle in the
complex $C^i(L_+\g,M_{\g,+})$ is uniquely determined by its
restriction to $\bigwedge^i L_+\h$.
\end{lem}

\begin{proof}
It follows from the construction of the Serre-Hochshield spectral
sequence (see \cite{Fuchs}, \S~5.1) that the restriction map
$C^\bullet(L_+\bb_-,\C) \to C^\bullet(L_+\h,\C)$ induces an
isomorphism at the level of cohomologies. The statement of the lemma
now follows by combining this with \lemref{shapiro lemma}.
\end{proof}

Since $L_+\h$ is abelian, we have
\begin{equation*}
H^\bullet(L_+\h,\C) = \bigwedge{}^\bullet(L_+\h),
\end{equation*}
and so
\begin{equation*}
H^\bullet(L_+\g,M_{\g,+}) \simeq \bigwedge{}^\bullet(L_+\h).
\end{equation*}

It is clear that this isomorphism is compatible with the action of $T$
on both sides. The kernel of $T$ acting on the right hand side is
equal to the subspace of the scalars. Therefore we obtain
$$
H^\bullet_{\on{loc}}(L\g,\cA^\g_{0,\on{loc}}) \simeq
\bigwedge{}^\bullet(L_+\h)/(\on{Im} T + \C).
$$

\subsection{Restricting the cocycles}

Any cocycle in $C^i_{\on{loc}}(L\g,\cA^\g_{0,\on{loc}})$ is a cocycle in
$$C^i(L\g,\cA^\g_{0,\on{loc}}) =
\on{Hom}_{\on{cont}}(\bigwedge{}^i(L\g),\cA^\g_{0,\on{loc}}),$$
and as such it may be restricted to $\bigwedge{}^i(L\h)$.

\begin{lem}    \label{any two}
Any two cocycles in $C^i_{\on{loc}}(L\g,\cA^\g_{0,\on{loc}})$ whose
restrictions to $\bigwedge{}^i(L\h)$ coincide represent the same
cohomology class.
\end{lem}

\begin{proof}
We need to show that any cocycle $\phi$ in
$C^i_{\on{loc}}(L\g,\cA^\g_{0,\on{loc}})$, whose restriction to
$\bigwedge{}^i(L\h)$ is equal to zero, is equivalent to the zero
cocycle. According to the above computation, $\phi$ may be written as
$\int Y(A,z) dz$, where $A$ is a cocycle in $C^i(L_+\g,M_{\g,+}) =
\on{Hom}_{\on{cont}}(\bigwedge^i L_+\g,M_{\g,+})$. But then the
restriction of $A$ to $\bigwedge^i L_+\h \subset \bigwedge^i L_+\g$,
denoted by $\overline{A}$, must be in the image of the operator $T$,
and hence so is $p(\ol{A}) = \mu(A) \in \bigwedge^i(L_+\h)$ (here $p$
is the projection defined in \eqref{projection} and $\mu$ is the map
defined in \eqref{def of mu}). Thus, $\mu(A) = T(h)$ for some $h \in
\bigwedge^i(L_+\h)$. Since $T$ commutes with the differential and the
kernel of $T$ on $\bigwedge^i(L_+\h)$ consists of the scalars, we
obtain that $h$ is also a cocycle in $C^i(L_+\h,\C) =
\bigwedge^i(L_+\h)$.

According to \lemref{map mu}, the map $\mu$ induces an isomorphism on
the cohomologies. Hence $h$ is equal to $\mu(B)$ for some cocycle $B$
in $C^i(L_+\g,M_{\g,+})$. It is clear from the definition that $\mu$
commutes with the action of the translation operator $T$. Therefore it
follows that the cocycles $A$ and $T(B)$ are equivalent in
$C^i(L_+\g,M_{\g,+})$. But then $\phi$ is equivalent to the zero
cocycle.
\end{proof}

Now we obtain

\begin{thm}    \label{two eq}
The cocycles $\omega$ and $i_{*}(\sigma)$ represent the same
cohomology class in
$H^2_{\on{loc}}(L\g,\cA^\g_{0,\on{loc}})$. Therefore there exists a
lifting of the homomorphism $L\g \to \T^\g_{\on{loc}}$ to a
homomorphism $\ghat_{\ka_c} \to \cA^\g_{\leq 1,\on{loc}}$ such that $K
\mapsto 1$. Moreover, this homomorphism may be chosen in such a way
that
\begin{equation}    \label{where to map}
J^a(z) \mapsto Y(P_{a}(\astar_{\alpha,0},a_{\beta,-1}) \vac,z) +
Y(B_a,z),
\end{equation}
where $P_a$ is a polynomial introduced in formula \eqref{Pa} and $B_a$
is a polynomial in the $a^*_{\al,n}$'s of degree one.
\end{thm}

\begin{proof}
By \lemref{any two}, it suffices to check that the restrictions of
$\omega$ and $i_{*}(\sigma)$ to $\bigwedge^2(L\h)$ coincide. We have
$$
i_{*}(\sigma)(h_n,h'_m) = n\ka_c(h,h') \delta_{n,-m}
$$
for all $h,h' \in \h$. Now let us compute the restriction of
$\omega$. We find from the formulas given in \secref{sec on formulas}
and \secref{action of Lg} that
$$
\imath(h(z)) = - \sum_{\al \in \Delta_+} \al(h) \; \Wick a^*_\al(z)
a_\al(z) \Wick
$$
(recall that $h(z) = \sum_{n\in\Z} h_n z^{-n-1}$). Therefore we find
that
$$
\omega(h_n,h'_m) = - n \delta_{n,-m} \sum_{\al \in \De_+} \al(h)
\al(h') = n \kappa_c(h,h') \delta_{n,-m},
$$
because by definition $\ka_c(\cdot,\cdot) = -\frac{1}{2}
\ka_{_K}(\cdot,\cdot)$ and for the Killing form $\kappa_{_K}$ we
have
$$
\ka_{_K}(h,h') = 2 \sum_{\al \in \De_+} \al(h) \al(h').
$$
Therefore the cocycles $\omega$ and $i_{*}(\sigma)$ represent the same
class in $H^2_{\on{loc}}(L\g,\cA^\g_{0,\on{loc}})$.

Hence there exists $\gamma \in
C^1_{\on{loc}}(L\g,\cA^\g_{0,\on{loc}})$ such that $\omega + d\gamma
= i_*(\sigma)$. We may write $\gamma$ as
$$
\gamma = \sum_{a} \int \psi^*_a(z) Y(B_a,z) \; dz, \qquad B_a \in
M_{\g,+}.
$$
It follows from the computations made in the proof of \lemref{repr in
loc} that both $\omega$ and $i_*(\sigma)$ are homogeneous elements of
$C^2_{\on{loc}}(L\g,\cA^\g_{0,\on{loc}})$ of degree $0$. Therefore
$\gamma$ may be chosen to be of degree $0$, i.e., $B_a$ may be chosen
to be of degree $1$.

Then by \lemref{L:diag} formulas \eqref{where to map} define a
homomorphism of Lie algebras $\ghat_{\ka_c} \to \cA^\g_{\leq
1,\on{loc}}$ such that $K \mapsto 1$. This completes the proof.
\end{proof}

In the next section we will interpret the above homomorphism in the
vertex algebra language.

\section{Wakimoto modules of critical level}

\subsection{Homomorphism of vertex algebras}    \label{hom of VA}

Recall the definition of the vertex algebra associated to the affine
Kac-Moody algebra $\ghat_\ka$ (see \cite{vertex}, \S~2.4). Define the
vacuum $\G_\ka$--module $V_{\kappa}(\g)$ corresponding to $\kappa$ as
the induced representation
$$
V_{\kappa}(\g) \overset{\on{def}}{=} \on{Ind}_{L_+\g
\oplus \C K}^{\ghat_\ka} \C_1,
$$
where $\C_1$ is the one-dimensional representation on which $L_+\g$
acts by $0$ and $K$ acts as the identity. We denote the vector $1
\otimes 1 \in V_{\kappa}(\g)$ by
$v_{\kappa}$.

Recall that $\kappa$ is unique up to a scalar. Therefore it is often
convenient to fix a particular invariant inner product $\kappa_0$ and
write an arbitrary one as $\kappa = k \kappa_0, k \in \C$. This is the
point of view taken in \cite{vertex}, where as $\kappa_0$ we take the
inner product with respect to which the square of the maximal root is
equal to $2$, and denote the corresponding vacuum module by $V_k(\g)$.
In particular, since the Killing form $\ka_{_K}(\cdot,\cdot)$ is equal
to $2h^\vee$ times $\kappa_0$ (see \cite{Kac}), where $\kk$ is the
dual Coxeter number of $\g$, we obtain that for the inner product
$\ka_c(\cdot,\cdot) = -\frac{1}{2} \ka_{_K}(\cdot,\cdot)$ introduced
above we have $k=-\kk$.  Thus, $V\cri(\g)$ is $V_{-\kk}(\g)$ in the
notation of \cite{vertex}.  We will call this module the vacuum module
of {\em critical level}.

The vertex algebra structure on $V_{\kappa}(\g)$ is
defined as follows:

\begin{enumerate}\label{Vacuum VA}

\item[$\bullet$] Gradation: $\on{deg}J^{a_1}_{n_1}\dots
J^{a_m}_{n_m}v_{\kappa}=-\sum_{i=1}^m {n_i}$.

\item[$\bullet$] Vacuum vector: $\vac=v_{\kappa}$.

\item[$\bullet$] Translation operator: $Tv_{\kappa}=0$,
$[T,J^a_n]=-nJ^a_{n-1}$.

\item[$\bullet$] Vertex operators:
$$Y(J^a_{-1}v_{\kappa},z)=J^a(z)=\sum_{n\in\Z}{J^a_n
z^{-n-1}},$$
\begin{equation*}
Y(J^{a_1}_{n_1}\dots J^{a_m}_{n_m}v_{\kappa},z) = \prod_{i=1}^m
\frac{1}{(-n_i-1)!}  \Wick
\delz^{-n_1-1}J^{a_1}(z)\dots\delz^{-n_m-1}J^{a_m}(z) \Wick
\end{equation*}
\end{enumerate}
(see Theorem 2.4.4 of \cite{vertex}).

\begin{lem}    \label{hom from V}
Defining a homomorphism of vertex algebras $V_{\kappa}(\g) \arr V$ is
equivalent to choosing vertex operators $\wt{J}^a(z), a=1,\ldots,\dim
\g$, of conformal dimension $1$ in the vertex algebra $V$, whose
Fourier coefficients satisfy the relations \eqref{KM rel} with $K=1$.
\end{lem}

\begin{proof}
Given a homomorphism $\rho: V_{\kappa}(\g) \arr V$, take
the images of the generating vertex operators $J^a(z)$ of
$V_{\kappa}(\g)$ under $\rho$. The fact that $\rho$ is
a homomorphism of vertex algebras implies that the OPEs between these
fields, and hence the commutation relations between their Fourier
coefficients, are preserved.

Conversely, suppose that we are given vertex operators $\wt{J}^a(z)$
satisfying the condition of the lemma.  Denote by $\wt{J}^a_n$ the
corresponding Fourier coefficients. Define a linear map $\rho:
V_{\kappa}(\g) \arr V$ by the formula
$$
J^{a_1}_{n_1} \ldots J^{a_m}_{n_m} v_{\kappa} \mapsto
\wt{J}^{a_1}_{n_1} \ldots \wt{J}^{a_m}_{n_m} \vac.
$$
It is easy to check that this map is a vertex algebra homomorphism.
\end{proof}

The complexes $C^\bullet(L_+\g,M_{\g,+})$ and
$C^\bullet_{\on{loc}}(L\g,\cA^\g_{0,\on{loc}})$ carry natural
$\Z$--gradations defined by the formulas $\deg a^*_{\al,n} = \deg
\psi^*_{a,n} = -n$, $\deg \vac = 0$, and $\deg \int Y(A,z) dz = \deg A
- 1$. The differentials preserve these gradations. The following is a
corollary of \thmref{two eq}.

\begin{cor}    \label{exists}
There exists a homomorphism of vertex algebras $V\cri(\g) \to
M_\g$.
\end{cor}

\begin{proof}
According to \thmref{two eq}, there exists a homomorphism
$\ghat_{\ka_c} \to \cA^\g_{\leq 1,\on{loc}}$ such that $K \mapsto 1$
and
$$
J^a(z) \mapsto Y(P_{a}(\astar_{\alpha,0},a_{\beta,-1}) \vac,z) +
Y(B_a,z),
$$
where the vectors $P_a$ and $B_a$ have degree one. Therefore by
\lemref{hom from V} we obtain a homomorphism of vertex algebras
$V\cri(\g) \to M_\g$ such that
$$
J^a_{-1} v_\ka \mapsto
P_{a}(\astar_{\alpha,0},a_{\beta,-1})\vac + B_a.
$$
\end{proof}

In addition to the $\Z$--gradations described above, the complex
$C^\bullet_{\on{loc}}(L\g,\cA^\g_{0,\on{loc}})$ also carries a weight
gradation with respect to the root lattice of $\g$ such that $\on{wt}
a^*_{\al,n} = - \al$, $\on{wt} \psi^*_{a,n} = - \on{wt} J^a$, $\on{wt}
\vac = 0$, and $\on{wt} \int Y(A,z) dz = \on{wt} A$. The differentials
preserve this gradation, and it is clear that the cocycles $\omega$
and $i_*(\sigma)$ are homogeneous. Therefore the element $\gamma$
introduced in the proof of \corref{exists} may be chosen in such a way
that it is also homogeneous with respect to the weight gradation.

For such $\gamma$ we necessarily have $B_a = 0$ for all $J^a \in \h
\oplus \n_+$. Furthermore, the term $B_a$ corresponding to $J^a = f_i$
must be proportional to $a^*_{\al_i,-1} \vac$. Using the formulas of
\secref{sec on formulas} and the discussion of \secref{action of Lg},
we therefore obtain a more explicit description of the homomorphism
$V\cri(\g) \to M_\g$:

\begin{thm}    \label{ci}
There exist constants $c_i \in \C$ such that the Fourier coefficients
of the vertex operators
\begin{align*}
e_i(z) &= a_{\al_i}(z) + \sum_{\beta \in \Delta_+}
\Wick P^i_\beta(a^*_\al(z)) a_\beta(z) \Wick , \\
h_i(z) &= - \sum_{\beta \in \Delta_+} \beta(h_i) \Wick a^*_\beta(z)
a_\beta(z) \Wick \, , \\
f_i(z) &= \sum_{\beta \in \Delta_+} \Wick
Q^i_\beta(a^*_\al(z)) a_\beta(z) \Wick + c_i \pa_z a^*_{\al_i}(z),
\end{align*}
where the polynomials $P^i_\beta, Q^i_\beta$ are introduced in
formulas \eqref{formulas1}--\eqref{formulas3}, generate an action of
$\ghat_{\ka_c}$ on $M_\g$.
\end{thm}

\begin{remark}    \label{right action}
In addition to the above homomorphism of Lie algebras $\ww_{\ka_c}:
\ghat_{\kappa_c} \to \cA^\g_{\leq 1,\on{loc}}$, there is also a Lie
algebra anti-homomorphism $$w^R: L\n_+ \to \cA^\g_{\leq 1,\on{loc}}$$
which is induced by the right action of $\n_+$ on $N_+$ (see
\secref{sec on formulas}). By construction, the images of $L\n_+$
under $\ww_{\ka_c}$ and $w^R$ commute. We have
\begin{equation}    \label{e right}
w^R(e_i(z)) = a_{\al_i}(z) + \sum_{\beta \in \Delta_+}
P^{R,i}_\beta(a^*_\al(z)) a_\beta(z),
\end{equation}
where the polynomials $P^{R,i}_\beta$ are defined in \secref{sec
on formulas}. More generally, we have
\begin{equation}    \label{e al right}
w^R(e_\al(z)) = a_{\al}(z) + \sum_{\beta \in \Delta_+}
P^{R,\al}_\beta(a^*_\ga(z)) a_\beta(z),
\end{equation}
for some polynomials $P^{R,\al}_\beta$.\qed
\end{remark}

\subsection{Other $\ghat$--module structures on $M_\g$}

In \thmref{ci} we constructed the structure of a
$\ghat_{\ka_c}$--module on $M_\g$. To obtain other
$\ghat_{\ka_c}$--module structures of critical level on $M_\g$ we need
to consider other homomorphisms $\ghat_{\ka_c} \to \cA^\g_{\leq
  1,\on{loc}}$ lifting the homomorphism $L\g \to \T^\g_{\on{loc}}$.
According to \lemref{L:diag}, the set of isomorphism classes of such
liftings is a torsor over $H^1(L\g,\cA^\g_{0,\on{loc}})$. Since we
have chosen a particular lifting in \thmref{ci}, we may identify this
set with $H^1(L\g,\cA^\g_{0,\on{loc}})$.

Recall that our complex
$C^\bullet_{\on{loc}}(L\g,\cA^\g_{0,\on{loc}})$ has a weight
gradation, and our cocycle $\omega$ has weight zero. Therefore among
all liftings we consider those which have weight zero. The set of such
liftings is in bijection with the weight zero homogeneous component of
$H^1(L\g,\cA^\g_{0,\on{loc}})$.

\begin{lem}    \label{family}
The weight $0$ component of $H^1(L\g,\cA^\g_{0,\on{loc}})$
is isomorphic to the (topological) dual space $(L\h)^*$ to $L\h$.
\end{lem}

\begin{proof}

First we show that the weight $0$ component of
$H^1(L\g,\cA^\g_{0,\on{loc}})$ is isomorphic to the space of weight
$0$ cocycles in $C^1(L\g,\cA^\g_{0,\on{loc}})$. Indeed, since $\on{wt}
a^*_{\al,n} = - \al$, the weight $0$ part of
$C^0_{\on{loc}}(L\g,\cA^\g_{0,\on{loc}}) = \cA^\g_{0,\on{loc}}$ is
one-dimensional and consists of the constants. If we apply the
differential to a constant we obtain $0$, and so there are no
coboundaries of weight $0$ in $C^1(L\g,\cA^\g_{0,\on{loc}})$.

Next, we show that any weight $0$ one-cocycle $\phi$ is uniquely
determined by its restriction to $L\h \subset L\g$, which may be an
arbitrary (continuous) linear functional on $L\h$. Indeed, since the
weights occurring in $\cA^\g_{0,\on{loc}}$ are less than or equal to
$0$, the restriction of $\phi$ to $L_+\n_-$ is equal to $0$, and the
restriction to $L\h$ takes values in the constants $\C \subset
\cA^\g_{0,\on{loc}}$. Now let us fix $\phi|_{L\h}$. We identify
$(L\h)^*$ with $\h^*((t)) dt$ using the residue pairing, and write
$\phi|_{L\h}$ as $\chi(t) dt$ using this identification. Here
$$\chi(t) = \sum_{n \in \Z} \chi_n t^{-n-1}, \qquad \chi_n \in \h^*,$$
where $\chi_n(h) = \phi(h_n)$. We denote $\lan \chi(t),h_i \ran$ by
$\chi_i(t)$.

We claim that for any $\chi(t) dt \in \h^*((t)) dt$, there is a unique
one-cocycle $\phi$ of weight $0$ in
$C^1_{\on{loc}}(L\g,\cA^\g_{0,\on{loc}})$ such that
\begin{align}    \label{ez}
\phi(e_i(z)) &= 0, \\ \label{hz}
\phi(h_i(z)) &= \chi_i(z), \\ \label{fz}
\phi(f_i(z)) &= \chi_i(z) a^*_{\al,n}(z).
\end{align}
Indeed, having fixed $\phi(e_i(z))$ and $\phi(h_i(z))$ as in
\eqref{ez} and \eqref{hz}, we obtain using the formula
$$\phi\left([e_{i,n},f_{j,m}]\right) =
e_{i,n} \cdot \phi(f_{j,m}) - f_{j,m} \cdot \phi(e_{i,n})
$$
that
$$
\delta_{i,j} \chi_{i,n+m} = e_{i,n} \cdot \phi(f_{j,m}).
$$
This equation on $\phi(f_{j,m})$ has a unique solution in
$\cA^\g_{0,\on{loc}}$ of weight $-\al_i$, namely, the one given by
formula \eqref{fz}. The cocycle $\phi$, if exists, is uniquely
determined once we fix its values on $e_{i,n}, h_{i,n}$ and $f_{i,n}$.
Let us show that it exists. This is equivalent to showing that the
Fourier coefficients of the fields
\begin{align}    \label{act1}
e_i(z) &= a_{\al_i}(z) + \sum_{\beta \in \Delta_+}
\Wick P^i_\beta(a^*_\al(z)) a_\beta(z) \Wick , \\ \label{act2}
h_i(z) &= - \sum_{\beta \in \Delta_+} \beta(h_i) \Wick a^*_\beta(z)
a^*_\beta(z) \Wick + \chi_i(z), \\ \label{act3}
f_i(z) &= \sum_{\beta \in \Delta_+} \Wick
Q^i_\beta(a^*_\al(z)) a_\beta(z) \Wick + c_i \pa_z a^*_{\al_i}(z) +
\chi_i(z) a^*_{\al,n}(z),
\end{align}
satisfy the relations of $\ghat_{\ka_c}$ with $K=1$.

Let us remove the normal ordering and set $c_i=0,
i=1,\ldots,\ell$. Then the corresponding Fourier coefficients are no
longer well-defined as linear operators on $M_\g$. But they are
well-defined linear operators on the space $\on{Fun} LU$. The
resulting $L\g$--module structure is easy to describe. Indeed, the Lie
algebra $L\g$ acts on $\on{Fun} LU$ by vector fields. More generally,
for any $L\bb_-$--module $R$ we obtain a natural action of $L\g$ on
the tensor product $\on{Fun} LU \wh{\otimes} R$ (this is just the
topological $L\g$--module induced from the $L\bb_-$--module $R$). If
we choose as $R$ the one-dimensional representation on which all
$f_{i,n}$ act by $0$ and $h_{i,n}$ acts by multiplication by
$\chi_{i,n}$ for all $i=1,\ldots,\ell$ and $n \in \Z$, then the
corresponding $L\g$--action on $\on{Fun} LU \wh{\otimes} R \simeq
\on{Fun} LU$ is given by the Fourier coefficients of the above
formulas, with the normal ordering removed. Hence if we remove the
normal ordering and set $c_i=0$, then these Fourier coefficients do
satisfy the commutation relations of $L\g$.

When we restore the normal ordering, these commutation relations may
in general be distorted, due to the double contractions, as explained
in \secref{double contractions}. But we know from \thmref{ci} that
when we restore normal ordering and set all $\chi_{i,n}=0$, then there
exist the numbers $c_i$ such that these Fourier coefficients satisfy
the commutation relations of $\ghat_{\ka_c}$ with $K=1$. The terms
\eqref{ez}--\eqref{fz} will not generate any new double contractions
in the commutators. Therefore computing these commutators with
non-zero values of $\chi_{i,n}$, we find that the relations remain the
same. This completes the proof.
\end{proof}

\begin{cor}    \label{wchit}
For each $\chi(t) \in \h^*((t))$ there is a $\ghat$--module structure
of critical level on $M_\g$, with the action given by formulas
\eqref{act1}--\eqref{act3}.
\end{cor}

We call these modules the {\em Wakimoto modules of critical level} and
denote them by $W_{\chi(t)}$.

Let $\pi_0$ be the commutative algebra $\C[b_{i,n}]_{i=1,\ldots,\ell;
n<0}$ with the derivation $T$ given by the formula
$$
T \cdot b_{i_1,n_1} \ldots b_{i_m,n_m} = - \sum_{j=1}^m n_j
b_{i_1,n_1} \ldots b_{i_j,n_j-1} \cdots b_{i_m,n_m}.
$$
Then $\pi_0$ is naturally a commutative vertex algebra (see
\cite{vertex}, \S~2.3.9). In particular, we have
$$
Y(b_{i,-1},z) = b_i(z) = \sum_{n<0} b_{i,n} z^{-n-1}.
$$
Using the same argument as in the proof of \lemref{family}, we now
obtain a stronger version of \corref{exists}.

\begin{thm}    \label{exist hom}
There exists a homomorphism of vertex algebras $$\ww_{\ka_c}: V\cri(\g)
\to M_\g \otimes \pi_0$$ such that
\begin{align*}
e_i(z) &\mapsto a_{\al_i}(z) + \sum_{\beta \in \Delta_+}
\Wick P^i_\beta(a^*_\al(z)) a_\beta(z) \Wick , \\
h_i(z) &\mapsto - \sum_{\beta \in \Delta_+} \beta(h_i) \Wick a^*_\beta(z)
a_\beta(z) \Wick + b_i(z) , \\
f_i(z) &\mapsto \sum_{\beta \in \Delta_+} \Wick
Q^i_\beta(a^*_\al(z)) a_\beta(z) \Wick + c_i \pa_z a^*_{\al_i}(z) +
b_i(z) a^*_{\al_i}(z),
\end{align*}
where the polynomials $P^i_\beta, Q^i_\beta$ are introduced in
formulas \eqref{formulas1}--\eqref{formulas3}.
\end{thm}

Thus, any module over the vertex algebra $M_\g \otimes \pi_0$ becomes
a $V\cri(\g)$--module, and hence a $\ghat$--module of critical
level. We will not require that the module is necessarily
$\Z$--graded. In particular, for any $\chi(t) \in \h^*((t))$ we have a
one-dimensional $\pi_0$--module $\C_{\chi(t)}$, on which $b_{i,n}$
acts by multiplication by $\chi_{i,n}$. The corresponding
$\ghat_{\ka_c}$--module is the Wakimoto module $W_{\chi(t)}$
introduced above.

\subsection{Application: proof of the Kac--Kazhdan conjecture}
\label{KK}

As an application of the above construction of Wakimoto modules, we
give a proof of the Kazhdan--Kazhdan conjecture from \cite{KK},
following \cite{F:thesis}.

Let us recall the notion of a character of a $\ghat_\ka$--module.
Suppose that we have a $\ghat_\ka$--module $M$ equipped with an action
of the gradation operator $L_0 = - t\pa_t$. Note that if $\ka \neq
\ka_c$, then any smooth $\ghat_\ka$--module (i.e., such that any
vector is annihilated by the Lie subalgebra $\g \otimes t^N \C[[t]]$
for sufficiently large $N$) carries an action of the Virasoro algebra
obtained via the Segal--Sugawara construction (see \secref{conformal
str} below), and so in particular an $L_0$--action. However, smooth
$\ghat_{\ka_c}$--modules do not necessarily carry an $L_0$--action.

Suppose in addition that $L_0$ and $\h \otimes 1 \subset \ghat_\ka$
act diagonally on $M$ with finite-dimensional common eigenspaces. Then
we define the character of $M$ as the formal series
\begin{equation}    \label{character}
\on{ch} M = \sum_{\wh{\la} \in (\h \oplus \C L_0)^*} \dim M(\wh\la) \;
e^{\wh{\la}},
\end{equation}
where $M(\wh\la)$ is the generalized eigenspace of $L_0$ and $\h
\otimes 1$ corresponding to $\wh\la: (\h \oplus \C L_0)^* \to \C$.

The direct sum $(\h \otimes 1) \oplus \C L_0 \oplus \C K$ is the
Cartan subalgebra of the extended Kac-Moody algebra $\ghat_\ka \oplus
\C L_0$. Elements of the dual space $\wt{\h}^*$ are called weights. We
will consider weights occurring in modules on which $K$ acts as the
identity. Therefore without loss of generality we may view weights as
elements of the dual space to $\wt{\h} = (\h \otimes 1) \oplus \C
L_0$. The set of positive roots of $\ghat$ is naturally a subset of
$\wt{\h}^*$, and this determines a natural partial order on the set of
weights. In what follows we will often denote $e^{-\delta} =
e^{-(0,1)}$, where $\delta$ is the imaginary root of $\g$, by $q$.

For $\wh{\la} = (\la,a)$, let $L_{\wh{\la},\ka}$ be the irreducible
$\ghat_\ka$--module with the highest weight $\wh{\la} = (a,\la)$, on
which the central element $K$ acts as the identity. It is the unique
irreducible subquotient of the Verma module $M_{\wh{\la},\ka}$. In
\cite{KK} a certain subset $H^\ka_{\beta,m} \in \wt{\h}^*$ is defined
for any pair $(\beta,m)$, where $\beta$ is a positive root of $\ghat
\oplus \C L_0$ and $m$ is a positive integer. If $\beta$ is a real
root, then $H^\ka_{\beta,m}$ is a hyperplane in $\wt{\h}^*$ and if
$\beta$ is an imaginary root, then $H^{\ka_c}_{\beta,m} = \wt{\h}^*$
and $H^{\ka}_{\beta,m} = \emptyset$ for $\ka \neq \ka_c$. It is shown
in \cite{KK} that $L_{\wh{\la},\ka}$ is a subquotient of
$M_{\wh{\mu},\ka}$ if and only if the following condition is
satisfied: there exists a finite sequence of weights
$\wh{\mu}_0,\ldots,\wh{\mu}_n$ such that $\wh{\mu}_0 = \wh{\la},
\wh{\mu}_n = \wh{\mu}$, $\wh{\mu}_{i+1} = \wh{\mu}_i - m_i \beta_i$
for some positive roots $\beta_i$ and positive integers $m_i$, and
$\wh{\mu}_i \in H^\ka_{\beta_i,m_i}$ for all $i=1,\ldots,n$.

Denote by $\wh\De^{\on{re}}_+$ the set of positive real roots of
$\ghat$. Let us call a weight $\wh{\la}$ a {\em generic weight of
critical level} if $\wh{\la}$ does not belong to any of the
hyperplanes $H^\ka_{\beta,m}, \beta \in \wh\De^{\on{re}}_+$, where
$\wh\De^{\on{re}}_+$ is the set of positive real roots. It is easy to
see from the above condition that $\wh{\la}$ is a generic weight of
critical level if and only if the only irreducible subquotients of
$M_{\wh{\la},\ka_c}$ have highest weights $\wh{\la} - n\delta$, where
$n$ is a non-negative integer (i.e., their $\h^*$ components are equal
to the $\h^*$ component of $\wh{\la}$). The following assertion is the
Kac--Kazhdan conjecture for the untwisted affine Kac-Moody algebras.

\begin{thm}    \label{KK conj}
For generic weight $\wh{\la}$ of critical level
$$
\on{ch} L_{\wh{\la},\ka_c} = e^{\wh{\la}} \prod_{\al \in
  \wh\De^{\on{re}}_+} (1-e^{-\al})^{-1}.
$$
\end{thm}

\begin{proof}
Without loss of generality, we may assume that $\wh{\la} = (\la,0)$.

  Introduce the gradation operator $L_0$ on the Wakimoto module
  $W_{\chi(t)}$ by using the vertex algebra gradation on $M_\g$. It is
  clear from the formulas defining the $\ghat_{\ka_c}$--action on
  $W_{\chi(t)}$ given in \thmref{exist hom} that this action is
  compatible with the gradation if and only if $\chi(t) = \la/t$,
  where $\la \in \h^*$. In that case
$$
\on{ch} W_{\la/t} = e^{\wh{\la}} \prod_{\al \in
  \wh\De^{\on{re}}_+} (1-e^{-\al})^{-1},
$$ where $\wh{\la} = (\la,0)$. Thus, in order to prove the theorem we
need to show that if $\wh{\la}$ is a generic weight of critical level,
then $W_{\la/t}$ is irreducible. Suppose that this is not so. Then
either $W_{\la/t}$ contains a singular vector, i.e., a vector
annihilated by the Lie subalgebra $\wt\n_+ = (\g \otimes t\C[[t]])
\oplus (\n_+ \otimes 1)$, other than the multiples of the highest
weight vector, or $W_{\la/t}$ is not generated by its highest weight
vector.

Suppose that $W_{\la/t}$ contains a singular vector other than a
multiple of the highest weight vector. Such a vector must then be
annihilated by the Lie subalgebra $L_+\n_+ = \n_+[[t]]$. We have
introduced in \remref{right action} the right action of $\n_+((t))$ on
$M_\g$, which commutes with the left action. It is clear from formula
\eqref{e al right} that the monomials
\begin{equation} \label{monn}
\prod_{n_{r(\al)}<0} e^R_{\al,n_{r(\al)}}
\prod_{m_{s(\al)}\leq 0} a^*_{\al,m_{s(\al)}} \vac
\end{equation}
form a basis of $M_\g$. Therefore the space of $L_+\n_+$--invariants
of $W_{\ka_c,0}$ is equal to the tensor product of the subspace
$M_{\g,-}$ of $W_{0,\ka_c}$ spanned by all monomials \eqref{monn} not
containing $a^*_{\al,n}$, and the space of $L_+\n_+$--invariants in
$M_{\g,+} = \C[a^*_{\al,n}]_{\al \in \De_+,n\leq 0}$. According to
\secref{M+}, $M_{\g,+}$ is an $L_+\g$--module isomorphic to
$\on{Coind}_{L_+\bb_-}^{L_+\g} \C$.  Therefore the action of $L_+\n_+$
on it is co-free, and the space of $L_+\n_+$--invariants is
one-dimensional, spanned by constants.  Thus, we obtain that the space
of $L_+\n_+$--invariants in $W_{\la/t}$ is equal to $M_{\g,-}$. In
particular, we find that the weight of any singular vector of
$W_{\la/t}$ which is not equal to the highest weight vector has the
form $(\la,0) - \sum_j (n_j\delta - \beta_j)$, where $n_j > 0$ and
each $\beta_j$ is a positive root of $\g$. But then $W_{\la/t}$
contains an irreducible subquotient of such a weight.

Now observe that
$$
\on{ch} M_{(\la,0),\ka_c} = \prod_{n>0} (1-q^n)^{-\ell} \cdot \on{ch}
W_{\la/t},
$$ where $q=e^{-\delta}$. If an irreducible module
$L_{\wh{\mu},\ka_c}$ appears as a subquotient of $W_{\la/t}$, then it
appears in the decomposition of $\on{ch} W_{\la/t}$ into the sum of
characters of irreducible representations and hence in the
decomposition of $\on{ch} M_{(\la,0),\ka_c}$. Since the characters of
irreducible representations are linearly independent, this implies
that $L_{\wh{\mu},\ka_c}$ is an irreducible subquotient of
$M_{(\la,0),\ka_c}$. But this contradicts our assumption that
$(\la,0)$ is a generic weight of critical level.  Therefore
$W_{\la/t}$ does not contain any singular vectors other than the
multiples of the highest weight vector.

Now suppose that $W_{\la/t}$ is not generated by its highest weight
vector. But then there exists a homogeneous linear functional on
$W_{\la/t}$, whose weight is less than the highest weight and which is
invariant under $\wt\n_- = (\g \otimes t^{-1} \C[t^{-1}]) \oplus (\n_-
\otimes 1)$, and in particular, under its Lie subalgebra $L_- \n_+ =
\n_+ \otimes t^{-1} \C[t^{-1}]$. Therefore this functional factors
through the space of coinvariants of $W_{\la/t}$ by $L_- \n_+$. But
$L_- \n_+$ acts freely on $W_{\la/t}$, and the space of coinvariants
is isomorphic to the subspace $\C[a_{\al,n}^*]_{\al\in\De_+,n\leq 0}$
of $W_{\la/t}$. Hence we obtain that the weight of this functional has
the form $(\la,0) - \sum_j (n_j\delta + \beta_j)$, where $n_j\geq 0$
and each $\beta_j$ is a positive root of $\g$.  In the same way as
above, it follows that this contradicts our assumption that $\la$ is a
generic weight. Therefore $W_{\la/t}$ is generated by its highest
weight vector. This completes the proof.
\end{proof}

\subsection{Coordinate--independent version}    \label{coord-indep
version}

In the above constructions we considered representations of the Lie
algebra $\ghat_{\ka_c}$, which is the central extension of $L\g =
\g((t))$. In applications, it is important to develop a theory which
applies to the central extension of the Lie algebra $\g(\K) = \g
\otimes \K$, where $\K$ is a topological algebra which is isomorphic
to $\C((t))$, but non-canonically. For instance, we can take as $\K$
the field of fractions of the completed local ring of a point $x$ of a
smooth curve $X$. Then if we choose a formal coordinate $t$ at $x$, we
may identify $\K$ with $\C((t))$, but this identification is
non-canonical as there is usually no preferred choice of coordinate
$t$. Formula \eqref{KM rel} defining the central extension of
$\g((t))$ is independent of the choice of $t$, and so the central
extension is well-defined for any algebra $\K$ as above. We will
denote it by $\ghat_{\ka_c}(\K)$ to emphasize this fact.

In order to recast the above construction of Wakimoto modules in a
coordinate--independent way, we need to understand how to incorporate
the action of the group of changes of coordinates into this
construction. Let $\AutO$ be the group of continuous automorphisms of
$\C[[t]]$. Any such automorphism is determined by its action on the
topological generator $t$ of $\C[[t]]$, $t \mapsto \rho(t)$. This
allows us to identify $\AutO$ with the group of formal power series
$\rho(t) = \rho_1 t + \rho_2 t^2 + \ldots$, where $\rho_1 \neq 0$ (see
\cite{vertex}, \S~5.1, for more details). The Lie algebra of $\AutO$
is $\DerpO = t\C[[t]] \pa_t$. Denote by $\DerO$ the Lie algebra
$\C[[t]] \pa_t$. Then $(\DerO,\AutO)$ is an example of a
Harish-Chandra pair (see \cite{vertex}, Ch.~5 and \S~16.2).

Now let $\OO \subset \K$ be a pair consisting of a complete local ring
$\OO$ isomorphic to $\C[[t]]$ and its field of fractions $\K$. Let
${\mf m}$ be the unique maximal ideal of $\OO$, and $\Au$ the set
of topological generators of ${\mf m}$ (these are the coordinates on
the disc $\on{Spec} \OO$). Then $\Au$ is naturally an $\AutO$--torsor,
with the (right) action given by the formula $z \mapsto \rho(z)$ for
each $\rho(t) \in \AutO$. Given any $\AutO$--module $V$, we may form
its twist by $\Au$,
$$
\V \overset{\on{def}}{=} \Au \underset{\AutO}\times V.
$$
Let $V$ be the module $V_\kappa(\g)$ defined in \secref{hom of VA}
(where we identify $L\g$ with $\g((t))$). Since $\g((t))$ and
$\g[[t]]$ carry natural $\AutO$--actions, so does $V_\kappa(\g)$. The
corresponding $\Au$--twist $\V_\kappa(\g)$ may alternatively be
described as follows:
\begin{equation}    \label{mcV}
\V_\kappa(\g) = \Au \underset{\AutO}\times V_\ka(\g) =
\on{Ind}^{\ghat_{\ka_c}(\K)}_{\g \otimes \OO \oplus \C K} \C_1,
\end{equation}
where $\ghat_{\ka_c}(\K)$ is now the central extension of $\g \otimes
\K$ rather than $\g((t))$.

We now want to describe an action of $\AutO$ on $M_\g$ and $\pi_0$
such that the homomorphism $V\cri(\g) \to M_\g \otimes \pi_0$ of
\thmref{exist hom} commutes with the action of $\AutO$. Then we will
obtain a map $\V_\kappa(\g)$ to the corresponding $\Au$--twist of
$M_\g \otimes \pi_0$ and this will enable us to describe the Wakimoto
modules in a coordinate--independent fashion.

Before describing the $\AutO$--action on $M_\g$ and $\pi_0$, we
discuss in the next section how to deform the homomorphism $V\cri(\g)
\to M_\g \otimes \pi_0$ away from the critical level.

\section{Deforming to other levels}

\subsection{Homomorphism of vertex algebras}    \label{pi0}

As before, we denote by $\h$ the Cartan subalgebra of $\g$. Let
$\hh_\ka$ be the one-dimensional central extension of the loop algebra
$L\h = \h((t))$ with the two-cocycle obtained by restriction of the
two-cocycle on $L\g$ corresponding to the inner product $\kappa$. Then
according to formula \eqref{KM rel}, $\hh_\ka$ is a Heisenberg Lie
algebra. We will consider a copy of this Lie algebra with generators
$b_{i,n}, i=1,\ldots,\ell, n \in {\mathbb Z}$, and $K$ with the
commutation relations
$$
[b_{i,n},b_{j,m}] = n \ka(h_i,h_j) K \delta_{n,-m}.
$$
Thus, the $b_{i,n}$'s satisfy the same relations as the
$h_{i,n}$'s.  Let $\pi^\kappa_{0}$ denote the $\hh_\ka$--module
induced from the one-dimensional representation of the abelian Lie
subalgebra of $\hh_\ka$ spanned by $b_{i,n}, i=1,\ldots,\ell, n \geq
0$, and $K$, on which $K$ acts as the indentity and all other
generators act by $0$. We denote by $\vac$ the generating vector of
this module. It satisfies: $b_{i,n} \vac = 0, n\geq 0$. Then
$\pi^\ka_{0}$ has the following structure of a vertex algebra (see
Theorem 2.3.7 of \cite{vertex}):

\begin{enumerate}

\item[$\bullet$] Gradation: $\on{deg} b_{i_1,n_1} \ldots b_{i_m,n_m}
\vac = -\sum_{i=1}^m {n_i}$.

\item[$\bullet$] Vacuum vector: $\vac$.

\item[$\bullet$] Translation operator: $T \vac = 0$,
$[T,b_{i,n}]=-n b_{i,n-1}$.

\item[$\bullet$] Vertex operators:
$$Y(b_{i,-1} \vac,z) = b_i(z) = \sum_{n\in\Z} b_{i,n}
z^{-n-1},$$
\begin{equation*}
Y(b_{i_1,n_1}\dots b_{i_m,n_m} \vac,z) = \prod_{j=1}^n
\frac{1}{(-n_j-1)!} \; \Wick \delz^{-n_1-1}
b_{i_1}(z)\dots\delz^{-n_m-1} b_{i_m}(z) \Wick \; .
\end{equation*}
\end{enumerate}

The tensor product $M_\g \otimes \pi^{\ka-\ka_c}_0$ also acquires a
vertex algebra structure.

\begin{thm}    \label{noncrit}
There exists a homomorphism of vertex algebras $$\ww_\ka: V_\kappa(\g)
\to M_\g \otimes \pi^{\kappa-\kappa_c}_0$$ such that
\begin{align} \notag
e_i(z) &\mapsto a_{\al_i}(z) + \sum_{\beta \in \Delta_+} \Wick
P^i_\beta(a^*_\al(z)) a_\beta(z) \Wick \; , \\ \label{in the stat}
h_i(z) &\mapsto - \sum_{\beta \in \Delta_+} \beta(h_i) \Wick
a^*_\beta(z) a_\beta(z) \Wick + b_i(z), \\ \notag f_i(z) &\mapsto
\sum_{\beta \in \Delta_+} \Wick Q^i_\beta(a^*_\al(z)) a_\beta(z) \Wick
+ \left( c_i+(\kappa-\kappa_c)(e_i,f_i) \right) \pa_z a^*_{\al_i}(z) +
b_i(z) a^*_{\al_i}(z),
\end{align}
where the polynomials $P^i_\beta, Q^i_\beta$ are introduced in
formulas \eqref{formulas1}--\eqref{formulas3}.
\end{thm}

\begin{proof}
  Denote by $\wt{\cA}^\g_{\on{loc}}$ the Lie algebra $U(M_\g \otimes
  \pi^{\ka-\ka_c}_0)$. By \lemref{hom from V}, in order to prove the
  theorem, we need to show that formulas \eqref{in the stat} define a
  homomorphism of Lie algebras $\ghat_\kappa \to
  \wt{\cA}^\g_{\on{loc}}$ sending the central element $K$ to the
  identity.
  
  Formulas \eqref{in the stat} certainly define a linear map
  $\ol{w}_\ka: L\g \to \wt{\cA}^\g_{\on{loc}}$. Denote by $\omega_\ka$
  the linear map $\bigwedge^2 L\g \to \wt{\cA}^\g_{\on{loc}}$ defined
  by the formula
$$
\omega_\ka(f,g) = [\ol{w}_\ka(f),\ol{w}_\ka(g)] - \ol{w}_\ka([f,g]).
$$
Evaluating it explicitly in the same way as in the proof of
\lemref{repr in loc}, we find that $\omega_\ka$ takes values in
$\cA^\g_{0,\on{loc}} \subset \wt{\cA}^\g_{\on{loc}}$. Furthermore, by
construction of $\ol{w}_\ka$, for any $X \in \cA^\g_{0,\on{loc}}$ and
$f \in L\g$ we have $[\ol{w}_\ka(f),X] = f \cdot X$, where in the
right hand side we consider the action of $f$ on the $L\g$--module
$\cA^\g_{0,\on{loc}}$.  This immediately implies that $\omega_\ka$ is a
two-cocycle of $L\g$ with coefficients in $\cA^\g_{0,\on{loc}}$. By
construction, $\omega_\ka$ is local, i.e., belongs to
$C^2_{\on{loc}}(L\g,\cA^\g_{0,\on{loc}})$.

Let us compute the restriction of $\omega_\ka$ to $\bigwedge^2 L\h$. The
calculations made in the proof of \lemref{two eq} imply that
$$
\omega_\ka(h_n,h'_m) = n(\ka_c(h,h') + (\ka-\ka_c)(h,h')) = n\ka(h,h').
$$
Therefore this restriction is equal to the restriction of the
Kac-Moody two-cocycle $\sigma_\ka$ on $L\g$ corresponding to
$\kappa$. Now \lemref{any two} implies that the two-cocycle
$\omega_\ka$ is cohomologically equivalent to $i_*(\sigma_\ka)$. We
claim that it is actually equal to $i_*(\sigma_\ka)$.

Indeed, the difference between these cocycles is the coboundary of
some element $\gamma \in C^1_{\on{loc}}(L\g,\cA^\g_{0,\on{loc}})$. The
discussion before \thmref{ci} implies that $\gamma(e_i(z)) =
\gamma(h_i(z)) = 0$ and $\gamma(f_i(z)) = c'_i \pa_z a^*_{\al_i}(z)$
for some constants $c'_i \in \C$. In order to find the constants
$c'_i$ we compute the value of the corresponding two-cocycle
$\omega_\ka + d\gamma$ on $e_{i,n}$ and $f_{i,-n}$. We find that it is
equal to $n\sigma_\ka(e_{i,n},f_{i,-n}) + c'_i n (e_i,f_i)$. Therefore
$c'_i=0$ for all $i=1,\ldots,\ell$, and so $\gamma=0$ and $\omega_\ka
= i_*(\sigma_\ka)$. This implies that formulas \eqref{in the stat}
indeed define a homomorphism of Lie algebras $\ghat_\kappa \to
\cA^\g_{\leq 1,\on{loc}}$ sending the central element $K$ to the
identity. This completes the proof.
\end{proof}

The following result is useful in applications.

\begin{prop}    \label{injective}
The homomorphism $\ww_\ka$ of \thmref{noncrit} is injective for any
$\ka$.
\end{prop}

\begin{proof}
  We will introduce filtrations on $V_\kappa(\g)$ and $W_{0,\ka} =
  M_\g \otimes \pi^{\kappa-\kappa_c}_0$ which are preserved by
  $\ww_\ka$, and then show that the induced map $\gr \ww_\ka: \gr
  V_\kappa(\g) \to \gr W_{0,\ka}$ (which turns out to be independent
  of $\ka$) is injective.

In order to define a filtration on $V_\kappa(\g)$ we observe that
$$
V_\kappa(\g) \simeq U(\g \otimes t^{-1}\C[t^{-1}]) v_\ka
$$
and use the Poincar\'e--Birkhoff--Witt filtration on $U(\g \otimes
t^{-1}\C[t^{-1}])$. Now we define a filtration $\{ W^{\leq p}_{0,\ka}
\}$ on $W_{0,\ka}$ by defining $W^{\leq p}_{0,\ka}$ to be the span of
monomials in the $a_{\al,n}$'s, $a^*_{\al,n}$'s and $b_{i,n}$'s whose
combined degree in the $a_{\al,n}$'s and $b_{i,n}$'s is less than or
equal to $p$ (this is analogous to the filtration by the order of
differential operator). It is clear from the construction of the
homomorphism $\ww_\ka$ that it preserves these filtrations.

Now we describe the corresponding operator $\gr \ww_\ka: \gr
V_\kappa(\g) \to \gr W_{0,\ka}$. Let $\wt{\g}$ be the variety of pairs
$(\bb,x)$, where $\bb$ is a Borel subalgebra in $\g$ and $x \in \bb$.
The natural morphism $\wt{\g} \to \on{Fl}$, where $\on{Fl}$ is the
flag variety of $\g$, mapping $(\bb,x)$ to $\bb \in \on{Fl}$
identifies $\wt{\g}$ with a vector bundle over the flag variety
$\on{Fl}$, whose fiber over $\bb \in \on{Fl}$ is the vector space
$\bb$.  There is also a morphism $\wt{\g} \to \g$ sending $(\bb,x)$ to
$x$.  Now let $U$ be again the big cell, i.e., the open $N_+$--orbit
of $\on{Fl}$, and $\wt{U}$ its preimage in $\wt{\g}$. In other words,
$\wt{U}$ consists of those pairs $(\bb,x)$ for which $\bb$ is in
generic relative position with $\bb_+$. It is isomorphic to an affine
space of dimension equal to $\dim \g$. The induced morphism $p: \wt{U}
\to \g$ is dominant and generically one-to-one.

Now let $J\wt{U}$ and $J\g$ be the infinite jet schemes of $\wt{U}$
and $\g$, defined as in \cite{vertex}, \S~8.4.4, and $Jp$ the
corresponding morphism $J\wt{U} \to J\g$. Then $Jp$ is clearly a
dominant morphism and so the corresponding homomorphism of rings of
functions $Jp: \C[J\g] \to \C[J\wt{U}]$ is injective.

Since $J\g \simeq \g[[t]]$, it follows that $$\gr V_\kappa(\g) \simeq
\on{Sym} \g((t))/\g[[t]] \simeq \C[J\g].$$ Moreover, $J\wt{U}$ is
isomorphic to $\gr W_{0,\ka}$, and it follows from our construction of
$\ww_\ka$ that $\gr \ww_\ka = Jp$. This implies the statement of the
proposition.
\end{proof}

\subsection{Wakimoto modules away from the critical level}
\label{pi00}

Any module over the vertex algebra $M_\g \otimes \pi^{\ka-\ka_c}_0$
now becomes a $V_\kappa(\g)$--module and hence a
$\ghat_\kappa$--module (with $K$ acting as $1$). For $\la \in \h^*$,
let $\pi^{\ka-\ka_c}_{\la}$ be the Fock representation of $\hh_\ka$
generated by a vector $|\la\ran$ satisfying
$$
b_{i,n}|\la\ran = 0, \quad n>0, \qquad b_{i,0}|\la\ran =
\la(h_i)|\la\ran, \qquad K |\la\ran = |\la\ran.
$$
Then $$W_{\la,\ka} \overset{\on{def}}{=} M_\g \otimes
\pi^{\ka-\ka_c}_\la$$ is an $M_\g \otimes \pi^{\ka-\ka_c}_0$--module,
and hence a $\ghat_\ka$--module. We call it the {\em Wakimoto module of
level $\ka$ and highest weight $\la$}.

\subsection{Conformal structures at non-critical levels}
\label{conformal str}

In this section we show that the homomorphism $\ww_\ka$ of
\thmref{noncrit} is a homomorphism of conformal vertex algebras when
$\ka \neq \ka_c$. This will allow us to obtain a
coordinate--independent version of this homomorphism. By taking the
limit $\ka \to \ka_c$, we will also obtain a coordinate--independent
version of the homomorphism $\ww_{\ka_c}$.

The vertex algebra $V_\ka(\g), \ka \neq \ka_c$, has the standard
conformal structure given by the Segal--Sugawara vector
\begin{equation}    \label{omegaka}
{\mb s}_\ka = \frac{1}{2} \sum_a J^a_{-1} J_{a,-1} v_\ka,
\end{equation}
where $\{ J_a \}$ is the basis of $\g$ dual to the basis $\{ J^a
\}$ with respect to the inner product $\ka-\ka_c$ (see \cite{vertex},
\S~2.5.10, for more details). We need to calculate the image of
${\mb s}_\ka$ under $\ww_\ka$.

\begin{lem}    \label{Ska}
The image of ${\mb s}_\ka$ under $\ww_\ka$ is equal to
\begin{equation}    \label{formula for Ska}
\left( \sum_{\al \in \De_+} a_{\al,-1} a^*_{\al,-1} + \frac{1}{2}
  \sum_{i=1}^\ell b_{i,-1} b^i_{-1} - \rho_{-2} \right) \vac,
\end{equation}
where $\{ b^i \}$ is a dual basis to $\{ b_i \}$ and $\rho$ is the
element of $\h$ corresponding to $\rho \in \h^*$ under the isomorphism
induced by the inner product $(\ka-\ka_c)|_{\h}$.
\end{lem}

\begin{proof}
The vector $\ww_\ka({\mb s}_\ka)$ is of degree $2$ with respect to the vertex
algebra gradation on $M_\g \otimes \pi^{\ka-\ka_c}_0$ and of weight
$0$ with respect to the gradation by the root system such that
$\on{wt} a_{\al,n} = - \on{wt} a^*_{\al,n} = \al, \on{wt} b_{i,n} =
0$. The basis in the corresponding subspace of $M_\g \otimes
\pi^{\ka-\ka_c}_0$ consists of the monomials of the form
\begin{equation}    \label{monom1}
b_{i,-1} b_{j,-1}, \qquad b_{i,-2},
\end{equation}
\begin{equation}    \label{monom2}
a_{\al,-1} a^*_{\al,-1}, \qquad a_{\al,-2} a^*_{\al,0},
\end{equation}
\begin{equation}    \label{monom3}
a_{\al,-1} a_{\beta,-1} a^*_{\al,0} a^*_{\beta,0}, \qquad a_{\al,-1}
a^*_{\al,0} b_{i,-1},
\end{equation}
applied to the vacuum vector $\vac$.

The Fourier coefficients $L_n, n \in \Z$, of the vertex operator
$\ww_\ka({\mb s}_\ka)$ preserve the weight gradation on $M_\g \otimes
\pi^{\ka-\ka_c}_0$ and $\deg L_n = -n$ with respect to the vertex
algebra gradation on $M_\g \otimes \pi^{\ka-\ka_c}_0$. This implies
that the vectors $a_{\al,-1}\vac$ and $a^*_{\al,0} \vac, \al \in
\De_+$, are primary vectors, i.e., they are annihilated by $L_n, n>0$,
and are eigenvectors of $L_0$. We claim that $L_0$ acts by $1$ on
$a_{\al,-1}\vac$ and by $0$ on $a^*_{\al,0} \vac$.

Indeed, the vectors $J^a_{-1} v_\ka \in V_\ka(\g)$ are primary of
degree $1$. Therefore the same is true for their images in $M_\g
\otimes \pi^{\ka-\ka_c}_0$. Let $\{ e_\al \}_{\al \in \De_+}$ be a
root basis of $\n_+ \subset \g$ such that $e_{\al_i} = e_i$. Then we
have
\begin{equation}    \label{eal}
\ww_\ka(e_{\al,-1}v_\ka) = \left( a_{\al,-1} + \sum_{\beta \in \Delta_+}
P^\al_\beta(a^*_{\al,0}) a_{\beta,-1} \right) \vac,
\end{equation}
where the polynomials $P^\al_\beta$ are found from the formula of the
action of $e_\al$ on $U$:
$$
e_\al = \frac{\pa}{\pa y_{\al}} + \sum_{\beta \in \Delta_+}
P^\al_\beta(y_\al) \frac{\pa}{\pa y_\beta}.
$$
In particular, $P^{\al_i}_\beta = P^i_\beta$ considered above.  In
addition we have the following formula for $\ww_\ka(h_{i,-1}v_\ka)$
which follows from \thmref{noncrit}:
\begin{equation}    \label{wkah}
\ww_{\ka}(h_{i,-1} v_\ka) = \left( - \sum_{\beta \in \Delta_+} \beta(h_i)
a^*_{0,\beta} a_{\beta,-1} + b_{i,-1} \right) \vac.
\end{equation}
Using these formulas, we obtain that the vectors
$a_{\al,-1}\vac$ (resp., $a^*_{\al,0}\vac$) are primary vectors of
degree $1$ (resp., $0$) with respect to $\ww_\ka({\mb s}_\ka)$.  This readily
implies that the basis vectors \eqref{monom3} do not appear in the
decomposition of $\ww_\ka({\mb s}_\ka)$ and that the vectors \eqref{monom2}
enter in the combination
\begin{equation}    \label{conf vec a}
\sum_{\al \in \De_+} a_{\al,-1} a^*_{\al,-1} \vac.
\end{equation}
It remains to determine the coefficients with which the monomials
\eqref{monom1} enter the formula for $\ww_\ka({\mb s}_\ka)$.

Using formula \eqref{wkah} and our knowledge of the $aa^*$ component
of $\ww_\ka({\mb s}_\ka)$ we find the following action of the Fourier
coefficients $L_n$ of $Y(\ww_\ka({\mb s}_\ka),z)$ on the first summand of
$\ww_{\ka}(h_{i,-1} v_\ka)$:
\begin{align*}
L_0 \cdot - \sum_{\beta \in \Delta_+} \beta(h_i) a^*_{0,\beta}
a_{\beta,-1} \vac &= - \sum_{\beta \in \Delta_+} \beta(h_i)
a^*_{0,\beta} a_{\beta,-1} \vac, \\ L_1 \cdot - \sum_{\beta \in
\Delta_+} \beta(h_i) a^*_{0,\beta} a_{\beta,-1} \vac & = - \sum_{\beta
\in \Delta_+} \beta(h_i)\vac = - 2\rho(h_i) \vac.
\end{align*}
and $L_n, n>1$, act by $0$. There is a unique combination of monomials
\eqref{monom1} which, when added to \eqref{conf vec a}, makes a
conformal vector with respect to which $\ww_{\ka}(h_{i,-1} v_\ka)$ is a
vector of degree $1$ annihilated by $L_n, n>0$, namely,
\begin{equation}    \label{conf vec b}
\frac{1}{2} \sum_{i=1}^\ell b_{i,-1} b^i_{-1} - \rho_{-2}.
\end{equation}
This completes the proof.
\end{proof}

\subsection{Quasi-conformal structures at the critical level}
\label{quasi-conf str}

Now we use this lemma to obtain additional information about the
homomorphism $w\cri$. Recall that a vertex algebra is called
quasi-conformal if it carries an action of the Lie algebra $\DerO$
satisfying the conditions of \cite{vertex}, \S~5.2.4. In particular, a
conformal vertex algebra is automatically quasi-conformal (with the
$\DerO$--action coming from the Virasoro action).

The Lie algebra $\DerO$ acts naturally on $\ghat\cri$ preserving
$\g[[t]]$, and hence it acts on $V\cri(\g)$. The
$\DerO$--action on $V\cri(\g)$ coincides with the limit $\ka \to
\ka_c$ of the $\DerO$--action on $V_\ka(\g), \ka \neq \ka_c$, obtained
from the Sugawara conformal structure.  Therefore this action defines
the structure of a quasi-conformal vertex algebra on $V\cri(\g)$.

Next, we define the structure of a quasi-conformal algebra on $M_\g
\otimes \pi_0$ as follows. The vertex algebra $M_\g$ is conformal with
the conformal vector \eqref{conf vec a}, and hence it is also
quasi-conformal. The commutative vertex algebra $\pi_0$ is the $\ka
\to \ka_c$ limit of the family of conformal vertex algebras
$\pi^{\ka-\ka_c}_0$ with the conformal vector \eqref{conf vec b}. The
induced action of the Lie algebra $\DerO$ on $\pi^{\ka-\ka_c}_0$ is
well-defined in the limit $\ka \to \ka_c$ and so it induces a
$\DerO$--action on $\pi_0$. Therefore it gives rise to the structure
of a quasi-conformal vertex algebra on $\pi_0$. The $\DerO$--action is
in fact given by derivations of the algebra structure on $\pi_0 \simeq
\C[b_{i,n}]$, and hence by Lemma 5.2.5 of \cite{vertex} it defines the
structure of a quasi-conformal vertex algebra on $\pi_0$. Explicitly,
the action of the basis elements $L_n = -t^{n+1}\pa_t, n\geq -1$, of
$\DerO$ on $\pi_0$ is determined by the following formulas:
\begin{align} \notag
L_n \cdot b_{i,m} &= - m b_{i,n+m}, \qquad -1\leq n<-m, \\
\label{action of Ln}
L_n \cdot b_{i,-n} &= n, \qquad n>0, \\ \notag
L_n \cdot b_{i,m} &= 0, \qquad n>-m
\end{align}
(note that $\rho(h_i)=1$ for all $i$). Now we obtain a quasi-conformal
vertex algebra structure on $M_\g \otimes \pi_0$ by taking the sum of
the above $\DerO$--actions.

Since the quasi-conformal structures on $V\cri(\g)$ and $M_\g \otimes
\pi_0$ both came as the limits of conformal structures as $\ka \to
\ka_c$, we obtain the following corollary of \lemref{Ska}:

\begin{cor}    \label{pres cri}
The homomorphism $w\cri: V\cri \to M_\g \otimes \pi_0$ preserves
quasi-conformal structures. In other words, it commutes with the
$\DerO$--action on both sides.
\end{cor}

\subsection{Transformation formulas for the fields}    \label{trans
  form for fields}

We can now obtain the transformation formulas for the fields
$a_\al(z), a^*_\al(z)$ and $b_i(z)$. According to the computations we
made in the proof of \lemref{Ska}, we have the following action of the
basis elements $L_n, n\geq 0$, of $\DerpO$ on the vectors
$a_{\al,-1}\vac$ and $a^*_{\al,0}\vac$:
$$
L_0 \cdot a_{\al,-1}\vac = a_{\al,-1}\vac, \qquad L_n \cdot
a_{\al,-1}\vac = 0, \quad n>0,
$$
$$
L_n \cdot a^*_{\al,0}\vac = 0, \quad n\geq 0.
$$
According to Proposition 5.3.8 of \cite{vertex}, this implies that
the field $a_\al(z) = Y(a_{\al,-1}\vac,z)$ transforms as a one-form on
the punctured disc $D^\times = \on{Spec} \C((z))$, while the field
$a^*_\al(z) = Y(a^*_{\al,0}\vac,z)$ transforms as a function on
$D^\times$. In particular, we obtain the following description of the
module $M_\g$. Consider the Heisenberg Lie algebra $\Gamma$ which is a
central extension of the commutative Lie algebra $U((t)) \oplus
U^*((t)) dt$ with the cocycle given by the formula
$$
f(t),g(t)dt \mapsto \int\lan f(t),g(t) \ran dt.
$$
This cocycle is coordinate--independent, and therefore $\Gamma$
carries natural actions of $\AutO$ and $\DerO$, which preserve the Lie
subalgebra $\Gamma_{+} = U[[t]] \oplus U[[t]] dt$. We identify the
completed Weyl algebra $\wt{\cA}^\g$ with a completion of
$U(\Gamma)/({\mb 1} - 1)$, where ${\mb 1} $ is the central
element. The module $M_\g$ is then identified with the
$\Gamma$--module induced from the one-dimensional representation of
$\Gamma_{+} \oplus \C {\mb 1}$, on which $\Gamma_{+}$ acts by $0$, and
${\mb 1}$ acts as the identity.  The $(\DerO,\AutO)$--action on $M_\g$
considered above is nothing but the natural action on the induced
module.

Now we consider the fields $b_i(z)$. We have
\begin{equation}
L_0 \cdot b_{i,-1}\vac = b_{i,-1}\vac, \qquad L_1 \cdot b_{i,-1}\vac =
2 \vac, \qquad L_n \cdot b_{i,-1}\vac = 0, \quad n>1.
\end{equation}
Recall that $b_i(z) = Y(b_{i,-1}\vac,z)$. According to \cite{vertex},
\S~7.1.11, these formulas imply that $\pa_z + b_i(z)$ transforms as a
connection on the line bundle $\Omega^{-\rho(h_i)}$ over $D^\times$
with values in $\on{End} \pi^{\ka-\ka_c}_0$.

Consider the $\h^*$--valued field ${\mb b}(z) = \sum_{i=1}^\ell b_i(z)
\omega_i$ such that $\lan {\mb b}(z),h_i \ran = b_i(z)$. The space
$\h^*$ is the Lie algebra of the group $^L H$ dual to $H$ (and it is
the Cartan subalgebra of the Langlands dual Lie algebra $^L\g$ of
$\g$). Denote by $\Omega^{-\rho}$ the unique principal $^L H$--bundle
on a smooth curve or a (punctured) disc such that the line bundle
associated any character $\check{\la}: {}^L H \to \C^\times$
(equivalently, a cocharacter of $H$) is $\Omega^{-\lan \rho,\check\la
\ran}$. The space of connections on this bundle is a torsor over $\h^*
\otimes \Omega$. The above statement about $b_i(z)$ may be
reformulated as follows: the operator $\pa_z + {\mb b}(z)$ transforms
as a connection on the $^L H$--bundle $\Omega^{-\rho}$.

Now we can answer the question posed in \secref{coord-indep version}.
Let $\K$ be a complete local field isomorphic to $\C((t))$, $\OO$ its
ring of integers isomorphic to $\C[[t]]$, and $\Au$ the corresponding
$\AutO$--torsor of formal coordinates on the disc $D = \on{Spec}
\OO$. Let $\ghat_\kappa(\K)$ be the affine algebra obtained as the
central extension of $\g \otimes \K$ (corresponding to the inner
product $\ka$ on $\g$). Let $\Gamma(\K)$ the Heisenberg Lie algebra
obtained as the central extension of $(U \otimes \K) \oplus (U^*
\otimes \Omega_\K)$, and $\Gamma_+(\K)$ its commutative Lie subalgebra
$(U \otimes \OO) \oplus (U^* \otimes \Omega_{\OO})$. Let $\V_\ka(\g)$
be the $\ghat_\kappa(\K)$--module defined by formula \eqref{mcV} and
${\mc M}_\g(\K)$ be the $\Gamma(\K)$--module
$$
{\mc M}_\g(\K) = \Au \underset{\AutO}\times M_\g.
$$
According to the above description of $M_\g$ as the induced module, we
have
\begin{align*}
  {\mc M}_\g(\K) &= \on{Ind}_{\Gamma_+(\K) \oplus \C{\mb 1}}^{\Gamma(\K)}
  \C_1 \\ & \simeq \on{Sym} (U \otimes \K \oplus U^* \otimes
  \Omega_\K)/(U \otimes \OO \oplus U^* \otimes \Omega_{\OO}) \\
  &\simeq \on{Fun}(U \otimes \Omega_{\OO} \oplus U^* \otimes
  \OO).
\end{align*}
The last isomorphism is obtained using the residue pairing.

Now consider the twist
$$
\Pi_0 = \Au \underset{\AutO}\times \pi_0,
$$ This is a module over the Lie algebra $$\hh(\K) = \Au
\underset{\AutO}\times \hh,$$ which is a central extension of $\h
\otimes \K$. According to the transformation formula for ${\mb b}(z)$
obtained above, $\hh(\K)$ may be described in the following
way. Consider the vector space
$\on{Conn}_{\{\la\}}(\Omega^{-\rho})_{D^\times}$ of $\la$--connections
on the $^L H$--bundle $\Omega^{-\rho}$ on $D^\times = \on{Spec} \K$
(for all possible complex values of $\la$). It fits into the exact
sequence
$$
0 \to \h^* \otimes \Omega_{D^\times} \to
\on{Conn}_{\{\la\}}(\Omega^{-\rho})_{D^\times} \to \C \pa \to 0.
$$
Then $\hh(\K)$ is by definition the topological dual vector space to
$\on{Conn}_{\{\la\}}(\Omega^{-\rho})_{D^\times}$ with the zero Lie
bracket. Let ${\mb 1}$ be the element dual to $\pa$. Any connection
$\nabla$ on $\Omega^{-\rho}$ over the punctured disc $\on{Spec} \K$
defines a one-dimensional representation $\C_\nabla$ of $\hh(\K)$
taking the value $1$ on ${\mb 1}$. If we choose an isomorphism $\K
\simeq \C((z))$, then the connection is given by the formula $\nabla =
\pa_z + \chi(z)$, where $\chi(z) \in \h^*((z))$. The action of the
generators $b_{i,n}$ is then given by the formula
$$
b_{i,n} \mapsto \int \lan \chi(z),h_i \ran z^n dz.
$$

Now $\Pi_0$ is identified as an $\hh(\K)$--module with the space of
functions on the subspace $\on{Conn}(\Omega^{-\rho})_{D}$ of
connections on the disc $D=\on{Spec} \OO$. The annihilator of $\Pi_0$
in $\hh(\K)$ is the subspace $\h \otimes \OO$, which is the orthogonal
complement of $\on{Conn}(\Omega^{-\rho})_{D}$.

\begin{prop}    \label{param of wak}
For any connection on the $^L H$--bundle $\Omega^{-\rho}$ over the
punctured disc $\on{Spec} \K$, there is a canonical
$\ghat\cri(\K)$--module structure on ${\mc M}_\g(\K)$.
\end{prop}

\begin{proof}
By \thmref{exist hom}, there exists a homomorphism of vertex algebras
$\ww_{\ka_c}: V\cri(\g) \to M_\g \otimes \pi_0$. According to
\corref{pres cri}, it commutes with the action of $\DerO$ and $\AutO$
on both sides. Therefore the corresponding homomorphism of Lie
algebras $\ghat \to U(M_\g \otimes \pi_0)$ also commutes with the
action of $\DerO$ and $\AutO$. Hence we may twist this homomorphism
with the $\AutO$--torsor $\Au$. Then we obtain a homomorphism of Lie
algebras
$$
\ghat_{\ka_c}(\K) \to U(M_\g \otimes \pi_0)(\K) = \Au
\underset{\AutO}\times U(M_\g \otimes \pi_0).
$$
Let us call a $\Gamma(\K) \oplus \hh(\K)$--module smooth if any vector
in this module is annihilated by the Lie subalgebra
$$(U \otimes \mm^N) \oplus (U^* \otimes \mm^N \Omega_{\OO}) \oplus (\h
\otimes \mm^N),$$ where $\mm$ is the maximal ideal of $\OO$, for
sufficiently large $N$. Then clearly any smooth $\Gamma(\K) \oplus
\hh(\K)$--module is automatically a $U(M_\g \otimes
\pi_0)(\K)$--module and hence a $\ghat_{\ka_c}(\K)$--module. Note that
${\mc M}_\g(\K)$ is a smooth $\Gamma(\K)$--module, and $\C_\nabla$ is
a smooth $\hh(\K)$--module for any connection $\nabla$ on the $^L
H$--bundle $\Omega^{-\rho}$ over the punctured disc $\on{Spec} \K$.
Taking the tensor product of these two modules we obtain a
$\ghat_{\ka_c}(\K)$--module, which is isomorphic to ${\mc M}_\g(\K)$
as a vector space. If we choose an isomorphism $\K \simeq \C((z))$,
this is nothing but the Wakimoto module $W_{\chi(z)}$ introduced in
\corref{wchit}.
\end{proof}

Thus, we obtain a family of $\ghat\cri(\K)$--modules parameterized by
flat connections on the $^L H$--bundle $\Omega^{-\rho}$ over the
  punctured disc $\on{Spec} \K$.

\section{Semi-infinite parabolic induction}    \label{parabolic}

\subsection{Wakimoto modules as induced representations}

The construction of the Wakimoto modules presented above may be
summarized as follows: for each representation $N$ of the Heisenberg
Lie algebra $\hh_\ka$, we have constructed a
$\ghat_{\ka+\ka_c}$--module structure on $M_\g \otimes N$. The
procedure consists of the extension the $\hh_\ka$--module by $0$ to
$\wh\bb_{-,\ka}$ followed by what may be viewed as a semi-infinite
analogue of induction from $\wh\bb_{-,\ka}$ to $\ghat_{\ka+\ka_c}$. An
important feature of this construction, as opposed to the ordinary
induction, is that the level gets shifted by $\ka_c$.  In particular,
if we start with an $\hh_0$--module, or equivalently, a representation
of the commutative Lie algebra $L\h$, then we obtain a
$\ghat_{\ka_c}$--module of critical level, rather than of level $0$.
Irreducible smooth representations of $L\h$ are one-dimensional and
are in one-to-one correspondence with the elements $\chi(t)$ of the
(topological) dual space $(L\h)^* \simeq \h^*((t)) dt$. Thus, we
obtain the Wakimoto modules $W_{\chi(t)}$ of critical level. But as we
have seen above, if we look at the transformation properties of
$\chi(t)$ under the action of the group $\AutO$ of changes of the
coordinate $t$, we find that $\chi(t)$ actually transforms not as a
one-form, but as a connection on a specific $^L H$--bundle.

In contrast, if $\ka \neq 0$, the irreducible smooth
$\hh_\ka$--modules are just the Fock representations $\pi^\ka_\la$. To
each of them we attach a $\ghat_{\ka+\ka_c}$--module
$W_{\la,\ka+\ka_c}$.

Now we want to generalize this construction replacing the Borel
subalgebra $\bb_-$ and its Levi quotient $\h$ by an arbitrary
parabolic subalgebra $\pp$ and its Levi quotient $\mm$. Then we wish
to attach to a module over a central extension of the loop algebra
$L\mm$ a $\ghat$--module. It turns out that this is indeed possible
provided that we pick a suitable central extension of $L\mm$. We call
the resulting $\ghat$--modules the generalized Wakimoto modules
corresponding to $\pp$. Thus, we obtain a functor from the category of
smooth $\wh{\mm}$--modules to the category of smooth
$\ghat$--modules. It is natural to call it the functor of {\em
semi-infinite parabolic induction} (by analogy with a similar
construction for representations of reductive groups).

\subsection{The main result}    \label{smooth}

Let $\pp$ be a parabolic Lie subalgebra of $\g$. We will assume that
$\pp$ contains the lower Borel subalgebra $\bb_-$ (and so in
particular, $\pp$ contains $\hh$). Let
$$
\pp = \mm \oplus \rr
$$
be the direct sum decomposition of $\pp$, where $\mm$ is the Levi
subgroup containing $\hh$ and $\rr$ is the nilpotent radical of
$\pp$. Further, let
$$
\mm = \bigoplus_{i=1}^s \mm_i \oplus \mm_0
$$
be the decomposition of $\mm$ into the direct sum of simple Lie
subalgebras $\mm_i, i=1,\ldots,s$, and an abelian subalgebra $\mm_0$
such that these direct summands are mutually orthogonal with respect
to the inner product on $\g$. We denote by $\ka_{i,c}$ the critical
inner product on $\mm_i, i=1,\ldots,s$, defined as in \secref{hom of
VAs}. We also set $\ka_{0,c}=0$.

Given a set of inner products $\ka_i$ on $\mm_i, 0=1,\ldots,s$, we
obtain an inner product on $\mm$. Let $\wh\mm_{(\ka_i)}$ be the
corresponding affine Kac-Moody algebra, i.e., the one-dimensional
central extension of $L\mm$ with the commutation relations given by
formula \eqref{KM rel}. We denote by $V_{\ka_i}(\mm_i), i=1,\ldots,s$,
the vacuum module over $\wh\mm_i$ with the vertex algebra defined as
in \secref{hom of VA}. We also denote by $V_{\ka_0}(\mm_0)$ the Fock
representation $\pi^{\ka_0}_0$ of the Heisenberg Lie algebra
$\mm_{\ka_0}$ with its vertex algebra structure defined as in
\secref{pi0}. Let
$$
V_{(\ka_i)}(\mm) \overset{\on{def}}{=} \bigotimes_{i=0}^s
V_{\ka_i}(\mm_i)
$$
be the vacuum module over $\wh\mm_{(\ka_i)}$ with the tensor product
vertex algebra structure.

Denote by $\De'_+$ the set of positive roots of $\g$ which do not
belong to $\pp$. Let $\cA^{\g,\pp}$ be the Weyl algebra with
generators $a_{\al,n}, a^*_{\al,n}, \al \in \De'_+, n \in \Z$,
and relations \eqref{commina}. Let $M_{\g,\pp}$ be the Fock
representation of $\cA^{\g,\pp}$ generated by a vector $\vac$ such
that
$$
a_{\al,n} \vac = 0, \quad n\geq 0; \qquad a^*_{\al,n} \vac = 0, \quad
n>0.
$$
Then $M_{\g,\pp}$ carries a vertex algebra structure defined as in
\secref{hva}.

We have the following analogue of \thmref{noncrit}.

\begin{thm}    \label{noncrit1}
Suppose that $\ka_i, i=0,\ldots,s$, is a set of inner products such
that there exists an inner product $\ka$ on $\g$ whose restriction to
$\mm_i$ equals $\ka_i-\ka_{i,c}$ for all $i=0,\ldots,s$. Then there
exists a homomorphism of vertex algebras $$w^{\pp}_\ka:
V_{\kappa+\kappa_c}(\g) \to M_{\g,\pp} \otimes V_{(\ka_i)}(\mm).$$
\end{thm}

\begin{proof}
The proof is a generalization of the proof of \thmref{noncrit} (in
fact, \thmref{noncrit} is a special case of \thmref{noncrit1} when
$\pp=\bb_-$). Let $P$ be the Lie subgroup of $G$ corresponding to
$\pp$, and consider the homogeneous space $G/P$. It has an open dense
subset $U_\pp = N_\pp \cdot [1]$, where $N_\pp$ is the subgroup of
$N_+$ corresponding to the subset $\De'_+ \subset \De_+$. We identify
$U_\pp$ with $N_\pp$ and with its Lie algebra $\n_\pp$ using the
exponential map.

Set $LU_\pp = U_\pp((t))$. We define functions and vector fields on
$LU_\pp$, denoted by $\on{Fun} LU_\pp$ and $\on{Vect} LU_\pp$,
respectively, in the same way as in \secref{infdim Weyl}. The action
of $L\g$ on $U_\pp((t))$ gives rise to a Lie algebra homomorphism
$$
\wh\rho_\pp: L\g \to \on{Vect} LU_\pp \oplus \on{Fun} U_\pp \wh\otimes
L\mm.
$$
Moreover, the image of this homomorphism is contained in the
``local part'', i.e., the direct sum of the local part
$\T^{\g,\pp}_{\on{loc}}$ of $\on{Vect} LU_\pp$ defined as in
\secref{local extension} and the local part ${\mc
  I}^{\g,\pp}_{\on{loc}}$ of $\on{Fun} U_\pp \wh\otimes L\mm$. By
definition, ${\mc I}^{\g,\pp}_{\on{loc}}$ is the span of the Fourier
coefficients of the formal power series $P(\pa_z^n a^*_\al(z))
J^b(z)$, where $P$ is a differential polynomial in $a^*_\al(z), \al
\in \De'_+$, and $J^a \in \mm$.

Let $\cA^{\g,\pp}_{0,\on{loc}}$ and $\cA^\g_{\leq 1,\on{loc}}$ be the
zeroth and the first terms of the natural filtration on the local
completion of the Weyl algebra $\cA^{\g,\pp}$, defined as in
\secref{local filtr}. We have a non-split exact sequence
\begin{equation}    \label{local ext1}
0 \to \cA^{\g,\pp}_{0,\on{loc}} \to \cA^{\g,\pp}_{\leq 1,\on{loc}} \to
\T^{\g,\pp}_{\on{loc}} \to 0.
\end{equation}
Set
$$
{\mc J}^{\g,\pp}_{\on{loc}} \overset{\on{def}}{=} \cA^{\g,\pp}_{\leq
1,\on{loc}} \oplus {\mc I}^{\g,\pp}_{\on{loc}},
$$
and note that ${\mc J}^{\g,\pp}_{\on{loc}}$ is naturally a Lie
subalgebra of the local Lie algebra $U(M_{\g,\pp} \otimes
V_{(\ka_i)}(\mm))$. Using the splitting of the sequence \eqref{local
ext1} as a vector space via the normal ordering, we obtain a linear
map $\ol{w}_{(\ka_i)}: L\g \to {\mc J}^{\g,\pp}_{\on{loc}}$.

We need to compute the failure of $\ol{w}_{(\ka_i)}$ to be a Lie
algebra homomorphism. Thus, we consider the corresponding linear map
$\omega_{(\ka_i)}: \bigwedge^2 L\g \to {\mc J}^{\g,\pp}_{\on{loc}}$
defined by the formula
$$
\omega_{(\ka_i)}(f,g) = [\ol{w}_\ka(f),\ol{w}_\ka(g)] -
\ol{w}_\ka([f,g]).
$$
Evaluating it explicitly in the same way as in the proof of
\lemref{repr in loc}, we find that $\omega_{(\ka_i)}$ takes values in
$\cA^{\g,\pp}_{0,\on{loc}} \subset {\mc
J}^{\g,\pp}_{\on{loc}}$. Furthermore, $\cA^{\g,\pp}_{0,\on{loc}}$ is
naturally an $L\g$--module, and by construction of $\ol{w}_{(\ka_i)}$,
for any $X \in \cA^{\g,\pp}_{0,\on{loc}}$ and $f \in L\g$ we have
$[\ol{w}_{(\ka_i)}(f),X] = f \cdot X$.  This implies that
$\omega_{(\ka_i)}$ is a two-cocycle of $L\g$ with coefficients in
$\cA^\g_{0,\on{loc}}$. By construction, it is local, i.e., belongs to
$C^2_{\on{loc}}(L\g,\cA^{\g,\pp}_{0,\on{loc}})$.

Following the argument used in the proof of \lemref{any two}, we show
that any two cocycles in
$C^2_{\on{loc}}(L\g,\cA^{\g,\pp}_{0,\on{loc}})$, whose restrictions to
$\bigwedge{}^2(L\mm)$ coincide, represent the same cohomology class.

Let us compute the restriction of $\omega_{(\ka_i)}$ to $\bigwedge^2
L\mm$. For that we evaluate $\ol{w}_{(\ka_i)}$ on elements of $L\mm$.
Let $\{ J^i \}_{i=1,\ldots,\dim \mm}$ be a basis of $\mm$. The adjoint
action of the Lie algebra $\mm$ on $\g$ preserves $\n_\pp$. Under this
action, which we denote by $\rho_{\n_\pp}$, an element $A \in \mm$
acts by the formula
$$
\rho_{\n_\pp}(A) \cdot J^\al = \sum_{\beta \in \De'_+}
c^{\al}_\beta(A) J^\beta
$$
for some $c^{\al}_\beta(A) \in \C$. Therefore
$$
\ol{w}_{(\ka_i)}(A(z)) = - \sum_{\beta \in \De'_+} c^{\al}_\beta(A)
\Wick a^*_\beta(z) a_\beta(z) \Wick + \wt{A}(z),
$$
where $\wt{A}(z)$ is the field $\sum_{n \in \Z} (A \otimes t^n)
z^{-n-1}$, considered as a generating series of elements of ${\mc
I}^{\g,\pp}_{\on{loc}}$. Let $\ka_{\n_\pp}$ be the inner product on
$\mm$ defined by the formula
$$
\ka_{\n_\pp}(A,B) = \on{Tr}_{\n_\pp} \rho_{\n_\pp}(A)
\rho_{\n_\pp}(B).
$$
Moreover, we find that for $A \in \mm_i, B \in \mm_j$,
\begin{equation}    \label{formula1}
\omega_{(\ka_i)}(A_n,B_m) = n(- \ka_{\n_\pp}(A,B) + \ka_i(A,B))
\delta_{n,-m},
\end{equation}
if $i=j$, and
\begin{equation}    \label{formula2}
\omega_{(\ka_i)}(A_n,B_m) = - n \ka_{\n_\pp}(A,B) \delta_{n,-m},
\end{equation}
if $i \neq j$. Thus, the restriction of $\omega_{(\ka_i)}$ to
$\bigwedge^2 L\mm$ takes values in the constants $\C \subset
\cA^{\g,\pp}_{0,\on{loc}}$.

Let $\ka_{_K}$ be the Killing form on $\g$ and $\ka_{i,K}$ be the
Killing form on $\mm_i$ (in particular, $\ka_{0,K} = 0$). Then
\begin{align*}
\ka_{_K}(A,B) &= \ka_{i,K}(A,B) + 2 \ka_{\n_\pp}(A,B), \qquad \on{if}
\quad i=j, \\
\ka_{_K}(A,B) &= 2 \ka_{\n_\pp}(A,B), \qquad \on{if} \quad i \neq j.
\end{align*}
Since $\mm_i$ and $\mm_j$ are orthogonal with respect to $\ka_{_K}$
for all $i \neq j$ by construction, we obtain that $\ka_{\n_\pp}(A,B)
= 0$, if $i \neq j$. Recall that by definition $\ka_c = -
\frac{1}{2} \ka_K$, and $\ka_{i,c} = -
\frac{1}{2} \ka_{i,K}$. Hence $\ka_{\n_\pp}|_{\mm_i} = - \ka_{c} +
\ka_{i,c}$, and so if $\ka_{i}$ equals $\ka|_{\mm_i} + \ka_{i,c}$ for
some invariant inner product $\ka$ on $\g$, then
$-\ka_{\n_\pp}|_{\mm_i} + \ka_i = (\ka+\ka_c)|_{\mm_i}$.

Thus, we obtain that if $\ka$ is an invariant innder product on $\g$
whose restriction to $\mm_i$ equals $\ka_i-\ka_{i,c}$ for all
$i=0,\ldots,s$, then the restriction of the two-cocycle
$\omega_{(\ka_i)}$ to $\bigwedge{}^2(L\mm)$ is equal to the
restriction to $\bigwedge{}^2(L\mm)$ of the two-cocycle
$\sigma_{\ka+\ka_c}$ on $L\g$ (corresponding to the one-dimensional
central extension with respect to the inner product $\ka+\ka_c$ on
$\g$). Therefore the $\omega_{(\ka_i)}$ is equivalent to
$\sigma_{\ka+\ka_c}$ in this case and we obtain that the linear map
$\ol{w}_{(\ka_i)}: L\g \to {\mc J}^{\g,\pp}_{\on{loc}}$ may be
modified by the addition of an element of
$C^1_{\on{loc}}(L\g,\cA^{\g,\pp}_{0,\on{loc}})$ to give us a Lie
algebra homomorphism
$$
\ghat_{\ka+\ka_c} \to {\mc J}^{\g,\pp}_{\on{loc}} \subset U(M_{\g,\pp}
\otimes V_{(\ka_i)}(\mm)).
$$
Now \lemref{hom from V} implies that there exists a homomorphism of
vertex algebras
$$
w^{\pp}_\ka: V_{\kappa+\kappa_c}(\g) \to M_{\g,\pp} \otimes
V_{(\ka_i)}(\mm).
$$
This completes the proof.
\end{proof}

We will call an $\wh\mm_{(\ka_i)}$--module {\em smooth} if any vector
in it is annihilated by the Lie subalgebra $\mm \otimes t^N\C[[t]]$
for sufficiently large $N$.

\begin{cor}    \label{functorial}
  For any smooth $\wh\mm_{(\ka_i)}$--module $R$ with the
  $\ka_i$'s satisfying the conditions of \thmref{noncrit1}, the tensor
  product $M_{\g,\pp} \otimes R$ is naturally a smooth
  $\ghat_{\ka+\ka_c}$--module.  There is a functor from the category
  of smooth $\wh\mm_{(\ka_i)}$--modules to the category of smooth
  $\ghat_{\ka+\ka_c}$--module sending a module $R$ to $M_{\g,\pp}
  \otimes R$ and $\wh\mm_{(\ka_i)}$--homomorphism $R_1 \to R_2$ to the
  $\ghat_{\ka+\ka_c}$--homomorphism $M_{\g,\pp} \otimes R_1 \to
  M_{\g,\pp} \otimes R_2$.
\end{cor}

We call the $\ghat_{\ka+\ka_c}$--module $M_{\g,\pp} \otimes R$ the
{\em generalized Wakimoto module corresponding to} $R$.

Consider the special case when $R$ is the tensor product of the
Wakimoto modules $W_{\la_i,\ka_i}$ over $\wh{\mm}_i, i=1,\ldots,s$,
and the Fock representation $\pi^{\ka_0}_{\la_0}$ over the Heisenberg
Lie algebra $\wh{\mm}_0$. In this case it follows from the
construction that the corresponding $\ghat_{\ka+\ka_c}$--module
$M_{\g,\pp} \otimes R$ is isomorphic to the Wakimoto module
$W_{\la,\ka+\ka_c}$ over $\ghat_{\ka+\ka_c}$, where $\la = (\la_i)$.

\subsection{General parabolic subalgebras}

So far we have worked under the assumption that the parabolic
subalgebra $\pp$ contains $\bb_-$. It is also possible to construct
Wakimoto modules associated to other parabolic subalgebras. Let us
explain how to do this in the case when $\pp = \bb_+$. We have the
involution of $\g$ sending $e_i$ to $f_i$ and $h_i$ to $-h_i$. Under
this involution $\bb_-$ goes to $\bb_+$.

Let $N$ be any module over the vertex algebra $M_\g \otimes
\pi^{\kappa-\kappa_c}_0$. Then \thmref{noncrit} implies that the
following formulas define a $\ghat_\ka$--structure on $N$ (with $K$
acting as the identity):
\begin{align*} \notag
f_i(z) &\mapsto a_{\al_i}(z) + \sum_{\beta \in \Delta_+} \Wick
P^i_\beta(a^*_\al(z)) a_\beta(z) \Wick \; , \\
h_i(z) &\mapsto \sum_{\beta \in \Delta_+} \beta(h_i) \Wick
a^*_\beta(z) a_\beta(z) \Wick - b_i(z), \\ \notag f_i(z) &\mapsto
\sum_{\beta \in \Delta_+} \Wick Q^i_\beta(a^*_\al(z)) a_\beta(z) \Wick
+ \left( c_i+(\kappa-\kappa_c)(e_i,f_i) \right) \pa_z a^*_{\al_i}(z) +
b_i(z) a^*_{\al_i}(z),
\end{align*}
where the polynomials $P^i_\beta, Q^i_\beta$ are introduced in
formulas \eqref{formulas1}--\eqref{formulas3}.

In order to make $N$ into a module with highest weight, we choose
$N$ as follows. Let $M'_{\g}$ be the Fock representation of the Weyl
algebra $\cA^{\g}$ generated by a vector $\vac'$ such that
$$
a_{\al,n} \vac' = 0, \quad n> 0; \qquad a^*_{\al,n} \vac' = 0, \quad
n \geq 0.
$$
We take as $N$ the module $M'_\g \otimes
\pi^{\ka-\ka_c}_{-2\rho-\la}$, where $\pi^{\ka-\ka_c}_{-2\rho-\la}$ is
the $\pi^{\ka-\ka_c}_0$--module defined in \secref{pi00}. It is easy to
see that the corresponding $\ghat_\ka$--module has a highest weight
vector on which $\h \otimes 1$ acts through the weight $\la$. We
denote this module by $W^+_{\la,\ka}$. This is the Wakimoto module
corresponding to $\bb_+$.

The following result will be used in \secref{char of zz}.

\begin{prop}    \label{Wak and Verma}
The Wakimoto module $W^+_{0,\ka_c}$ is isomorphic to the Verma module
$M_{0,\ka_c}$.
\end{prop}

\begin{proof}
The proof is analogous to the proof of \thmref{KK conj}. Since
$W^+_{0,\ka_c}$ has a highest weight vector of the same weight as that
of the highest weight vector of $M_{0,\ka_c}$, there is a non-zero
homomorphism $M_{0,\ka_c} \to W^+_{0,\ka_c}$. But the characters of
$M_{0,\ka_c}$ and $W^+_{0,\ka_c}$ are equal. Therefore the theorem
will follow if we show that $W^+_{0,\ka_c}$ is generated by its
highest weight vector.

Suppose that $W^+_{0,\ka_c}$ is not generated by its highest weight
vector. Then there exists a homogeneous linear functional on
$W^+_{0,\ka_c}$, whose weight is less than the highest weight
$(\la,0)$ and which is invariant under the Lie subalgebra
$$
\wt\n_- = (\n_- \otimes 1) \oplus (\g \otimes t^{-1} \C[t^{-1}]),
$$ and in particular, under its Lie subalgebra $L_- \n_- = \n_-
\otimes \C[t^{-1}]$. Therefore this functional factors through the
space of coinvariants of $W_{0,\ka_c}$ by $L_- \n_-$. But it follows
from the construction that $L_- \n_-$ acts freely on $W^+_{0,\ka_c}$,
and the space of coinvariants is isomorphic to the subspace
$\C[a_{\al,n}^*]_{\al\in\De_+,n<0} \otimes
\C[b_{i,n}]_{i=1,\ldots,\ell;n<0}$ of $W^+_{0,\ka_c}$. Hence we obtain
that the weight of this functional has the form
\begin{equation}    \label{form of weight}
- \sum_j (n_j\delta - \beta_j), \qquad n_j > 0, \quad \beta_j \in
\De_+.
\end{equation}
But then $M_{0,\ka_c}$ must have an irreducible subquotient of highest
weight of this form.

Now recall the Kac--Kazhdan theorem \cite{KK} describing the set of
highest weights of irreducible subquotients of Verma modules (see
\secref{KK}). In the case at hand the statement is as follows. A
weight $\wh\mu = (\mu,n)$ appears in the decomposition of
$M_{0,\ka_c}$ if and only $n\leq 0$ and either $\mu=0$ or there
exists a finite sequence of weights $\mu_0,\ldots,\mu_m \in \h^*$
such that $\mu_0 = \mu, \mu_m = 0$, $\mu_{i+1} = \mu_i \pm m_i
\beta_i$ for some positive roots $\beta_i$ and positive integers
$m_i$ which satisfy
\begin{equation}    \label{sing vec}
2(\mu_i+\rho,\beta_i) = m_i (\beta_i,\beta_i)
\end{equation}
(here $(\cdot,\cdot)$ is the inner product on $\h^*$ induced by an
arbitrary non-degenerate invariant inner product on $\g$).

Now observe that the equations \eqref{sing vec} coincide with the
equations appearing in the analysis of irreducible subquotients of the
Verma module over $\g$ of highest weight $0$. This implies that the
above statement is equivalent to the following: a weight $\wh\mu =
(\mu,n)$ appears in the decomposition of $M_{0,\ka_c}$ if and only
$n\leq 0$ and $\mu = w(\rho)-\rho$ for some element $w$ of the Weyl
group of $\g$. But for any $w$, the weight $w(\rho)-\rho$ equals the
sum of negative simple roots of $\g$. Hence the weight of any
irreducible subquotient of $M_{0,\ka_c}$ has the form $-n \delta -
\sum_i m_i \al_i, m\geq 0$. Such a weight cannot be of the form
\eqref{form of weight}. Therefore $W^+_{0,\ka_c}$ is generated by the
highest weight vector and hence is isomorphic to $M_{0,\ka_c}$.
\end{proof}

\begin{remark}
The same argument as in the proof of \propref{Wak and Verma} shows
that the Wakimoto module $W^+_{\la,\ka_c}$ is isomorphic to
$M_{\la,\ka_c}$ if $\la$ is such that all weights $(\mu,n)$ of
irreducible subquotients of $M_{\la,\ka_c}$ satisfy $\mu = \la -
\sum_i m_i \al_i, m_i\geq 0$. Likewise, the Wakimoto module
$W_{\la,\ka_c}$ is isomorphic to $M_{\la,\ka_c}$ if we have $\mu = \la
+ \sum_i m_i \al_i, m_i\geq 0$ for all such $\mu$.\qed
\end{remark}

In \secref{char of zz} we will need one more result on the structure
of $W^+_{0,\ka_c}$. Consider the Lie algebra $\wt{\bb}_+ = (\bb_+
\otimes 1) \oplus (\g \otimes t\C[[t]])$.

\begin{lem}    \label{sing vect imag}
The space of $\wt{\bb}_+$--invariants of $W^+_{0,\ka_c}$ is equal to
$\pi_0 \subset W^+_{0,\ka_c}$.
\end{lem}

\begin{proof}
It follows from the above formulas for the action of $\ghat_{\ka_c}$
on $W^+_{0,\ka_c}$ that all vectors in $\pi_0$ are annihilated by
$\wt{\bb}_+$. Let us show that there are no other
$\wt{\bb}_+$--invariant vector in $W^+_{0,\ka_c}$.

A $\wt{\bb}_+$--invariant vector is in particular annihilated by the
Lie subalgebra $L_+\n_- = \n_+ \otimes t\C[[t]]$. In the same way as
in the proof of \thmref{KK conj} we show that the space of
$L_+\n_+$--invariants of $W_{\ka_c,0}$ is equal to the tensor product
of $\pi_0$ and the subspace of $W^+_{0,\ka_c}$ spanned by all
monomials of the form $\ds \prod_{m_{s(\al)}\leq 0}
a^*_{\al,m_{s(\al)}} \vac'$. But a $\wt{\bb}_+$--invariant vector is
also annihilated by $\h \otimes 1$. A monomials of the above form is
annihilated by $\h \otimes 1$ only if it belongs to $\pi_0$. Hence the
space of $\wt{\bb}_+$--invariants of $W^+_{0,\ka_c}$ is equal to
$\pi_0$.
\end{proof}

\section{Wakimoto modules over $\su$}

In this section we describe explicitly the Wakimoto modules over $\su$
and some intertwining operators between them.

Let $\{ e,h,f \}$ be the standard basis of the Lie algebra
$\sw_2$. Let $\ka_0$ be the invariant inner product on $\sw_2$
normalized in such a way that $\ka_0(h,h)=2$. We will write an
arbitrary invariant inner product $\ka$ on $\sw_2$ as $k\ka_0$, where
$k \in \C$, and will use $k$ in place of $\ka$ in our notation. In
particular, $\ka_c$ corresponds to $k=-2$. The set $\De_+$ consists of
one element in the case of $\sw_2$, so we will drop the index $\al$ in
$a_\al(z)$ and $a^*_\al(z)$. Likewise, we will drop the index $i$ in
$b_i(z)$, etc., and will write $M$ for $M_{\sw_2}$ throughout this
section. We will also identify the dual space to the Cartan subalgebra
$\h^*$ with $\C$ by sending $\chi \in \h^*$ to $\chi(h)$.

\subsection{The homomorphism $w_k$}

The Weyl algebra ${\mc A}_{\sw_2}$ has generators $a_n, a^*_n, n \in
\Z$, with the commutation relations
$$
[a_n,a^*_m] = \delta_{n,-m}.
$$
Its Fock representation is denoted by $M$. The Heisenberg Lie
algebra $\hh_k$ has generators $b_n, n \in \Z$, with the commutation
relations
$$
[b_n,b_m] = 2 k n K \delta_{n,-m},
$$
and $\pi^k_0$ is its Fock representation generated by a vector $\vac$
such that
$$
b_n \vac = 0, \quad n \geq 0; \qquad K\vac = \vac.
$$
The tensor product $M \otimes \pi^k_0$ is a vertex algebra.

The homomorphism $w_k: V_k(\sw_2) \to M \otimes \pi^{k+2}_0$ of vertex
algebras is given by the following formulas:
\begin{eqnarray}
e(z) &\mapsto& a(z) \notag \\ h(z) &\mapsto& -2\Wick a^*(z) a(z)\Wick
+b(z) \label{Wakimoto formulas} \\ f(z) &\mapsto& -\Wick a^*(z)^2
a(z)\Wick +k\delz a^*(z) +a^*(z)b(z). \notag
\end{eqnarray}

The Wakimoto module $M \otimes \pi^{k+2}_\la$ will be denoted by
$W_{\la,k}$ and its highest weight vector will be denoted by
$v_{\la,k}$.

Recall that under the homomorphism $w_k$ the generators $e_n$ of $\su$
are mapped to $a_n$ (from now on, by abuse of notation, we will
identify the elements of $\su$ with their images under $w_k$). The
commutation relations of $\su$ imply the following formulas
$$[e_n,a_{-1}] = 0, \qquad [h_n,a_{-1}] = 2 a_{n-1},
\quad [f_n,a_{-1}] = - h_{n-1} + k \delta_{n,1}.$$ In addition,
it follows from formulas \eqref{Wakimoto formulas} that $$e_n v_{-2,k}
= a_n v_{-2,k} = 0, \quad n\geq 0; \quad \quad h_n v_{-2,k} = f_n
v_{-2,k} = 0, \quad n>0,$$ $$h_0 v_{-2,k} = -2 v_{-2,k}.$$ Therefore $$e_n
\cdot a_{-1} v_{-2,k} = h_n \cdot a_{-1} v_{-2,k} = 0, \quad
\quad n\geq 0,$$ and $$f_n \cdot a_{-1} v_{-2,k} = 0, \quad \quad
n>1.$$ We also find that $$f_1 \cdot a_{-1} v_{-2,k} = (k+2)
v_{-2,k},$$ and
\begin{align*}
f_0 \cdot a_{-1} v_{-2,k} &= a_{-1} f_0 v_{-2,k} - h_0
v_{-2,k} \\ &= - b_{-1} v_{-2,k} =
(k+2) T v_{-2,k},
\end{align*}
where $T$ is the translation operator.

\subsection{Vertex operators associated to a module over a vertex
  algebra}    \label{generalities}

We need to recall some general results on the vertex operators
associated to a module over a vertex algebra, following \cite{FHL},
\S~5.1. Let $V$ be a conformal vertex algebra $V$ and $M$ a
$V$--module. Then there exists a linear map
$$
Y_{V,M}: M \to \on{Hom}(V,M)[[z^{\pm 1}]],
$$
satisfying the following property: for any $A \in V, B \in M, C \in
V$, there exists an element $$f \in
M[[z,w]][z^{-1},w^{-1},(z-w)^{-1}]$$ such that the formal power series
$$Y_M(A,z) Y_{V,M}(B,w) C, \qquad Y_{V,M}(B,w) Y(A,z) C,$$
$$Y_{V,M}(Y_{V,M}(B,w-z)A,z)C, \qquad Y_{V,M}(Y(A,z-w) B,w) C$$ are
expansions of $f$ in $$M((z))((w)), \qquad M((w))((z)), \qquad
M((z))((z-w)), \qquad M((w))((z-w)),$$ respectively (compare with
Corollary 3.2.3 of \cite{vertex}). Abusing notation, we will write
$$
Y_M(A,z) Y_{V,M}(B,w) = Y_{V,M}(Y(A,z-w) B,w),
$$
and call this formula the operator product expansion (OPE), as in the
case $M=V$ when $Y_{V,M} = Y$ (see \cite{vertex}, \S~3.3).

This property implies the following commutation relations between the
Fourier coefficients of $Y(A,z)$ and $Y_{V,M}(B,w)$, which is proved
in the same way as in the case $M=V$ (see \cite{vertex}, \S~3.3.6). If
we write
$$
Y_M(A,z) = \sum_{n \in \Z} A_{(n)} z^{-n-1}, \qquad Y_{V,M}(B,w) =
\sum_{n \in \Z} B_{(n)} w^{-n-1},
$$
then we have
\begin{equation}    \label{OPE}
[B_{(m)},A_{(k)}]= \sum_{n\geq 0} \left(
                        \begin{array}{c} m \\ n \end{array} \right)
                        (B_{(n)} \cdot A)_{(m+k-n)}.
\end{equation}

In particular, we obtain that
\begin{equation}    \label{kern}
\left[ \int Y_{V,M}(B,z) dz,Y(A,w) \right] = Y_{V,M}\left( \int
  Y_{V,M}(B,z) dz \cdot A,w \right). 
\end{equation}

\subsection{The screening operator}    \label{first screening}

Let us apply the results of \secref{generalities} in the case of the
vertex algebra $W_{0,k}$ and its module $W_{-2,k}$ for $k \neq -2$.
Then we find, in the same way as in \cite{vertex}, \S~4.2.6, that
$$
Y_{W_{0,k},W_{-2,k}}(a_{-1} v_{-2,k}) = S_k(z)
\overset{\on{def}}{=} a(z) V_{-2}(z),
$$
where $V_{-2}(z)$ is the bosonic vertex operator acting from
$\pi^{k+2}_0$ to $\pi^{k+2}_{-2}$ given by the formula
\begin{equation}
V_{-2}(z) = T_{-2} \exp \left( \frac{1}{k+2} \sum_{n<0}
\frac{b_n}{n} z^{-n} \right) \exp \left( \frac{1}{k+2} \sum_{n>0}
\frac{b_n}{n} z^{-n} \right).
\end{equation}
Here we denote by $T_{-2}$ the translation operator $\pi^{k+2}_0 \to
\pi^{k+2}_{-2}$ sending the highest weight vector to the highest
weight vector and commuting with all $b_n, n \neq 0$.

Now formula \eqref{OPE} implies that the operator $S_k(z)$ has the
following OPEs:
\begin{align*}
e(z) S_k(w) &= \on{reg.}, \qquad h(z) S_k(w) = \on{reg.}, \\
f(z) S_k(w) &= \frac{(k+2) V_{-2}(w)}{(z-w)^2} + \frac{(k+2) \pa_w
V_{-2}(w)}{z-w} + \on{reg.} \\ &= (k+2) \pa_w
\frac{V_{-2}(w)}{z-w} + \on{reg.}
\end{align*}
Since the residue of a total derivative is equal to $0$, we obtain the
following:

\begin{prop}    \label{screening first}
The residue $S_k = \ds \int S_k(w) dw$ is an intertwining
operator between the $\su$--modules $W_{0,k}$ and $W_{-2,k}$.
\end{prop}

We call $S_k$ the {\em screening operator of the first kind} for $\su$.

\begin{prop}    \label{seq is exact}
For $k \not\in -2+{\mathbb Q}_{\geq 0}$ the sequence
$$
0 \to V_k(\sw_2) \to W_{0,k} \overset{S_k}{\to} W_{-2,k} \to 0
$$
is exact.
\end{prop}

\begin{proof}
By \propref{injective}, $V_k(\sw_2)$ is naturally an $\su$--submodule
of $W_{0,k}$ for any value of $k$. The module $V_k(\sw_2)$ is
generated from the vector $v_k$, whose image in $W_{0,k}$ is the
highest weight vector. This vector is clearly in the kernel of
$S_k$. Hence $V_k(\sw_2)$ belongs to the kernel of $S_k$.

In order to prove that $V_k(\sw_2)$ coincides with the kernel of
$S_k$, we compare their characters. See \secref{KK} for the definition
of the character. We will use the notation $q=e^{-\delta}, u=e^{\al}$.

Since $V_k(\g)$ is isomorphic to the universal enveloping algebra of
the Lie algebra spanned by $e_n, h_n$ and $f_n$ with $n<0$, we find
that
\begin{equation}    \label{char vac}
\on{ch} V_k(\g) = \prod_{n>0}
(1-q^n)^{-1}(1-uq^n)^{-1}(1-u^{-1}q^n)^{-1}.
\end{equation}
Similarly, we obtain that
$$
\on{ch} W_{\la,k} = u^{\la/2} \prod_{n>0}
(1-q^n)^{-1}(1-uq^n)^{-1}(1-u^{-1}q^{n-1})^{-1}.
$$
Thus, $\on{ch} W_{\la,k} = \on{ch} M_{\la,k}$, where $M_{\la,k}$ is
the Verma module over $\su$ with highest weight $(\la,k)$. Therefore
if $M_{\la,k}$ is irreducible, then so is $W_{\la,k}$. The set of
values $(\la,k)$ for which $M_{\la,k}$ is irreducible is described in
\cite{KK}. It follows from this description that if $k \not\in
-2+{\mathbb Q}_{\geq 0}$, then $M_{-2,k}$, and hence $W_{-2,k}$, is
irreducible. It is easy to check that $S_k(a^*_0 v_{0,k}) = v_{-2,k}$,
so that $S_k$ is a non-zero homomorphism. Hence it is surjective for
such values of $k$. Therefore the character of its kernel for such
$k$ is equal to $\on{ch} W_{0,k} - \on{ch} W_{-2,k} = \on{ch}
V_k(\sw_2)$. This completes the proof.
\end{proof}

Next, we will describe the second screening operator for $\su$. For
that we need to recall the Friedan--Martinec--Shenker bosonization.

\subsection{Friedan--Martinec--Shenker bosonization}
\label{bosonization}

Consider the Heisenberg Lie algebra with the generators $p_n, q_n, n \in
\Z$, and the central element $K$ with the commutation relations
$$
[p_n,p_m] = n \delta_{n,-m} K, \qquad [q_n,q_m] = - n \delta_{n,-m} K,
\qquad [p_n,q_m] = 0.
$$
We set
$$
p(z) = \sum_{n \in \Z} p_n z^{-n-1}, \qquad q(z) = \sum_{n \in \Z} q_n
z^{-n-1}.
$$

For $\la \in \C, \mu \in \C$, let $\Pi_{\la,\mu}$ be the Fock
representation of this Lie algebra generated by a highest weight
vector $|\la,\mu\ket$ such that
$$
p_n |\la,\mu\ket = \la \delta_{n,0} |\la,\mu\ket, \quad q_n
|\la,\mu\ket = \mu \delta_{n,0} |\la,\mu\ket, \quad n \geq 0; \quad K
|\la,\mu\ket = |\la,\mu\ket.
$$

Consider the vertex operators $V_{\la,\mu}(z): \Pi_{\la',\mu'} \to
\Pi_{\la+\la',\mu+\mu'}$ given by the formula
$$
V_{\la,\mu}(z) = T_{\la,\mu} z^{\la\la' - \mu\mu'} \exp \left( -
\sum_{n<0} \frac{\la p_n + \mu q_n}{n} z^{-n} \right) \exp \left( -
\sum_{n>0} \frac{\la p_n + \mu q_n}{n} z^{-n} \right),
$$
where $T_{\la,\mu}$ is the translation operator $\Pi_{0,0} \to
\Pi_{\la,\mu}$ sending the highest weight vector to the highest
weight vector and commuting with all $p_n, q_n, n \neq 0$.

Abusing notation, we will write these operators as $e^{\la u + \mu
v}$, where $u(z)$ and $v(z)$ stand for the anti-derivatives
of $p(z)$ and $q(z)$, respectively, i.e., $p(z) = \pa_z u(z), q(z) =
\pa_z v(z)$.

For $\al \in \C$, set
$$
\Pi_{\al} = \bigoplus_{n \in \Z} \Pi_{n+\al,n+\al}.
$$
Using the vertex operators, one defines a vertex algebra structure on
the direct sum $\Pi_0$ (with the vacuum vector $|0,0\ket$) as in
\cite{vertex}, \S~4.2.6.

The following theorem is due to Friedan, Martinec and Shenker
\cite{FMS} (see also \cite{FF:weil}).

\begin{thm}
There is a (unique) embedding of vertex algebras $M \hookrightarrow
\Pi_0$ under which the fields $a(z)$ and $a^*(z)$ are mapped to the
fields
$$
\wt{a}(z) = e^{u+v}, \qquad \wt{a}^*(z) = (\pa_z e^{-u}) e^{-v} = -
\Wick p(z) e^{-u-v} \Wick \; .
$$
Further, the image of $M$ in $\Pi_0$ is equal to the kernel of the
operator $\ds \int e^u dz$.
\end{thm}

Equivalently, $\Pi_0$ may be described as the localization of $M$ with
respect to $a_{-1}$, i.e.,
\begin{equation} \label{localiz} \Pi_0
\simeq M[(a_{-1})^{-1}] = \C[a_n]_{n<-1} \otimes \C[a^*_n]_{n\leq 0}
\otimes \C[(a_{-1})^{\pm 1}].
\end{equation}
The vertex algebra structure on $\Pi_0$ is obtained by a natural
extension of the vertex algebra structure on $M$.

\medskip

\begin{remark}
Under the embedding $M \hookrightarrow \Pi_0$ the Virasoro vertex
operator of $M$ described in \secref{conformal str} maps to the
following conformal vector in $\Pi_0$:
$$
\frac{1}{2} \Wick p(z)^2 \Wick - \frac{1}{2} \pa_z p(z) - \frac{1}{2}
\Wick q(z)^2 \Wick + \frac{1}{2} \pa_z q(z).
$$
Thus, the map $M \hookrightarrow \Pi_0$ becomes a homomorphism of
conformal vertex algebras if we choose these conformal structures.\qed
\end{remark}

\medskip

%From this point of view, $\Pi_0$ may be identified with the vertex
%algebra of chiral differential operators corresponding to the variety
%$\C^\times = \on{Spec} \C[y,y^{-1}]$, while $M$ is its vertex
%subalgebra corresponding to $\C = \on{Spec} \C[y]$.

The above ``bosonization'' allows us to make sense of the field
$a(z)^\al$, where $\al$ is an arbitrary complex number. Namely, we
replace $a(z)^\al$ with the field
$$
\wt{a}(z)^\al = e^{\al(u+v)}: \Pi_0 \to \Pi_\al.
$$

Now we take the tensor product $\Pi_0 \otimes \pi^{k+2}_0$, where we
again assume that $k \neq -2$. This is a vertex algebra which contains
$M \otimes \pi^{k+2}_0$, and hence $V_k(\sw_2)$, as vertex
subalgebras. In particular, for any $\al, \la \in \C$, the tensor
product $\Pi_\al \otimes \pi^{k+2}_\la$ is a module over $V_k(\sw_2)$
and hence over $\su$. We denote it by $\wt{W}_{\al,\la,k}$. In
addition, we introduce the bosonic vertex operator
\begin{equation}    \label{second}
V_{2(k+2)}(z) = T_{2(k+2)} \exp \left( - \sum_{n<0}
\frac{b_n}{n} z^{-n} \right) \exp \left( - \sum_{n>0}
\frac{b_n}{n} z^{-n} \right),
\end{equation}
A straightforward computation similar to the one performed in
\secref{first screening} yields:

\begin{prop}    \label{screening second}
The residue
$$
\wt{S}_k = \int \wt{a}(z)^{-(k+2)} V_{2(k+2)}(z): \wt{W}_{0,0,k} \to
\wt{W}_{-(k+2),2(k+2),k}
$$
is an intertwining operator between the $\su$--modules
$\wt{W}_{0,0,k}$ and $\wt{W}_{-(k+2),2(k+2),k}$.
\end{prop}

We call $\wt{S}_k$, or its restriction to $W_{0,k} \subset
\wt{W}_{0,0,k}$, the {\em screening operator of the second kind} for
$\su$. This operator was first introduced by V. Dotsenko \cite{Do}.

\begin{prop}    \label{kernel}
For generic $k$ the $\su$--submodule $V_\ka(\sw_2) \subset W_{0,k}$ is
equal to the kernel of $\wt{S}_k: W_{0,k} \to
\wt{W}_{-(k+2),2(k+2),k}$.
\end{prop}

\begin{proof}
  We will show that for generic $k$ the kernel of $\wt{S}_k: W_{0,k}
  \to \wt{W}_{-(k+2),2(k+2),k}$ coincides with the kernel of the
  screening operator of the first kind, $S_k: W_{0,k} \to W_{-2,k}$.
  This, together with \propref{seq is exact}, will imply the statement
  of the proposition. The above two kernels are equal to the
  intersection of $W_{0,k} \subset \wt{W}_{0,0,k}$ and the kernels of
  the operators $\wt{S}_k$ and $S_k$ acting from $\wt{W}_{0,0,k}$ to
  $\wt{W}_{-(k+2),2(k+2),k}$ and $\wt{W}_{0,-2,k}$, respectively,
  (note that the operator $S_k$ has an obvious extension to an
  operator $\wt{W}_{0,0,k} \to \wt{W}_{0,-2,k}$). Therefore it is
  sufficient to show that the kernels of $\wt{S}_k$ and $S_k$ in
  $\wt{W}_{0,0,k}$ are equal for generic $k$.
  
  Now let $\phi(z)$ be the anti-derivative of $b(z)$, i.e., $b(z) =
  \pa_z \phi(z)$. By abusing notation we will write $V_{-2}(z) =
  e^{-(k+2)^{-1} \phi}$ and $V_{2(k+2)}(z) = e^{\phi}$. Then the
  screening currents $S_k(z)$ and $\wt{S}_k(z)$ become:
$$
S_k(z) = e^{u+v-(k+2)^{-1}\phi(z)}, \qquad \wt{S}_k(z) =
e^{-(k+2)u-(k+2)v+\phi}.
$$

Consider a more general situation: let $\h$ be an abelian Lie algebra
with a non-degenerate inner product $\ka$. Using this inner product,
we identify $\h$ with $\h^*$. Let $\hh_\ka$ be the Heisenberg Lie
algebra and $\pi^\ka_\la, \la \in \h^* \simeq \h$ be its Fock
representations, defined as in \secref{pi0}. Then for any $\chi \in
\h^*$ there is a vertex operator $V^\ka_\chi(z): \pi^\ka_0 \to
\pi^\ka_\chi$ given by the formula
\begin{equation}    \label{vertex op}
V^\ka_\chi(z) = T_\chi \exp \left( - \sum_{n<0}
\frac{\chi_n}{n} z^{-n} \right) \exp \left( - \sum_{n>0}
\frac{\chi_n}{n} z^{-n} \right),
\end{equation}
where the $\chi_n = \chi \otimes t^n$ are the elements of the
Heisenberg Lie algebra $\hh_\ka$ corresponding to $\chi \in \h^*$,
which we identify with $\h$ using the inner product $\ka$.

Suppose that $\ka(\chi,\chi) \neq 0$, and denote by $\chi^\vee$ the
element of $\h$ equal to $-2\chi/\ka(\chi,\chi)$. We claim that if
$\chi$ is generic (i.e., away from countably many hypersurfaces in
$\h$) then the kernels of $\ds \int V^\ka_\chi(z) dz$ and $\ds \int
V^\ka_{\chi^\vee}(z) dz$ in $\pi^\ka_0$ coincide. Indeed, each of
these kernels is equal to the tensor product of the subspace in
$\pi^\ka_0$ generated from the vacuum vector by elements of $\hh_\ka$
which commute with $\chi_n$ (they span a Heisenberg Lie subalgebra of
$\hh_\ka$ corresponding to the orthogonal complement of $\chi$ in
$\h$; note that by our assumption on $\chi$, this orthogonal
complement does not contain $\chi$) and the kernel of our operator on
the Fock representation of the Heisenberg Lie algebra generated by
$\chi_n, n \in \Z$.  But the latter two kernels coincide for generic
values of $\ka(\chi,\chi)$, as shown in \cite{vertex}, \S~ 14.4.15.
Hence the kernels of $\ds \int V^\ka_\chi(z) dz$ and $\ds \int
V^\ka_{\chi^\vee}(z) dz$ also coincide generically.

Our situation corresponds to the $3$--dimensional Lie algebra $\h$
with a basis $\ol{u},\ol{v},\ol{\phi}$ with the following non-zero
inner products of the basis elements:
$$
\ka(\ol{u},\ol{u}) = -\ka(\ol{v},\ol{v}) = 1, \qquad
\ka(\ol{\phi},\ol{\phi}) = 2(k+2).
$$
Our screening currents $S_k(z)$ and $\wt{S}_k(z)$ are equal to
$V^\ka_\chi(z)$ and $V^\ka_{\chi^\vee}(z)$, where
$$
\chi = \ol{u}+\ol{v}-(k+2)^{-1}\ol{\phi}, \qquad \chi^\vee = -
(k+2) \chi = -(k+2)\ol{u}-(k+2)\ol{v}+\ol{\phi}.
$$
Therefore for generic $k$ the kernels of the screening operators $S_k$
and $\wt{S}_k$ coincide.
\end{proof}

\section{Intertwining operators for an arbitrary $\g$}

In this section we construct screening operators between Wakimoto
modules over $\ghat_\ka$ for an arbitrary simple Lie algebra $\g$ and
use them to characterize $V_\ka(\g)$ inside $W_{0,\ka}$.

\subsection{Parabolic induction}    \label{parab}

Denote by $\sw_2^{(i)}$ the Lie subalgebra of $\g$, isomorphic to
$\sw_2$, which is generated by $e_i, h_i$, and $f_i$. Let $\pp^{(i)}$
be the parabolic subalgebra of $\g$ spanned by $\bb_-$ and $e_i$, and
$\mm^{(i)}$ its Levi subalgebra. Thus, $\mm^{(i)}$ is equal to the
direct sum of $\sw_2^{(i)}$ and the orthogonal complement $\h^\perp_i$
of $h_i$ in $\h$.

We apply to $\pp^{(i)}$ the results on semi-infinite parabolic
induction of \secref{parabolic}. According to \corref{functorial}, we
obtain a functor from the category of smooth representations of $\su
\oplus \hh^\perp_{i,\ka_0}$, with $k$ and $\ka_0$ satisfying the
conditions of \thmref{noncrit1}, to the category of smooth
$\ghat_{\ka+\ka_c}$--modules. The condition on $k$ and $\ka_0$ is that
the inner products on $\sw_2^{(i)}$ corresponding to $(k+2)$ and
$\ka_0$ are both restrictions of an invariant inner product
$(\ka-\ka_c)$ on $\g$. In other words, $(\ka-\ka_c)(h_i,h_i) = 2(k+2)$
and $\ka_0 = \ka|_{\h^\perp_i}$. By abuse of notation we will write
$\ka$ for $\ka_0$. If $k$ and $\ka$ satisfy this condition, then for
any smooth $\su$--module $R$ of level $k$ and any smooth
$\hh^\perp_{\ka}$--module the tensor product $M_{\g,\pp^{(i)}} \otimes
R \otimes L$ is a smooth $\ghat_{\ka+\ka_c}$--module.

In particular, if we choose $R$ to be the Wakimoto module $W_{\la,k}$
over $\sw_2$, and $L$ to be the Fock representation $\pi^\ka_{\la_0}$,
the corresponding $\ghat_\ka$--module will be isomorphic to the
Wakimoto module $W_{(\la,\la_0),\ka+\ka_c}$, where $(\la,\la_0)$ is
the weight of $\g$ built from $\la$ and $\la_0$. Under this
isomorphism the generators $a_{\al_i,n}, n \in \Z$, will have a
special meaning: they correspond to the {\em right} action of the
elements $e_{i,n}$ of $\ghat$, which was defined in \remref{right
action}. In other words, making the above identification of modules
forces us to choose a system of coordinates $\{ y_\al\}_{\al \in
\De_+}$ on $N_+$ such that $\rho^R(e_i) = \pa/\pa y_{\al_i}$ (in the
notation of \secref{sec on formulas}), and so $w^R(e_i(z)) =
a_{\al_i}(z)$ (in the notation of \remref{right action}). {}From now
on we will denote $w^R(e_i(z))$ by $e^R_i(z)$. For a general
coordinate system on $N_+$ we have
\begin{equation} \label{e right1}
e^R_i(z) = a_{\al_i}(z) + \sum_{\beta \in \Delta_+}
P^{R,i}_\beta(a^*_\al(z)) a_\beta(z)
\end{equation}
(see formula \eqref{e right}).

Now any intertwining operator between Wakimoto modules $W_{\la_1,k}$
and $W_{\la_2,k}$ over $\su$ gives rise to an intertwining operator
between the Wakimoto modules $W_{(\la_1,\la_0),\ka+\ka_c}$ and
$W_{(\la_2,\la_0),\ka+\ka_c}$ over $\ghat_{\ka+\ka_c}$ for any weight
$\la_0$ of $\h^\perp_i$. We will use this fact and the $\su$ screening
operators introduced in the previous section to construct intertwining
operators between Wakimoto modules over $\g$.

\subsection{Screening operators of the first kind}

Let $\ka$ be a non-zero invariant inner product on $\g$. We will use
the same notation for the restriction of $\ka$ to $\h$. Now to any
$\chi \in \h$ we associate a vertex operator $V^\ka_\chi(z): \pi^\ka_0
\to \pi^\ka_\chi$ given by formula \eqref{vertex op}. For $\ka \neq
\ka_c$, we set
$$
S_{i,\ka}(z) \overset{\on{def}}{=} e_i^R(z) V^{\ka-\ka_c}_{-\al_i}(z):
W_{0,\ka} \to W_{-\al_i,\ka},
$$
where $e^R_i(z)$ is given by formula \eqref{e right}. Note that
\begin{equation}    \label{YVM}
S_{i,\ka}(z) = Y_{W_{0,\ka},W_{-\al_i,\ka}}(e^R_{i,-1}
v_{-\al_i,\ka},z)
\end{equation}
in the notation of \secref{generalities}.

According to \propref{screening first} and the above discussion, the
operator
$$S_{i,\ka} = \ds \int S_{i,\ka}(z) dz: W_{0,\ka} \to W_{-\al_i,\ka}$$
is induced by the screening operator of the first kind $S_k$ for the
$i$th $\su$ subalgebra, where $k$ is determined from the formula
$(\ka-\ka_c)(h_i,h_i) = 2(k+2)$. Hence \propref{screening first}
implies:

\begin{prop}    \label{ith int}
The operator $S_{i,\ka}$ is an intertwining operator between the
$\ghat_\ka$--modules $W_{0,\ka}$ and $W_{-\al_i,\ka}$ for each
$i=1,\ldots,\ell$.
\end{prop}

We call $S_{i,\ka}$ the $i$th screening operator of the first kind for
$\ghat_\ka$.

Recall that by \propref{injective} $V_\ka(\g)$ is naturally a
$\ghat_\ka$--submodule and a vertex subalgebra of $W_{0,\ka}$. On the
other hand, the intersection of the kernels of $S_{i,\ka},
i=1,\ldots,\ell$, is a $\ghat_\ka$--submodule of $W_{0,\ka}$, by
\propref{ith int}, and a vertex subalgebra of $W_{0,\ka}$, due to
formula \eqref{YVM} and the commutation relations \eqref{kern}. The
following proposition is proved in \cite{FF:wak}.

\begin{prop}    \label{intersection1}
For generic $\ka$, $V_\ka(\g)$ is equal to the intersection of the
kernels of the screening operators $S_{i,\ka}, i=1,\ldots,\ell$.
\end{prop}

Furthermore, in \cite{FF:wak}, \S~3, a complex $C^\bullet(\g)$ of
$\ghat_\ka$--modules is constructed for generic $\ka$. Its $i$th
degree term is the direct sum of the Wakimoto modules
$W_{w(\rho)-\rho,\ka}$, where $w$ runs over all elements of the Weyl
group of $\g$ of length $i$. Its $0$th cohomology is isomorphic to
$V_\ka(\g)$, and all other cohomologies vanish. For $\g=\sw_2$ this
complex has length one and coincides with the one appearing in
\propref{seq is exact}. In general, the degree zero term of the
complex is $W_{0,\ka}$, the degree one term is $\bigoplus_{i=1}^\ell
W_{-\al_i,\ka}$, and the zeroth differential is the sum of the
screening operators $S_{i,\ka}$.

In \cite{FF:wak} it is also explained how to construct other
intertwining operators as compositions of the screening operators
$S_{i,\ka}$. Roughly speaking, the screening operators $S_{i,\ka}$
satisfy the $q$--Serre relations, i.e., the defining relations of the
quantized enveloping algebra $U_q(\n_+)$ with appropriate parameter
$q$. Then for generic $\ka$ we attach to a singular vector in the
Verma module $M_\la$ over $U_q(\g)$ of weight $\mu$ an intertwining
operator $W_{\la,\ka} \to W_{\mu,\ka}$. This operator is equal to the
integral of a product of the screening currents $S_{i,\ka}(z)$ over a
certain cycle on the configuration space with coefficients in a local
system that is naturally attached to the above singular vector.

\subsection{Screening operators of the second kind}

In order to define the screening operators of the second kind, we
need to make sense of the series $(e^R_i(z))^\al$. So we choose a
system of coordinates on $N_+$ in such a way that $e^R_i(z) =
a_{\al_i}(z)$ (note that this cannot be achieved for all
$i=1\ldots,\ell$ simultaneously). This is automatically so if we
define the Wakimoto modules over $\ghat$ via the semi-infinite
parabolic induction from Wakimoto modules over the $i$th subalgebra
$\su$ (see \secref{parab}).

Having chosen such a coordinate system, we define the series
$a_{\al_i}(z)^\al$ using the Friedan--Martinec--Shenker bosonization
of the Weyl algebra generated by $a_{\al_i,n}, a^*_{\al_i,n}$, $n \in
\Z$, as explained in \secref{bosonization}. Namely, we have a vertex
algebra
$$
\Pi_0^{(i)} = \C[a_{\al_i,n}]_{n<-1} \otimes
\C[a^*_{\al_i,n}]_{n\leq 0} \otimes \C[a_{\al_i,-1}^{\pm 1}]
$$
containing
$$
M^{(i)}_\g = \C[a_{\al_i,n}]_{n\leq -1} \otimes
\C[a^*_{\al_i,n}]_{n\leq 0}
$$
and a $\Pi_0^{(i)}$--module $\Pi_{\ga}^{(i)}$ defined as in
\secref{bosonization}. We then set
$$
\wt{W}^{(i)}_{\ga,\la,\ka} = W_{\la,\ka} \underset{M^{(i)}_\g}\otimes
\Pi_{\ga}^{(i)}.
$$
This is a $\ghat_\ka$--module, which contains $W_{\la,\ka}$ if
$\al=0$. Note that $\wt{W}^{(i)}_{0,0,\ka_c}$ is the $\ghat$--module
obtained by the semi-infinite parabolic induction from the
$\su$--module $\wt{W}_{0,0,-2}$.

Now let $\beta = \frac{1}{2} (\ka-\ka_c)(h_i,h_i)$ and define the
field
$$
\wt{S}_{\ka,i}(z) \overset{\on{def}}{=} (e^R_i(z))^{-\beta}
V_{\al_i^\vee}(z): \wt{W}^{(i)}_{0,0,\ka} \to
\wt{W}^{(i)}_{-\beta,\beta\al_i^\vee,\ka}.
$$
Here $\al_i^\vee = h_i \in \h$ denotes the $i$th coroot of $\g$.
Then the operator $\wt{S}_{\ka,i} = \ds \int \wt{S}_{\ka,i}(z) dz$ is
induced by the screening operator of the second kind $\wt{S}_k$ for
the $i$th $\su$ subalgebra, where $k$ is determined from the formula
$(\ka-\ka_c)(h_i,h_i) = 2(k+2)$. Hence \propref{screening first}
implies (see also \cite{PRY}):

\begin{prop}
The operator $\wt{S}_{\ka,i}$ is an intertwining operator between the
$\ghat_\ka$--modules $\wt{W}^{(i)}_{0,0,\ka}$ and
$\wt{W}^{(i)}_{-\beta,\beta\al_i^\vee,\ka}$.
\end{prop}

Note that $\wt{S}_{i,\ka}$ is the residue of
$Y_{V,M}((e^R_{i,-1})^{-\beta} v_{\al_i^\vee,\ka},z)$, where
$V=\wt{W}^{(i)}_{0,0,\ka}$ and
$M=\wt{W}^{(i)}_{-\beta,\beta\al_i^\vee,\ka}$ (see
\secref{generalities}). Therefore, according to the commutation
relations \eqref{kern}, the intersection of kernels of
$\wt{S}_{i,\ka}, i=1,\ldots,\ell$, is naturally a vertex subalgebra of
$\wt{W}^{(i)}_{0,0,\ka}$ or $W_{0,\ka}$. Combining
\propref{intersection1} and \propref{kernel}, we obtain:

\begin{prop}    \label{intersection}
For generic $\ka$, $V_\ka(\g)$ is equal, as a $\ghat_\ka$--module and
as a vertex algebra, to the intersection of the kernels of the
screening operators
$$
\wt{S}_{i,\ka}: W_{0,\ka} \to \wt{W}^{(i)}_{-\beta,\beta\al_i^\vee,\ka},
\qquad i=1,\ldots,\ell.
$$
\end{prop}

One can use the screening operators $\wt{S}_{i,\ka}$ to construct more
general intertwining operators following the procedure of
\cite{FF:wak}.

\subsection{Screening operators of second kind at the critical level}

Now we define the limits of the intertwining operators
$\wt{S}_{i,\ka}$ as $\ka \to \ka_c$.

First we consider the case when $\g=\sw_2$.  Denote $\beta=k+2$. In
order to define the limit of $\wt{S}_k$ as $\beta \to 0$ (i.e., $k \to
-2$) we replace $\wt{W}_{0,0,k}$ and $\wt{W}_{-\beta,2\beta,k}$ by
their tensor products with $\C[\beta^{\pm 1}]$ (still denoted by the
same symbols) and define $\C[\beta]$--lattices in them which are
preserved by $\wt{S}_k$. These lattices are spanned over $\C[\beta]$
by the monomials in $b_n, p_n, q_n, n<0$, applied to the highest
weight vectors.

Consider first the expansion of $V_{2\beta}(z)$ in powers of $\beta$
with respect to the above lattice. Recall that
$$
[b_n,b_m] = 2\beta n \delta_{n,-m}.
$$
Let us write $V_{2\beta}(z) = \sum_{n\in\Z} V_{2\beta}[n] z^{-n}$. We
will identify $\pi^\beta_{2\beta}$ and $\pi^\beta_0$ with $\pi^0 =
\C[b_{-m}]_{m>0}$ in the limit $\beta \to 0$. Introduce the operators
$\ol{V}[n], n\leq 0$, via the formal power series
$$
\sum_{n\leq 0} \ol{V}[n] z^{-n} = \exp \left( \sum_{m>0}
\frac{b_{-m}}{m} z^m \right).
$$
Using formula \eqref{second}, we obtain the following expansion of
$V_{2\beta}[n]$:
\begin{equation*}
V_{2\beta}[n] = \begin{cases} \ol{V}[n] + \beta(\ldots), \quad n\leq
0, \\ - 2 \beta \sum_{m\leq 0} \ol{V}[m] \ds \frac{\pa}{\pa b_{m-n}} +
\beta^2(\ldots), \quad n>0. \end{cases}
\end{equation*}
We denote the $\beta$--linear term of $V_{2\beta}[n], n>0$, by
$\ol{V}[n]$.

Next, we consider the expansion of $\wt{a}(z)^{-\beta} =
e^{-\beta(u+v)}$ in powers of $\beta$. Let us write
$\wt{a}(z)^{-\beta} = \sum_{n\in\Z} \wt{a}(z)^{-\beta}_{[n]}
z^{-n}$. In the limit $\beta \to 0$ we will identify $\Pi_{-\beta}$
with $\Pi_0$. Then we find that
$$
\wt{a}(z)^{-\beta}_{[n]} = \begin{cases} 1 + \beta(\ldots), \quad n=0,
  \\ \ds  \beta \frac{p_n+q_n}{n} + \beta^2(\ldots), \quad n\neq 0.
\end{cases}
$$

The above formulas imply the following expansion of the screening
operator $\wt{S}_k$:
$$
\wt{S}_k = \beta \left( \ol{V}[1] + \sum_{n>0} \frac{1}{n}
\ol{V}[-n+1] (p_n+q_n) \right) + \beta^2(\ldots).
$$
Therefore the leading term of $\wt{S}_k$ at $k=-2$ is the operator
\begin{equation}    \label{olS}
\ol{S} \overset{\on{def}}{=} \ol{V}[1] + \sum_{n>0} \frac{1}{n}
\ol{V}[-n+1] (p_n+q_n),
\end{equation}
acting on $\wt{W}_{0,0,-2} = \Pi_0 \otimes \pi^0$. It follows from the
construction that $\ol{S}$ commutes with the $\su$--action on
$\wt{W}_{0,0,-2}$.

It is possible to express the operators
$$
\frac{p_n+q_n}{n} = -(u_n+v_n), \quad \quad n>0,
$$
in terms of the Heisenberg algebra generated by $a_m, a^*_m, m \in
\Z$. Namely, from the definition $\wt{a}(z) = e^{u+v}$ it follows that
$u(z)+v(z) = \log \wt{a}(z)$, so that the field $u(z)+v(z)$ commutes
with $a(w)$, and we have the following OPE with $a^*(w)$:
$$
(u(z) + v(z)) a^*(w) = \frac{1}{z-w} \wt{a}(w)^{-1} + \on{reg.}
$$ Therefore, writing $\wt{a}(z)^{-1} = \sum_{n \in \Z}
\wt{a}(z)^{-1}_{[n]} z^{-n}$, we obtain the following commutation
relations:
\begin{equation}    \label{com pn}
\left[ \frac{p_n+q_n}{n},a^*_m \right] = - \wt{a}(z)^{-1}_{[n+m-1]},
\qquad n>0.
\end{equation}
Using the realization \eqref{localiz} of $\Pi_0$, the series
$\wt{a}(z)^{-1}$ is expressed as follows:
\begin{equation}    \label{a inverse}
\wt{a}(z)^{-1} = (a_{-1})^{-1} \left( 1 + (a_{-1})^{-1} \sum_{n \neq
-1} a_n z^{-n-1} \right)^{-1},
\end{equation}
where the right hand side is expanded as a formal power series in
positive powers of $(a_{-1})^{-1}$. It is easy to see that each Fourier
coefficient of this power series is well-defined as a linear operator
on $\Pi_0$.

The above formulas completely determine the action of $(p_n+q_n)/n,
n>0$, and hence of $\ol{S}$, on any vector in $\wt{W}_{0,-2} = \Pi_0
\otimes \pi^0$. Namely, we use the commutation relations \eqref{com
  pn} to move $(p_n+q_n)/n$ through the $a^*_m$'s. As the result, we
obtain Fourier coefficients of $\wt{a}(z)^{-1}$, which are given by
formula \eqref{a inverse}. Applying each of them to any vector in
$\wt{W}_{0,-2}$, we always obtain a finite sum.

This completes the construction of the limit $\ol{S}$ of the screening
operator $\wt{S}_k$ as $k \to -2$ in the case of $\sw_2$. Now we
consider the case of an arbitrary $\g$.

The limit of $\wt{S}_{i,\ka}$ as $\ka \to \ka_c$ is by definition the
operator $\ol{S}_i$ on $\wt{W}^{(i)}_{0,0,\ka_c}$, obtained from
$\ol{S}$ via the functor of semi-infinite parabolic
induction. Therefore it is given by the formula
\begin{equation}    \label{olSi}
\ol{S}_i = \ol{V}_i[1] + \sum_{n>0} \frac{1}{n}
\ol{V}_i[-n+1] (p_{i,n}+q_{i,n}).
\end{equation}
Here $\ol{V}_i[n]: \pi^0 \to \pi^0$ are the linear operators given by
the formulas
\begin{equation}    \label{olVin}
\sum_{n\leq 0} \ol{V}_i[n] z^{-n} = \exp \left( \sum_{m>0}
\frac{b_{i,-m}}{m} z^m \right),
\end{equation}
$$
\ol{V}_i[1] = - \sum_{m\leq 0} \ol{V}_i[m] \; {\mb D}_{b_{i,m-1}},
$$
where ${\mb D}_{b_{i,m}}$ denotes the derivative in the direction
of $b_{i,m}$ given by the formula
\begin{equation}    \label{Db}
{\mb D}_{b_{i,m}} \cdot b_{j,n} = a_{ji} \delta_{n,m},
\end{equation}
and $(a_{kl})$ is the Cartan matrix of $\g$ (it is normalized so that
we have ${\mb D}_{b_{i,m}} \cdot b_{i,n} = 2 \delta_{n,m}$). The
operators $(p_{i,n}+q_{i,n})/n$ acting on $\Pi_0^{(i)}$ are defined in
the same way as above.

Thus, we obtain well-defined linear operators $\ol{S}_i: W_{0,\ka_c}
\to \wt{W}^{(i)}_{0,0,\ka_c}$. By construction, they commute with the
action of $\ghat_{\ka_c}$ on both modules. It is clear that the
operators $\ol{S}_i, i=1,\ldots,\ell$, annihilate the highest weight
vector of $W_{0,\ka_c}$. Therefore they annihilate all vectors
obtained from the highest weight vector under the action of
$\ghat_{\ka_c}$, i.e., all vectors in $V_{\ka_c}(\g) \subset
W_{0,\ka_c}$. Thus we obtain

\begin{prop}    \label{int crit}
The vacuum module $V_{\ka_c}(\g)$ is contained in the intersection of
the kernels of the operators $\ol{S}_i: W_{0,\ka_c} \to
\wt{W}^{(i)}_{0,0,\ka_c}, i=1,\ldots,\ell$.
\end{prop}

\section{Description of the center of $V_{\ka_c}(\g)$}    \label{descr
  of center}

In this section we use \propref{int crit} to describe the center of
$V_{\ka_c}(\g)$. Recall from \cite{vertex}, \S~4.6.2, that the center
of a vertex algebra $V$ is by definition its commutative vertex
subalgebra spanned by all vectors $v \in V$ such that $Y(A,z) v \in
V[[z]]$ for all $A \in V$. In the case when $V=V_\ka(\g)$, the center
is equal to the space of invariants of the Lie subalgebra $L_+\g =
\g[[t]]$ of $\ghat_\ka$. According to Lemma 17.4.1 of \cite{vertex},
it is trivial, i.e., equal to the span of the vacuum vector $v_\ka$ if
$\ka\neq \ka_c$. So let us consider the center of $V_{\ka_c}(\g)$,
which we denote by $\zz(\ghat)$.

\subsection{The center of $V_{\ka_c}(\g)$ and Wakimoto modules}

Recall that by \propref{injective} the homomorphism $\ww_{\ka_c}:
V_{\ka_c}(\g) \hookrightarrow W_{0,\ka_c}$ is injective. The vertex
algebra $V_{\ka_c}(\g)$ contains the commutative subalgebra
$\zz(\ghat)$, its center. On the other hand, $W_{0,\ka_c} = M_\g
\otimes \pi_0$ contains the commutative subalgebra $\pi_0$, which is
its center.

\begin{lem}    \label{image}
The image of $\zz(\ghat) \subset V_{\ka_c}(\g)$ in $W_{\ka_c,0}$ under
$\ww_{\ka_c}$ is contained in $\pi_0 \subset W_{\ka_c,0}$.
\end{lem}

\begin{proof}
Observe that the lexicographically ordered monomials of the form
\begin{equation} \label{mon} \prod_{n_{r(\al)}<0} e^R_{\al,n_{r(\al)}}
\prod_{m_{s(\al)}\leq 0} a^*_{\al,m_{s(\al)}} \prod_{n_{t(i)}<0}
b_{i,n_{t(i)}} \vac, \end{equation} where the $e^R_{i,n}$'s are given
by formula \eqref{e right}, form a basis of $W_{0,\ka_c}$. The image of
any element of $\zz(\ghat)$ in $W_{\ka_c,0}$ is an $L_+\g$--invariant
vector. In particular, it is annihilated by $L_+\n_+ = \n_+[[t]]$ and
by $\h$. Since $L_+\n_+$ commutes with $e^R_{\al,n}$ and $b_{i,n}$,
the space of $L_+\n_+$--invariants of $W_{\ka_c,0}$ is equal to the
tensor product of the subspace $M_{\g,-}$ of $W_{0,\ka_c}$ spanned by
all monomials \eqref{mon} not containing $a^*_{\al,n}$, and the space
of $L_+\n_+$--invariants in $M_{\g,+} = \C[a^*_{\al,n}]_{\al \in
\De_+,n\leq 0}$. According to \secref{M+}, $M_{\g,+}$ is an
$L_+\g$--module isomorphic to $\on{Coind}_{L_+\bb_-}^{L_+\g}
\C$. Therefore the action of $L_+\n_+$ on it is co-free, and the space
of $L_+\n_+$--invariants is one-dimensional, spanned by constants.

Thus, we obtain that the space of $L_+\n_+$--invariants in
$W_{0,\ka_c}$ is equal to $M_{\g,-}$. However, its subspace of
$\h$--invariants is the span of the monomials \eqref{mon}, which only
involve the $b_{i,n}$'s, i.e., $\pi_0 \subset W_{0,\ka_c}$. This
completes the proof.
\end{proof}

Using \propref{int crit}, we now obtain that $\zz(\ghat)$ is contained
in the intersection of the kernels of the operators $\ol{S}_i,
i=1,\ldots,\ell$, restricted to $\pi_0 \subset W_{0,\ka_c}$. But
according to formula \eqref{olSi}, the restriction of $\ol{S}_i$ to
$\pi_0$ is nothing but the operator $\ol{V}_i[1]: \pi_0 \to \pi_0$
given by the formula
\begin{equation}    \label{olVi}
\ol{V}_i[1] = - \sum_{m\leq 0} \ol{V}_i[m] \; {\mb D}_{b_{i,m-1}}
\end{equation}
where $\ol{V}_i[m], {\mb D}_{b_{i,m-1}}, m\leq 0$, are given by
formulas \eqref{olVin} and \eqref{Db}, respectively. Therefore we
obtain the following

\begin{prop}    \label{center}
  The center $\zz(\ghat)$ of $V_{\ka_c}(\g)$ is contained in the
  intersection of the kernels of the operators $\ol{V}_i[1],
  i=1,\ldots,\ell$, in $\pi_0$.
\end{prop}

We now show that $\zz(\ghat)$ is actually equal to the intersection of
  the kernels of the operators $\ol{V}_i[1], i=1,\ldots,\ell$.

\subsection{The associated graded of $\zz(\ghat)$}    \label{ass graded}

The $\ghat_{\ka_c}$--module $V_{\ka_c}(\g)$ has a natural filtration
induced by the Poincar\'e--Birkhoff--Witt filtration on the universal
enveloping algebra $U(\ghat_{\ka_c})$, and the action of
$\ghat_{\ka_c}$ preserves this filtration. Furthermore, it is clear
that the associated graded space $\on{gr} V_{\ka_c}(\g)$ is isomorphic
to
$$
\on{Sym} \g((t))/\g[[t]] \simeq \on{Fun} \g[[t]],
$$
where we use the following (coordinate-dependent) pairing
\begin{equation}    \label{pairing}
\langle A \otimes f(t),B \otimes g(t) \rangle = \ka(A,B) \on{Res} f(t)
g(t) \frac{dt}{t},
\end{equation}
and $\ka$ is an arbitrary fixed non-degenerate inner product on $\g$. 

Recall that $\zz(\ghat)$ is equal to the space of
$\g[[t]]$--invariants in $V_{\ka_c}(\g)$. The symbol of a
$\g[[t]]$--invariant vector in $V_{\ka_c}(\g)$ is a
$\g[[t]]$--invariant vector in $\on{gr} V_{\ka_c}(\g)$, i.e., an
element of the space of $\g[[t]]$--invariants in $\on{Fun}
\g[[t]]$. Hence we obtain an injective map
\begin{equation}    \label{injection}
\on{gr} \zz(\ghat) \hookrightarrow (\on{Fun} \g[[t]])^{\g[[t]]}.
\end{equation}

In order to describe the space $(\on{Fun} \g[[t]])^{\g[[t]]}$, recall
that $(\on{Fun} \g)^\g$ is a polynomial algebra in generators $P_i,
i=1,\ldots,\ell$, of degrees $d_i+1$, where $d_i$ is the $i$th
exponent of $\g$. Each $P_i$ is a polynomial in the basis elements
$J^a$ of $\g$ (considered as linear functionals of $\g$ via the inner
product $\ka$). Substituting the generating functions
$$
J^a(t) = \sum_{n<0} J^a_n t^{-n-1}
$$
of the basis elements $J^a_n$ of $t^{-1}\g[t^{-1}]$ into $P_i$ instead
of the $J^a$'s, we obtain a generating series
$$
P_i(J^a(t)) = \sum_{n\geq 0} P_{i,n} t^n
$$
of elements of $\on{Fun} \g[[t]]$. Clearly, its coefficients
$P_{i,n}$ is $\g[[t]]$--invariant. Furthermore we have the following
result due to Beilinson and Drinfeld, \cite{BD}, \S~2.4.1. For
reader's convenience, we reproduce the proof (see also \cite{Mus},
Appendix, Prop. A.1):

\begin{prop}    \label{free pol}
The algebra $(\on{Fun} \g[[t]])^{\g[[t]]}$ is the free polynomial
algebra in the generators $P_{i,n}, i=1,\ldots,\ell; n \geq 0$.
\end{prop}

\begin{proof}
Let $G$ be the connected simply-connected algebraic group with the Lie
algebra $\g$. Then the space of $\g[[t]]$--invariants in $\on{Fun}
\g[[t]]$ coincides with the space of $G[[t]]$--invariants in $\on{Fun}
\g[[t]]$.

Recall that given an algebraic variety $X$, we denote by $JX$ its
infinite jet scheme, as defined in \cite{vertex}, \S~ 8.4.4 (see also
\cite{Mus}, where the notation $X_\infty$ is used). In particular,
$J\g = \g[[t]]$ and $JG = G[[t]]$.

Let $\g_{\on{reg}}$ be the smooth open subset of $\g$ consisting of
regular elements. It is known that the morphism $$\chi: \g_{\on{reg}}
\to {\mc P} := \Spe (\on{Fun} \g)^\g = \Spe \C[P_1,\ldots,P_\ell]$$ is
smooth and surjective (see \cite{Ko}). Therefore the morphism
$$J\chi: J\g_{\on{reg}} \to J{\mc P} := \Spe
\C[P_{1,m},...,P_{\ell,m}]_{m\geq 0}$$ is also smooth and
surjective.

Consider the map $a: G \times \g_{\on{reg}} \to \g_{\on{reg}}
\underset{{\mc P}}\times \g_{\on{reg}}$ defined by the formula $a(g,x)
= (x,g \cdot x)$. The map $a$ is smooth, and since $G$ acts
transitively along the fibers of $\chi$, it is also surjective. Hence
the corresponding map of jet schemes $Ja: JG \times J\g_{\on{reg}} \to
J\g_{\on{reg}} \underset{J{\mc P}}\times J\g_{\on{reg}}$ is
surjective. Given two points $y_1,y_2$ in the same fiber of $J\chi$,
let $(h,y_1)$ be a point in the (non-empty) fiber
$(Ja)^{-1}(y_1,y_2)$. Then $y_2 = h \cdot y_1$. Hence $JG$ acts
transitively along the fibers of the map $J\chi$.  This implies that
the ring of $JG$--invariant functions on $J\g_{\on{reg}}$ is equal to
$\C[P_{1,m},...,P_{\ell,m}]_{m\geq 0}$. Because $\g_{\on{reg}}$ is an
open dense subset of $\g$, we obtain that $J\g_{\on{reg}}$ is dense in
$J\g$, and so any $JG$--invariant function on $J\g$ is determined by
its restriction to $J\g_{\on{reg}}$. This proves the proposition.
\end{proof}

We have an action of the operator $L_0 = - t \pa_t$ on $\on{Fun}
\g[[t]]$. It defines a $\Z$--gradation on $\on{Fun} \g[[t]]$ such
that $\deg J^a_n = -n$. Then $\deg P_{i,n} = d_i+n+1$, and we obtain
the following formula for the character of $(\on{Fun}
\g[[t]])^{\g[[t]]}$ (i.e., the generating function of its graded
dimensions):
\begin{equation}    \label{character of vac}
\on{ch} (\on{Fun} \g[[t]])^{\g[[t]]} = \prod_{i=1}^\ell \prod_{n_i
  \geq d_i+1} (1-q^{n_i})^{-1}.
\end{equation}

\subsection{Computation of the character of $\zz(\ghat)$}
  \label{char of zz}

We now show that the map \eqref{injection} is an isomorphism using
\propref{Wak and Verma}.

Consider the Lie subalgebra $\wt{\bb}_+ = (\bb_+ \otimes 1) \oplus (\g
\otimes t\C[[t]])$ of $\g[[t]]$. The natural surjective
homomorphism $M_{0,\ka_c} \to V_{\ka_c}(\g)$ gives rise to a map
of $\wt{\bb}_+$--invariants
$$
\phi: M_{0,\ka_c}^{\wt{\bb}_+} \to V_{\ka_c}(\g)^{\wt{\bb}_+}.
$$
Both $M_{0,\ka_c}$ and $V_{\ka_c}(\g)$ have natural filtrations
which are preserved by the homomorphism between them. Therefore we
have the corresponding map of associated graded
$$
\phi_{\on{cl}}: (\on{gr} M_{0,\ka_c})^{\wt{\bb}_+} \to (\on{gr}
V_{\ka_c}(\g))^{\wt{\bb}_+}.
$$
Since $V_{\ka_c}(\g)$ is a direct sum of finite-dimensional
representations of the constant subalgebra $\g \otimes 1$ of
$\g[[t]]$, we find that any $\wt{\bb}_+$--invariant in $V_{\ka_c}(\g)$
or $\on{gr} V_{\ka_c}(\g)$ is automatically a
$\g[[t]]$--invariant. Thus, we have
\begin{align*}
V_{\ka_c}(\g)^{\wt{\bb}_+} &= V_{\ka_c}(\g)^{\g[[t]]}, \\ (\on{gr}
V_{\ka_c}(\g))^{\wt{\bb}_+} &= (\on{gr}
V_{\ka_c}(\g))^{\g[[t]]} = \C[P_{1,m},...,P_{\ell,m}]_{m\geq
0}.
\end{align*}

We need to describe $(\on{gr} M_{0,\ka_c})^{\wt{\bb}_+}$. First,
observe that
$$
\on{gr} M_{0,\ka_c} = \on{Sym} \g((t))/\wt{\bb}_+ \simeq \on{Fun}
\wt{\n}^{(-1)}_+,
$$
where
$$
\wt{\n}^{(-1)}_+ = (\n_+ \otimes t^{-1}) \oplus \g[[t]],
$$
and we use the pairing \eqref{pairing} on $\g((t))$. In terms of this
identification, the map $\phi_{\on{cl}}$ becomes a ring homomorphism
$$
\on{Fun} \wt{\n}^{(-1)}_+ \to \on{Fun} \g[[t]]
$$
induced by the natural embedding $\g[[t]] \to \wt{\n}^{(-1)}_+$.

We construct $\wt{\bb}_+$--invariant functions on $\wt{\n}^{(-1)}_+$
in the same way as above, by substituting the generating functions
$J^a(t)$ into the $P_i$'s. However, now $J^a(t)$ has non-zero $t^{-1}$
coefficients if the inner product of $J^a$ with $\n_+$ is
non-zero. Therefore the resulting series $P_i(J^a(t))$ will have
non-zero coefficients $P_{i,m}$ in front of $t^m$ if and only if
$m\geq -d_i$. Thus, we obtain a natural inclusion
$$
\C[P_{i,m_i}]_{i=1,\ldots,\ell;m_i\geq -d_i} \subset (\on{Fun}
\wt{\n}^{(-1)}_+)^{\wt{\bb}_+}.
$$

\begin{lem}
This inclusion is an equality.
\end{lem}

\begin{proof}
Let $\wt{\n}_+ = (\n_+ \otimes 1) \oplus t\g[[t]] = t
\wt{\n}^{(-1)}_+$. Clearly, the spaces of $\wt{\bb}_+$--invariant
functions on $\wt{\n}^{(-1)}_+$ and $\wt{\n}_+$ are isomorphic, so we
will consider the latter space.

Denote by $\wt{\n}_+^{\on{reg}}$ the intersection of $\wt{\n}_+$ and
$J\g_{\on{reg}} = \g_{\on{reg}} \oplus t\g[[t]]$. Thus,
$\wt{\n}_+^{\on{reg}} = \wt{\n}_+^{\on{reg}} \oplus t\g[[t]]$, where
$\n_+^{\on{reg}} = \n_+ \cap \g_{\on{reg}}$ is an open dense subset of
$\n_+$, so that $\wt{\n}_+^{\on{reg}}$ is open and dense in
$\wt{\n}_+$.

Recall the morphism $J\chi: J\g_{\on{reg}} \to J{\mc P}$ introduced in
the proof of \propref{free pol}. It was shown there that the group $JG
= G[[t]]$ acts transitively along the fibers of $J\chi$. Let
$\wt{B}_+$ be the subgroup of $G[[t]]$ corresponding to the Lie
algebra $\wt{\bb}_+ \subset \g[[t]]$. Note that for any $x \in
\wt{\n}_+^{\on{reg}}$, the group $\wt{B}_+$ is equal to the subgroup
of all elements $g$ of $G[[t]]$ such that $g \cdot x \in
\wt{\n}_+^{\on{reg}}$. Therefore $\wt{B}_+$ acts transitively along
the fibers of the restriction of the morphism $J\chi$ to
$\wt{\n}_+^{\on{reg}}$.

This implies that the ring of $\wt{B}_+$--invariant (equivalently,
$\wt{\bb}_+$--invariant) polynomials on $\wt{\n}_+^{\on{reg}}$ is the
ring of functions on the image of $\wt{\n}_+^{\on{reg}}$ in $J{\mc P}$
under the map $J\chi$. But it follows from the construction that the
image of $\wt{\n}_+^{\on{reg}}$ in $J{\mc P}$ is the subspace
determined by the equations $P_{i,m} = 0, i=1,\ldots,\ell; m>0$. Hence
the ring of $\wt{\bb}_+$--invariant polynomials on
$\wt{\n}_+^{\on{reg}}$ is equal to $\C[P_{i,m}]_{i=1,\ldots,\ell;
m>0}$. Since $\wt{\n}_+^{\on{reg}}$ is dense in $\wt{\n}_+$, we obtain
that this is also the ring of invariant polynomials on
$\wt{\n}_+$. When we pass from $\wt{\n}_+$ to $\wt{\n}^{(-1)}_+$, we
obtain the statement of the lemma.
\end{proof}

\begin{cor}    \label{surj}
The map $\phi_{\on{cl}}$ is surjective.
\end{cor}

Since $\on{deg} P_{i,m} = d_i+m+1$, we obtain the character of
$(\on{gr} M_{0,\ka_c})^{\wt{\bb}_+} = (\on{Fun}
\wt{\n}^{(-1)}_+)^{\wt{\bb}_+}$:
$$
\on{ch} \, (\on{gr} M_{0,\ka_c})^{\wt{\bb}_+} = \prod_{m>0}
(1-q^m)^{-\ell}.
$$

Now recall that by \thmref{Wak and Verma}, the Verma module
$M_{0,\ka_c}$ is isomorphic to the Wakimoto module
$W^+_{0,\ka_c}$. Hence $(M_{0,\ka_c})^{\wt{\bb}_+} =
(W^+_{0,\ka_c})^{\wt{\bb}_+}$. In addition, according to \lemref{sing
vect imag} we have $(W^+_{0,\ka_c})^{\wt{\bb}_+} = \pi_0$, and so its
character is also equal to $\ds \prod_{m>0}
(1-q^m)^{-\ell}$. Therefore we find that the natural embedding
$$
\on{gr} (M_{0,\ka_c}^{\wt{\bb}_+}) \hookrightarrow (\on{gr}
M_{0,\ka_c})^{\wt{\bb}_+}
$$
is an isomorphism. Now consider the commutative diagram
$$
\begin{CD}
\on{gr} (M_{0,\ka_c}^{\wt{\bb}_+}) @>>> \on{gr}
(V_{\ka_c}(\g)^{\g[[t]]}) \\ @VVV @VVV \\ (\on{gr}
M_{0,\ka_c})^{\wt{\bb}_+} @>>> (\on{gr}
V_{\ka_c}(\g))^{\g[[t]]}.
\end{CD}
$$
It follows from the above discussion that the left vertical arrow is
an isomorphism. Moreover, by \corref{surj} the lower horizontal arrow
is surjective. Therefore we obtain an isomorphism
$$
\on{gr} (V_{\ka_c}(\g)^{\g[[t]]}) \simeq (\on{gr}
V_{\ka_c}(\g))^{\g[[t]]}.
$$
In particular, this implies that the character of $\on{gr} \zz(\ghat)
= \on{gr} (V_{\ka_c}(\g))^{\g[[t]]}$ is equal to that of $(\on{gr}
V_{\ka_c}(\g)^{\g[[t]]})$ given by formula \eqref{character of
  vac}. Thus we find that
\begin{equation}    \label{char of z}
\on{ch} \zz(\ghat) = \prod_{i=1}^\ell \prod_{n_i
  \geq d_i+1} (1-q^{n_i})^{-1}.
\end{equation}

\subsection{The center and the classical $\W$--algebra}
  \label{center and w}

According to \propref{center}, $\zz(\ghat)$ is contained in the
intersection of the kernels of the operators $\ol{V}_i[1],
i=1,\ldots,\ell$, on $\pi_0$. Now we compute the character of this
intersection and compare it with the character formula \eqref{char of
z} for $\zz(\ghat)$ to show that $\zz(\ghat)$ is equal to the
intersection of the kernels of the operators $\ol{V}_i[1]$.

First, we identify this intersection with a classical limit of a
one-parameter family of vertex algebras, called the
$\W$--algebras. The $\W$--algebra $\W_\nu(\g)$ associated to a simple
Lie algebra $\g$ and an invariant inner product $\nu$ on $\g$ is
defined via the quantum Drinfeld-Sokolov reduction in \cite{vertex},
Ch. 14. It is shown there that $\W_\nu(\g)$ is a vertex subalgebra of
$\pi_0^\nu$, which for generic values of $\nu$ is equal to the
intersection of the kernels of screening operators.

More precisely, consider another copy of the Heisenberg Lie algebra
$\hh_\nu$ introduced in \secref{pi0}.  To avoid confusion, we will
denote the generators of this Heisenberg Lie algebra by ${\mb
b}_{i,n}$. We have the vertex operator $V^\nu_{-\al_i}(z): \pi^\nu_0
\to \pi^\nu_{-\al_i}$ defined by formula \eqref{vertex op}, and let
$V^\nu_{-\al_i}[1] = \int V^\nu_{-\al_i}(z) dz$ be its residue.  We
call it the $\W$--algebra screening operator (to distinguish it from
the Kac-Moody screening operators defined above).  Since
$V^\nu_{-\al_i}(z) = Y_{\pi^\nu_0,\pi^\nu_{-\al_i}}(z)$ (in the
notation of \secref{generalities}), we obtain that the intersection of
the kernels of $V^\nu_{-\al_i}[1], i=1,\ldots,\ell$, in $\pi^\nu_0$ is
a vertex subalgebra of $\pi^\nu_0$. By Theorem 14.4.12 of
\cite{vertex}, for generic values of $\nu$ the $\W$--algebra
$\W_\nu(\g)$ is isomorphic to the intersection of the kernels of the
operators $V^\nu_{-\al_i}[1], i=1,\ldots,\ell$, in $\pi^\nu_0$.

Now we consider the limit when $\nu \to \infty$. We fix an invariant
inner product $\nu_0$ on $\g$ and denote by $\ep$ the ratio between
$\nu_0$ and $\nu$. We have the following formula for the $i$th simple
root $\al_i \in \h^*$ as an element of $\h$ using the identification
between $\h^*$ and $\h$ induced by $\nu = \nu_0/\ep$:
$$
\al_i = \ep \frac{2}{\nu_0(h_i,h_i)} h_i.
$$
Let
\begin{equation}    \label{b prime}
{\mb b}'_{i,n} = \ep \frac{2}{\nu_0(h_i,h_i)} {\mb b}_{i,n}.
\end{equation}
Consider the $\C[\ep]$--lattice in $\pi^{\nu}_0 \otimes \C[\ep]$
spanned by all monomials in ${\mb b}'_{i,n}, i =1,\ldots,\ell;
n<0$. We denote by $\pi_{0,\nu_0}$ the specialization of this lattice
at $\ep=0$; it is a commutative vertex algebra. In the limit $\ep \to
0$, we obtain the following expansion of the operator
$V^\nu_{-\al_i}[1]$:
$$
V^\nu_{-\al_i}[1] = \ep \frac{2}{\nu_0(h_i,h_i)} {\mb V}_i[1] +
\ldots,
$$
where the dots denote terms of higher order in $\ep$, and the operator
${\mb V}_i[1]$ acting on $\pi_{0,\nu_0}$ is given by the formula
\begin{equation}    \label{mbV}
{\mb V}_i[1] = \sum_{m\leq 0} {\mb V}_i[m] \; {\mb D}_{{\mb
    b}'_{i,m-1}},
\end{equation}
where
$$
{\mb D}_{{\mb b}'_{i,m}} \cdot b_{j,n} = a_{ij}
\frac{\pa}{\pa {\mb b}'_{i,n}} \delta_{n,m},
$$
$(a_{ij})$ is the Cartan matrix of $\g$, and 
$$
\sum_{n\leq 0} {\mb V}_i[n] z^{-n} = \exp \left( - \sum_{m>0}
\frac{{\mb b}'_{i,-m}}{m} z^m \right).
$$
The intersection of the kernels of the operators ${\mb V}_i[1],
i=1,\ldots,\ell$, is a commutative vertex subalgebra of
$\pi_{0,\nu_0}$, which we denote by ${\mb W}_{\nu_0}(\g)$ and call the
{\em classical $\W$--algebra} associated to $\g$ and $\nu_0$. Note
that the structures of commutative vertex algebra on $\pi_{0,\nu_0}$
and ${\mb W}_{\nu_0}(\g)$ and the operators ${\mb V}_i[1]$ are
independent of the choice of $\nu_0$. However, as we will see below,
both $\pi_{0,\nu_0}$ and ${\mb W}_{\nu_0}(\g)$ also carry vertex
Poisson algebra structures, and those structures do depend on $\nu_0$.

Now let us look at the $\W$--algebra ${\mb W}_{\nu_0}({}^L\g)$,
attached to the Langlands dual Lie algebra $^L\g$ whose Cartan matrix
is the transpose of the Cartan matrix of $\g$. Comparing formulas
\eqref{olVi} and \eqref{mbV}, we find that if we make the substitution
$b_{i,n} \mapsto - {\mb b}'_{i,n}$, then operator $\ol{V}_i[1]$
attached to $\g$ becomes the operator ${\mb V}_i[1]$ attached to
$^L\g$. Therefore the intersection of the kernels of the operators
$\ol{V}_i[1], i=1,\ldots,\ell$, attached to a simple Lie algebra $\g$
is isomorphic to the intersection of the kernels of the operators
${\mb V}_i[1], i=1,\ldots,\ell$, attached to $^L\g$, i.e., to ${\mb
W}_{\nu_0}({}^L\g)$.

We need to show that the character of ${\mb W}_{\nu_0}({}^L\g)$ is
given by the right hand side of formula \eqref{char of z}. This may be
done in several ways. One of the possible ways is presented in
\secref{ident with w and center} below. Another possibility is to use
the results of \cite{FF:laws}. According to \cite{FF:laws},
Prop. 2.2.8, the operators ${\mb V}_i[1], i=1,\ldots,\ell$, satisfy
the Serre relations of $\g$, and so they generate an action of the Lie
algebra $\n_+$ on $\pi_{0,\nu_0}$. Moreover, by \cite{FF:laws},
Prop. 2.4.6, the action of $\n_+$ on $\pi_{0,\nu_0}$ is co-free, and
this enables us to compute the character of ${\mb
W}_{\nu_0}({}^L\g)$. The result is that the character of ${\mb
W}_{\nu_0}({}^L\g)$ is given by the right hand side of formula
\eqref{char of z} (see \cite{FF:laws}, Prop. 2.4.7). Hence we obtain
that the character of the intersection of the kernels of the operators
$\ol{V}_i[1], i=1,\ldots,\ell$, is equal to the character of
$\zz(\ghat)$. Therefore $\zz(\ghat)$ is isomorphic to the intersection
of the kernels of the operators $\ol{V}_i[1], i=1,\ldots,\ell$. Thus,
we obtain the following result.

\begin{thm}    \label{isom with w}
The center $\zz(\ghat)$ is isomorphic to the intersection of the
kernels of the operators $\ol{V}_i[1], i=1,\ldots,\ell$, on
$\pi_{0,\nu_0}$, and to the classical $\W$--algebra ${\mb
  W}_{\nu_0}({}^L\g)$.

The natural map $\on{gr} \zz(\ghat) \to (\on{gr}
V_{\ka_c}(\g))^{\g[[t]]}$ is an isomorphism.
\end{thm}

The theorem implies in particular that for each generator $P_{i,0}$ of
$(\on{gr} V_{\ka_c}(\g))^{\g[[t]]}$ there exists an element $S_i$ of
$V_{\ka_c}(\g)^{\g[[t]]}$ whose symbol is equal to
$P_{i,0}$. Moreover, we have
$$
\zz(\ghat) = \C[S_{i}^{(n)}]_{i=1,\ldots,\ell;n\geq 0},
$$
where $S_{i}^{(n)} = T^n S_i$, and we use the commutative algebra
structure on $\zz(\ghat)$ which comes from its commutative vertex
algebra structure (see \cite{vertex}, \S~1.4). Note that the symbol of
$S_{i}^{(n)}$ is then equal to $n! P_{i,n}$.

For example, the vector $S_1$, which is unique up to a scalar, is
the Segal--Sugawara vector
\begin{equation}    \label{Segal}
S_1 = \frac{1}{2} \sum_a J^a_{-1} J_{a,-1} v_{\ka_c},
\end{equation}
where $\{ J_a \}$ is the basis of $\g$ dual to the basis $\{ J^a \}$
with respect to a non-degenerate invariant inner product. Explicit
formulas for other elements $S_i, i>1$, are unknown in general (see
however \cite{Ha,GW}).

We note that the center $\zz(\ghat)$ may be identified with the
algebra of $\ghat_{\ka_c}$--endomorphisms of $V_{\ka_c}(\g)$. Indeed,
a $\g[[t]]$--invariant vector $x \in V_{\ka_c}(\g)$ gives rise to a
non-trivial endomorphism of $V$, which maps the highest weight vector
$v_{\ka_c}$ to $x$. On the other hand, given an endomorphism $e$ of
$V_{\ka_c}(\g)$, we obtain a $\g[[t]]$--invariant vector $e(v)$. It is
clear that the two maps are inverse to each other. Now let $x_1$ and
$x_2$ be two $\g[[t]]$--invariant vectors in $V_{\ka_c}(\g)$. Let
$e_1$ and $e_2$ be the corresponding endomorphisms of
$V_{\ka_c}(\g)$. Then the image of $v$ under the composition $e_1 e_2$
equals $e_1(x_2) = x_2 x_1$, where in the right hand side we use the
product structure on $\zz(\ghat)$. Therefore we find that as an
algebra $\zz(\ghat)$ is isomorphic to the algebra
$\on{End}_{\ghat_{\ka_c}} V_{\ka_c}(\g)$ with the opposite
multiplication. But since $\zz(\ghat)$ is commutative, we obtain an
isomorphism of algebras
$$
\zz(\ghat) \simeq \on{End}_{\ghat_{\ka_c}} V_{\ka_c}(\g).
$$

\subsection{The vertex Poisson algebra structure}

In addition to the structures of commutative vertex algebras, both
$\zz(\ghat)$ and ${\mb W}_{\nu_0}(\g)$ also carry the structures of
vertex Poisson algebra. For the definition of vertex Poisson algebra,
see, e.g., \cite{vertex}, \S~15.2. In particular, we recall that to a
vertex Poisson algebra $P$ we attach a Lie algebra $$\on{Lie}(P) = P
\otimes \C((t))/\on{Im}(T \otimes 1 + 1 \otimes \pa_t)$$ (see
\cite{vertex}, \S~15.1.7).

According to Proposition 15.2.7 of \cite{vertex}, if $V_\ep$ is a
one-parameter family of vertex algebras, then the center ${\mc
  Z}(V_0)$ of $V_0$ acquires a natural vertex Poisson algebra
structure. Namely, at $\ep=0$ the polar part of the operation $Y$,
restricted to ${\mc Z}(V_0)$, vanishes, so we define the operation
$Y_-$ on ${\mc Z}(V_0)$ as the $\ep$--linear term of the polar part of
$Y$.

Let us fix a non-zero inner product $\ka_0$ on $\g$, and let $\ep$ be
the ratio between the inner products $\ka-\ka_c$ and $\ka_0$. Consider
the vertex algebras $V_\ka(\g)$ as a one-parameter family using $\ep$
as a parameter. Then we obtain a vertex Poisson structure on
$\zz(\ghat)$, the center of $V_{\ka_c}(\g)$ (corresponding to $\ep=0$).
We will denote $\zz(\ghat)$, equipped with this vertex Poisson
structure, by $\zz(\ghat)_{\ka_0}$.

Likewise, we obtain a vertex Poisson structure on $\pi_{0,\ka_0}$ by
considering the vertex algebras $\pi_0^{\ka-\ka_c}$ as a one-parameter
family with the same parameter $\ep$, in the same way as in
\secref{center and w}. The corresponding vertex Poisson structure on
$\pi_{0,\ka_0}$ is uniquely characterized by the Poisson brackets
$$
\{ b_{i,n},b_{j,m} \} = n \ka_0(h_i,h_j) \delta_{n,-m}.
$$

Recall that the homomorphism $\ww_{\ka_c}: V_{\ka_c}(\g) \to
W_{0,\ka_c} = M \otimes \pi_0$ may be deformed to a homomorphism
$\ww_\ka: V_\ka(\g) \to W_{0,\ka} = M \otimes
\pi_0^{\ka-\ka_c}$. Therefore the $\ep$--linear term of the polar part
of the operation $Y$ of $V_\ka(\g)$, restricted to $\zz(\ghat)$, which
is used in the definition of the vertex Poisson structure on
$\zz(\ghat)$, may be computed by restricting to $\zz(\ghat) \subset
\pi_0$ the $\ep$--linear term of the polar part of the operation $Y$
of $\pi_0^{\ka-\ka_c} \subset W_{0,\ka}$. But the latter defines the
vertex Poisson structure on $\pi_{0,\ka_0}$. Therefore we obtain the
following

\begin{lem}    \label{emb vpa}
  The embedding $\zz(\ghat)_{\ka_0} \hookrightarrow \pi_{0,\ka_0}$ is
  a homomorphism of vertex Poisson algebras and the corresponding map
  of local Lie algebras $\on{Lie}(\zz(\ghat)_{\ka_0}) \hookrightarrow
  \on{Lie}(\pi_{0,\ka_0})$ is a Lie algebra homomorphism.
\end{lem}

On the other hand, it follows from the definition of the classical
$\W$--algebra ${\mb W}_{\nu_0}(^L\g)$ that it also carries a vertex
Poisson algebra structure. The isomorphism of \thmref{isom with w}
preserves the vertex Poisson algebra structures in the following
sense.

Note that the restriction of a non-zero invariant
inner product $\ka_0$ on $\g$ to $\h$ defines a non-zero inner product
on $\h^*$, which is the restriction of an invariant inner product
$\ka^\vee_0$ on $^L\g$. To avoid confusion, let us denote by
$\pi_0(\g)_{\ka_0}$ and $\pi_0({}^L\g)_{\ka_0^\vee}$ the classical
limits of the Heisenberg vertex algebras for $\g$ and $^L\g$,
respectively, with their vertex Poisson structures corresponding to
$\ka_0$ and $\ka_0^\vee$, respectively. We have an isomorphism of
vertex Poisson algebras
\begin{align*}
  \imath: \pi_0(\g)_{\ka_0} & \overset{\sim}\to
  \pi_0({}^L\g)_{\ka_0^\vee}, \\
  b_{i,n} \in \hh & \mapsto - {\mb b}'_{i,n} \in \wh{^L\h},
\end{align*}
where ${\mb b}'_{i,n}$ is given by formula \eqref{b prime}. The
restriction of the isomorphism $\imath$ to the subspace
$\bigcap_{1\leq i\leq \ell} \on{Ker} \ol{V}_i[1]$ of
$\pi_0(\g)_{\ka_0}$ gives us an isomorphism of vertex Poisson algebras
$$
\bigcap_{1\leq i\leq \ell} \on{Ker} \ol{V}_i[1] \simeq \bigcap_{1\leq
  i\leq \ell} \on{Ker} {\mb V}^\vee_i[1],
$$
where by ${\mb V}^\vee_i[1]$ we denote the operator \eqref{mbV}
attached to $^L\g$. Recall that we have the following isomorphisms of
vertex Poisson algebras:
\begin{align*}
\zz(\ghat) & \simeq \bigcap_{1\leq i\leq \ell} \on{Ker} \ol{V}_i[1],
\\ {\mb W}_{\ka_0^\vee}({}^L\g) & \simeq \bigcap_{1\leq i\leq \ell}
\on{Ker} {\mb V}^\vee_i[1].
\end{align*}

Therefore we obtain the following strengthening of \thmref{isom with
  w}:

\begin{thm}    \label{odin}
The center $\zz(\ghat)_{\ka_0}$ is isomorphic, as a vertex Poisson
algebra, to the classical $\W$--algebra ${\mb W}_{\ka_0^\vee}({}^L\g)$.
The Lie algebras $\on{Lie}(\zz(\ghat)_{\ka_0})$ and
$\on{Lie}({\mb W}_{\ka_0^\vee}({}^L\g))$ are also isomorphic.
\end{thm} 

\subsection{$\AutO$--module structures}

Both $\zz(\ghat)_{\ka_0}$ and ${\mb W}_{\ka_0^\vee}({}^L\g)$ carry
actions of the group $\AutO$ and the isomorphism of \thmref{odin}
intertwines these actions. To see that, we describe the two actions as
coming from the vertex Poisson algebra structures.

In both cases the action of the group $\AutO$ is obtained by
exponentiation of the action of the Lie algebra $\DerO$. In the case
of the center $\zz(\ghat)_{\ka_0}$, the action of $\DerO$ is the
restriction of the natural action on $V_{\ka_c}(\g)$ which comes from
its action on $\ghat_{\ka_c}$ (and preserving its Lie subalgebra
$\g[[t]]$) by infinitesimal changes of variables. But away from the
critical level, i.e., when $\ka\neq \ka_c$, the action of $\DerO$ is
obtained thtough the action of the Virasoro algebra which comes from
the conformal vector ${\mb s}_\ka$ given by formula \eqref{formula for
Ska} which we rewrite as follows:
$${\mb s}_\ka = \frac{\ka_0}{\ka-\ka_c} S_1,$$ where $S_1$ is given by
formula \eqref{Segal} and $\ka_0$ is the inner product used in that
formula. Thus, the Fourier coefficients ${\mb s}_{\ka,(n)}, n\geq -1$,
of the vertex operator
$$
Y({\mb s}_\ka,z) = \sum_{n \in \Z} {\mb s}_{\ka,(n)} z^{-n-1}
$$
generate the $\DerO$--action on $V_\ka(\g)$ when $\ka \neq \ka_c$.

In the limit $\ka \to \ka_c$, we have ${\mb s}_\ka = \ep^{-1} S_1$,
where as before $\ep = \dfrac{\ka-\ka_c}{\ka_0}$. Therefore the action
of $\DerO$ is obtained through the vertex Poisson operation $Y_-$
(defined as the limit of $\ep^{-1}$ times the polar part of $Y$ when
$\ep \to 0$, see \cite{vertex}, \S~15.2) on $\zz(\ghat)$ applied to
$S_1 \in \zz(\ghat)$. In other words, the $\DerO$--action is generated
by the Fourier coefficients $S_{1,(n)}, n\geq 0$, of the series
$$
Y_-(S_1,z) = \sum_{n\geq 0} S_{1,(n)} z^{-n-1},
$$
namely, $L_n \mapsto S_{(n)}$. Thus, we see that the natural
$\DerO$--action (and hence $\AutO$--action) on $\zz(\ghat)$ is encoded
in the vector $S_1 \in \zz(\ghat)$ through the vertex Poisson algebra
structure on $\zz(\ghat)$. Note that this action endows $\zz(\ghat)$
with a quasi-conformal structure (see \secref{quasi-conf str}).

Likewise, there is a $\DerO$--action on ${\mb W}_{\ka_0^\vee}({}^L\g)$
coming from its vertex Poisson algebra structure. The vector
generating this action is also equal to the limit of a conformal
vector in the quantum $\W$--algebra $\W_\nu({}^L\g)$ as $\nu \to
\infty$. This conformal vector is unique because as we see from the
character formula for $\W_\nu({}^L\g)$ given in the right hand side of
\eqref{char of z} the homogeneous component of $\W_\nu({}^L\g)$ of
degree two (where all conformal vectors live) is one-dimensional. The
limit of this vector as $\nu \to \infty$ gives rise to a vector ${\mb
t}$ in ${\mb W}_{\ka_0^\vee}({}^L\g)$ such that the Fourier
coefficients of $Y_-(t,z)$ generate a $\DerO$--action on ${\mb
W}_{\ka_0^\vee}({}^L\g)$ (in the next section we will explain the
geometric meaning of this action). Since such a vector is unique, it
must be the image of $S_1 \in \zz(\ghat)$ under the isomorphism
$\zz(\ghat) \simeq {\mb W}_{\ka_0^\vee}({}^L\g)$. This implies the
following result.

\begin{prop}    \label{intertw}
The isomorphism \thmref{odin} intertwines the actions of $\DerO$ and
$\AutO$ on $\zz(\ghat)$ and ${\mb W}_{\ka_0^\vee}({}^L\g)$.
\end{prop}

We have already calculated in formula \eqref{formula for Ska} the
image of ${\mb s}_\ka$ in $W_{0,\ka}$ when $\ka \neq \ka_c$. By passing
to the limit $\ka \to \ka_c$ we find that the image of $S_1 = \ep
{\mb s}_\ka$ belongs to $\pi_0(\g)_{\ka_0} \subset W_{0,\ka}$ and is
equal to
$$
\frac{1}{2} \sum_{i=1}^\ell b_{i,-1} b^i_{-1} - \rho_{-2},
$$
where $\{ b_i \}$ and $\{ b^i \}$ are dual bases with respect to the
inner product $\ka_0$ used in the definition of $S_1$, restricted to
$\h$, and $\rho$ is the element of $\h$ corresponding to $\rho \in
\h^*$ under the isomorphism $\h^* \simeq \h$ induced by $\ka_0$. Under
the isomorphism of \thmref{odin}, this vector becomes the vector ${\mb
t} \in \pi_0({}^L\g)_{\ka_0^\vee} \simeq \pi_0(\g)_{\ka_0}$, which
automatically lies in ${\mb W}_{\ka_0^\vee}({}^L\g) \subset
\pi_0({}^L\g)_{\ka_0^\vee}$ and is responsible for the $\DerO$ action
on it.

The action of the corresponding operators $L_n \in \DerO, n \geq -1$,
on $\pi_0(\g)_{\ka_0}$ is given by derivations of the algebra
structure which are uniquely determined by formulas \eqref{action of
Ln}. Therefore the action of the $L_n$'s on ${\mb
W}_{\ka_0^\vee}({}^L\g)$ is as follows (recall that $b_{i,n} \mapsto
- {\mb b}'_{i,n}$):
\begin{align} \notag
L_n \cdot {\mb b}'_{i,m} &= - m {\mb b}'_{i,n+m}, \qquad -1\leq n<-m, \\
\label{new action of Ln}
L_n \cdot {\mb b}'_{i,-n} &= -n, \qquad n>0, \\ \notag
L_n \cdot {\mb b}'_{i,m} &= 0, \qquad n>-m.
\end{align}

\section{Opers and Miura opers}

In this section we introduce opers, (generic) Miura opers and a
natural map between them. We will then show in the next section that
the corresponding rings of functions are isomorphic to ${\mb
W}_{\nu_0}(\g)$ and $\pi_0(\g)_{\nu_0}$, respectively, and the
corresponding homomorphism ${\mb W}_{\nu_0}(\g) \to \pi_0(\g)_{\nu_0}$
is the embedding constructed in \secref{center and w}. In other words,
the ring of functions on opers is isomorphic to the subring of the
ring of functions on generic Miura opers equal to the intersection of
the kernels of the screening operators ${\mb V}_i[1],
i=1,\ldots,\ell$.

\subsection{Opers}    \label{opers}

Let $G$ be a simple algebraic group of adjoint type, $B$ a Borel
subgroup and $N = [B,B]$ its unipotent radical, with the corresponding
Lie algebras $\n \subset \bb\subset \g$. There is an open $B$--orbit
${\bf O}\subset \n^\perp/\bb \subset \g/\bb$, consisting of vectors
which are stabilized by the radical $N\subset B$, and such that all of
their negative simple root components, with respect to the adjoint
action of $H = B/N$, are non-zero. This orbit may also be described as
the $B$--orbit of the sum of the projections of simple root generators
$f_i$ of any nilpotent subalgebra $\n_-$, which is in generic position
with $\bb$, onto $\g/\bb$. The torus $H = B/N$ acts simply
transitively on ${\bf O}$, so ${\bf O}$ is an $H$--torsor.

Suppose we are given a principal $G$--bundle $\F$ on a smooth curve
$X$, or on a disc $D = \on{Spec} \OO$, or on a punctured disc $D^\times
= \on{Spec} \K$ (here we use the notation of \secref{coord-indep
version}), together with a connection $\nabla$ (automatically flat)
and a reduction $\F_B$ to the Borel subgroup $B$ of $G$. Then we
define the relative position of $\nabla$ and $\F_B$ (i.e., the failure
of $\nabla$ to preserve $\F_B$) as follows. Locally, choose any flat
connection $\nabla'$ on $\F$ preserving $\F_B$, and take the
difference $\nabla - \nabla'$.  It is easy to show that the resulting
local sections of $(\g/\bb)_{\F_B} \otimes \Omega$ are independent of
$\nabla'$, and define a global $(\g/\bb)_{\F_B}$--valued one-form on
$X$, denoted by $\nabla/\F_B$.

\medskip

Let $X$ be as above. A $G$--{\em oper} on $X$ is by definition a
triple $(\F,\nabla,\F_B)$, where $\F$ is a principal $G$--bundle $\F$
on $X$, $\nabla$ is a connection on $\F$ and $\F_B$ is a
$B$--reduction of $\F$, such that the one--form $\nabla/\F_B$ takes
values in ${\bf O}_{\F_B} \subset(\g/\bb)_{\F_B}$.

\medskip

This definition is due to A. Beilinson and V. Drinfeld \cite{BD} (in
the case when $X$ is the punctured disc opers were introduced in
\cite{DS}). Note that ${\bf O}$ is $\C^\times$--invariant, so that
${\bf O} \otimes \Omega$ is a well-defined subset of $(\g/\bb)_{\F_B}
\otimes \Omega$.

Equivalently, the above condition may be reformulated as saying that
if we choose a local trivialization of $\F_B$ and a local coordinate
$t$ then the connection will be of the form
\begin{equation}    \label{form of nabla}
\nabla = \pa_t + \sum_{i=1}^\ell \psi_i(t) f_i + {\mb v}(t),
\end{equation}
where each $\psi_i(t)$ is a nowhere vanishing function, and ${\mb
v}(t)$ is a $\bb$--valued function. If we change the trivialization of
$\F_B$, then this operator will get transformed by the corresponding
gauge transformation. This observation allows us to describe opers on
a disc $D = \on{Spec} \OO$ in a more explicit way.

Let us choose a coordinate $t$ on $D$, i.e., an isomorphism $\OO
\simeq \C[[t]]$. Then the space $\on{Op}_G(D)$ of $G$--opers on $D$ is
the quotient of the space of all operators of the form \eqref{form of
nabla}, where $\psi_i(t) \in \C[[t]], \psi_i(0) \neq 0,
i=1,\ldots,\ell$, and ${\mb v}(t) \in \bb[[t]]$, by the action of
the group $B[[t]]$ by gauge transformations:
$$
g \cdot (\pa_t + A(t)) = \pa_t + g A(t) g^{-1} - g^{-1}\pa_t g.
$$

Since the $B$--orbit ${\bf O}$ is an $H$--torsor, we can use the
$H$--action to make all functions $\psi_i(t)$ equal to $1$. Thus, we
obtain that $\on{Op}_G(D)$ is equal to the quotient of the space
$\wt{\on{Op}}_G(D)$ of operators of the form
\begin{equation}    \label{another form of nabla}
\nabla = \pa_t + \sum_{i=1}^\ell f_i + {\mb v}(t), \qquad {\mb v}(t) \in
\bb[[t]],
\end{equation}
by the action of the group $N[[t]]$ by gauge transformations. 

If we choose another coordinate $s$, such that $t = \varphi(s)$, then
the operator \eqref{another form of nabla} will become
$$
\nabla = \pa_s + \varphi'(s) \sum_{i=1}^\ell f_i + \varphi'(s) \cdot
       {\mb v}(\varphi(s)).
$$
In order to bring it back to the form \eqref{another form of nabla} we
need to apply the gauge transformation by $\rho^\vee(\varphi'(s))$,
where $\rho^\vee: \C^\times \to H$ is the one-parameter subgroup of
$H$ equal to the sum of the fundamental coweights of $G$. Choose a
splitting $\imath: H \to B$ of the homomorphism $B \to H$. Then,
considering $\rho^\vee$ as an element of the Lie algebra $\h =
\on{Lie} H$, we have $[\rho^\vee,e_i] = e_i$ and $[\rho^\vee,f_i] =
-f_i$. Therefore we find that
$$
\rho^\vee(\varphi'(s)) \cdot \left( \pa_s + \varphi'(s)
\sum_{i=1}^\ell f_i + \varphi'(s) \cdot {\mb v}(\varphi(s))\right)
$$
\begin{equation}    \label{change of var}
= \pa_s + \sum_{i=1}^\ell f_i + \varphi'(s) \rho^\vee(\varphi'(s))
  \cdot {\mb v}(\varphi(s)) \cdot \rho^\vee(\varphi'(s))^{-1}.
\end{equation}
Thus we obtain a well-defined action of the group $\AutO$ on the space
$\on{Op}_G(\Db)$ of opers on the standard disc, $$\Db = \on{Spec}
\C[[t]].$$ We may therefore define $\on{Op}_G(D)$ as the twist of
$\on{Op}_G(\Db)$ by the $\AutO$--torsor ${\mc A}ut$ (see
\secref{coord-indep version}).

In particular, the above formulas imply the following result. Consider
the $H$--bundle $\Omega^{\rho^\vee}$ on $D$,
defined in the same way as in \secref{trans form for fields} (where we
considered the $^L H$--bundle $\Omega^{-\rho}$).

\begin{lem}    \label{FH}
The $H$--bundle $\F_H = \F_B \underset{B}\times H = \F_B/N$ is
isomorphic to $\Omega^{\rho^\vee}$.
\end{lem}

\begin{proof}
It follows from formula \eqref{change of var} for the action of the
changes of variables on opers that if we pass from a coordinate $t$ on
$D$ to the coordinate $s$ such that $t = \varphi(s)$, then we obtain a
new trivialization of the $H$--bundle $\F_H$, which is related to the
old one by the transformation $\rho^\vee(\varphi'(s))$. This precisely
means that $\F_H \simeq \Omega^{\rho^\vee}$.
\end{proof}

\subsection{Canonical representatives}    \label{canon}

In this section we find canonical representatives in the
$N[[t]]$--gauge classes of connections of the form \eqref{another form
  of nabla}.

\begin{lem}[\cite{DS}]    \label{free}
The action of $N[[t]]$ on $\wt{\on{Op}}_G(D)$ is free.
\end{lem}

\begin{proof} The operator $\on{ad} \rho^\vee$ defines the principal
  gradation on $\bb$, with respect to which we have a direct sum
decomposition $\bb = \oplus_{i\geq 0} \bb_i$. Let
$$
p_{-1} = \sum_{i=1}^\ell f_i.
$$

The operator $\on{ad} p_{-1}$ acts from $\bb_{i+1}$ to $\bb_{i}$
injectively for all $i\geq 0$. Hence we can find for each $i\geq 0$ a
subspace $V_i \subset \bb_i$, such that $\bb_i = [p_{-1},\bb_{i+1}]
\oplus V_i$. It is well-known that $V_i \neq 0$ if and only if $i$ is
an exponent of $\g$, and in that case $\dim V_i$ is equal to the
multiplicity of the exponent $i$. In particular, $V_0=0$. Let $V =
\oplus_{i\in E} V_i \subset \n$, where $E$ is the set of exponents
of $\g$ counted with multiplicity.

We claim that each element of $\pa_t + p_{-1} + {\mb v}(t) \in
\wt{\on{Op}}_G(D)$ can be uniquely represented in the form
\begin{equation}    \label{gauge}
\pa_t + p_{-1} + {\mb v}(t) = \exp \left( \on{ad} U \right) \cdot
\left( \pa_t + p_{-1} + {\mb c}(t) \right),
\end{equation}
where $U \in \n[[t]]$ and ${\mb c}(t) \in V[[t]]$. To see that,
we decompose with respect to the principal gradation: $U=\sum_{j\geq
0} U_j$, ${\mb v}(t)=\sum_{j\geq 0} {\mb v}_j(t)$, ${\mb c}(t) =
\sum_{j\in E} {\mb c}_j(t)$. Equating the homogeneous components of
degree $j$ in both sides of \eqref{gauge}, we obtain that ${\mb c}_i +
[U_{i+1},p_{-1}]$ is expressed in terms of ${\mb v}_i,{\mb c}_j, j<i$,
and $U_j, j\leq i$. The injectivity of $\on{ad} p_{-1}$ then allows us
to determine uniquely ${\mb c}_i$ and $U_{i+1}$. Hence $U$ and ${\mb
c}$ satisfying equation \eqref{gauge} may be found uniquely by
induction, and the lemma follows.
\end{proof}

There is a special choice of the transversal subspace $V = \oplus_{i
\in E} V_i$ defined in the proof of \lemref{free}. Namely, let $p_1$
be the unique element of $\n$, such that $\{ p_{-1},2\rv,p_1 \}$ is an
$\sw_2$--triple. Let $V_{\can} = \oplus_{i \in E} V_{\can,i}$ be the
space of $\on{ad} p_1$--invariants in $\n$. Then $p_1$ spans
$V_{\on{can},1}$. Let $p_j$ be a linear generator of $V_{\can,d_j}$
(if the multiplicity of $d_j$ is greater than one, which happens only
in the case $\g=D^{(1)}_{2n}, d_j=2n$, then we choose linearly
independent vectors in $V_{\on{can},d_j}$).

According to \lemref{free}, each $G$--oper may be represented by a
unique operator $\nabla = \pa_t + p_{-1} + {\mathbf v}(t)$, where
${\mathbf v}(t) \in V_{\can}[[t]]$, so that we can write
$${\mathbf v}(t) = \sum_{j=1}^\ell v_j(t) \cdot p_j.$$

Suppose now that $t=\varphi(s)$, where $s$ is another coordinate on
$D$ such that $t=\varphi(s)$. With respect to the new coordinate $s$,
$\nabla$ becomes equal to $\pa_s + \wt{{\mathbf v}}(s)$, where
$\wt{{\mathbf v}}(s)$ is expressed via ${\mathbf v}(t)$ and
$\varphi(s)$ as in formula \eqref{change of var}. By \lemref{free},
there exists a unique operator $\pa_s + p_{-1} + \ol{{\mathbf v}}(s)$
with $\ol{{\mathbf v}}(s) \in V_{\can}[[s]]$ and $g \in B[[s]]$, such
that
\begin{equation}    \label{cano}
\pa_s + p_{-1} + \ol{{\mathbf v}}(s) = g
\cdot \left( \pa_s + \wt{{\mathbf v}}(s) \right).
\end{equation}
It is straightforward to find that
\begin{align*}
g &= \exp \left(\frac{1}{2} \frac{\varphi''}{\varphi'} \cdot p_1
\right) \rho^\vee(\varphi), \\
\ol{v}_1(s) &= v_1(\varphi(s)) \left( \varphi'
\right)^2 - \frac{1}{2} \{ \varphi,s \}, \\ \ol{v}_j(s) &=
v_j(\varphi(s)) \left(  \varphi' \right)^{d_j+1},
\quad \quad j>1,
\end{align*}
where $$\{ \varphi,s \} = \frac{\varphi'''}{\varphi'} - \frac{3}{2}
\left( \frac{\varphi''}{\varphi'} \right)^2$$ is the Schwarzian
derivative.

These formulas mean that under changes of variables, $v_1$ transforms
as a projective connection, and $v_j, j>1$, transforms as a
$(d_j+1)$--differential on $D$. Thus, we obtain an isomorphism
\begin{equation}    \label{repr}
\Op \simeq  {\mc P}roj \times \bigoplus_{j=2}^\el
\omega^{\otimes(d_j+1)},
\end{equation}
where $\omega^{\otimes n}$ is the space of $n$--differentials on $D$
and ${\mc P}roj$ is the $\omega^{\otimes 2}$--torsor of projective
connections on $D$.

The above formulas also show that the $G$--bundle $\F$ is defined by
the transition function $U(s)$. In particular, the $G$--bundles $\F$
underlying all opers are isomorphic to each other.

\subsection{Miura opers}    \label{miura opers}

A {\em Miura $G$--oper} on $X$, which is a smooth curve, or $D$, or
$D^\times$, is by definition a quadruple $(\F,\nabla,\F_B,\F'_B)$,
where $(\F,\nabla,\F_B)$ is a $G$--oper on $X$ and $\F'_B$ is another
$B$--reduction of $\F$ which is preserved by $\nabla$.

\medskip

Consider the space of Miura $G$--opers on the disc $D$. A
$B$--reduction of $\F$ which is preserved by the connection $\nabla$
is uniquely determined by a $B$--reduction of the fiber $\F_0$ of $\F$
at the origin $0 \in D$ (recall that the underlying $G$--bundles of
all $G$--opers are isomorphic to each other). The set of such
reductions is the $\F_0$--twist $(G/B)_{\F_0}$ of the flag manifold
$G/B$. Therefore the natural forgetful morphism from the space of all
Miura opers on $D$ to $\Op$ is a $G/B$--bundle.

A Miura $G$--oper is called {\em generic} if the $B$--reductions
$\F_B$ and $\F'_B$ are in generic relative position. We denote the
space of generic Miura opers on $D$ by $\MOp$. We have a natural
forgetful morphism $\MOp \to \Op$. The group $N_{\F_{B,0}}$, where
$\F_{B,0}$ is the fiber of $\F_B$ at $0$, acts on $(G/B)_{\F_0}$, and
the subset of generic reductions is the open $N_{\F_{B,0}}$--orbit of
$(G/B)_{\F_0}$. This orbit is in fact an
$N_{\F_{B,0}}$--torsor. Therefore we obtain that the morphism $\MOp
\to \Op$ is a principal $N_{\F_{B,0}}$--bundle. We denote it by $\mu$
and call it the {\em Miura transformation}.

Now we identify $\MOp$ with the space of $H$--connections. Consider
the $H$--bundles $\F_H = \F_B/N$ and $\F'_H = \F'_B/N$.

\begin{lem}    \label{H bundles isom}
The $H$--bundles $\F_H$ and $\F'_H$ are isomorphic.
\end{lem}

\begin{proof}
Consider the vector bundles $\g_\F = \F \underset{G}\times \g$,
$\bb_{\F_B} = \F_B \underset{B}\times \bb$ and $\bb_{\F'_B} = \F'_B
\underset{G}\times \bb$. We have the inclusions $\bb_{\F_B},
\bb_{\F'_B} \subset \g_\F$ which are in generic position. Therefore
the intersection $b_{\F_B} \cap b_{\F'_B}$ is isomorphic to
$\bb_{\F_B}/[\bb_{\F_B},\bb_{\F_B}]$, which is the trivial vector
bundle with the fiber $\h$. It naturally acts on the bundle $\g_\F$
and under this action $\g_\F$ decomposes into a direct sum of $\h$ and
the line subbundles $\g_{F,\al}, \al \in \De$. Furthermore,
$\bb_{\F_B} = \bigoplus_{\al \in \De_+} \g_{F,\al}, \bb_{\F'_B} =
\bigoplus_{\al \in \De_-} \g_{F,\al}$. Since the action of $B$ on
$\n/[\n,\n]$ factors through $H = B/N$, we find that
$$
\F_H \underset{H}\times \bigoplus_{i=1}^\ell \C_{\al_i} \simeq
\bigoplus_{i=1}^\ell \g_{\F,\al_i}, \qquad \F'_H \underset{H}\times
\bigoplus_{i=1}^\ell \C_{-\al_i} \simeq \bigoplus_{i=1}^\ell \g_{\F,-\al_i}.
$$
Clearly, the line bundle $\g_{\F,\al}$ is dual to
$\g_{\F,-\al}$. Therefore we obtain that
$$
\F_H \underset{H}\times \C_{\al_i} \simeq \F'_H
\underset{H}\times \C_{\al_i}, \qquad i=1,\ldots,\ell.
$$
Since $G$ is of adjoint type by our assumption, the above associated
line bundles completely determine $\F_H$ and $\F'_H$. Therefore
$\F_H \simeq \F'_H$.
\end{proof}

Combining \lemref{H bundles isom} with \lemref{FH}, we find that
$\F'_H \simeq \Omega^{\rho^\vee}$. Since the $B$--bundle $\F'_B$ is
preserved by the oper connection $\nabla$, we obtain a connection
$\ol{\nabla}$ on $\F'_H \simeq \Omega^{\rho^\vee}$. Therefore we
obtain a morphism $\beta$ from the space $\MOp$ of generic Miura opers
on $D$ to the space of connections $\on{Conn}(\Omega^{\rho^\vee})_D$
on the $H$--bundle $\Omega^{\rho^\vee}$ on $D$.

\begin{prop}    \label{map beta}
The map $\beta: \MOp \to \on{Conn}(\Omega^{\rho^\vee})_D$ is an
isomorphism.
\end{prop}

\begin{proof}
We define a map $\tau$ in the opposite direction. Suppose we are given
a connection $\ol\nabla$ on the $H$--bundle $\Omega^{\rho^\vee}$ on
$D$. We associate to it a generic Miura oper as follows. Let us choose
a splitting $H \to B$ of the homomorphism $B \to H$ and set $\F =
\Omega^{\rho^\vee} \underset{H}\times G, \F_B = \Omega^{\rho^\vee}
\underset{H}\times B$, where we consider the adjoint action of $H$ on
$G$ and on $B$ obtained through the above splitting. The choice of the
splitting also gives us the opposite Borel subgroup $B_-$, which is
the unique Borel subgroup in generic position with $B$ containing
$H$. Then we set $\F'_B = \Omega^{\rho^\vee} \underset{H}\times B_-$.

Observe that the space of connections on $\F$ is isomorphic to the
direct product
$$
\on{Conn}(\Omega^{\rho^\vee})_D \times \bigoplus_{\al \in \De}
\omega^{\al(\rho^\vee) + 1}.
$$
Its subspace corresponding to negative simple roots is isomorphic to
$\left( \bigoplus_{i=1}^\ell \g_{-\al_i} \right) \otimes \OO$. Having
chosen a basis element $f_i$ of $\g_{-\al_i}$ for each
$i=1,\ldots,\ell$, we now construct an element $p_{-1} =
\sum_{i=1}^\ell f_i$ of this space. Now we set $\nabla = \ol\nabla +
p_{-1}$. By construction, $\nabla$ has the correct relative position
with the $B$--reduction $\F_B$ and preserves the $B$--reduction
$\F'_B$. Therefore the quadruple $(\F,\nabla,\F_B,\F'_B)$ is a generic
Miura oper on $D$. We set $\tau(\ol{\nabla}) =
(\F,\nabla,\F_B,\F'_B)$.

This map is independent of the choice of splitting $H \to B$ and of
the generators $f_i, i=1,\ldots,\ell$. Indeed, changing the splitting
$H \to B$ amounts to a conjugation of the old splitting by an element
of $N$. This is equivalent to applying to $\nabla$ the gauge
transformation by this element. Therefore it will not change the
underlying Miura oper structure. Likewise, rescaling of the generators
$f_i$ may be achieved by a gauge transformation by a constant element
of $H$, and this again does not change the Miura oper structure.  It
is clear from the construction that $\beta$ and $\tau$ are mutually
inverse isomorphisms.
\end{proof}

\section{Identification with the $\W$--algebra and the center}

The Miura transformation $\mu: \MOp \to \Op$ gives rise to a
homomorphism of the corresponding rings of functions $\wt{\mu}:
\on{Fun} \Op \to \on{Fun} \MOp$. We will now identify $\on{Fun} \MOp$
with $\pi_0(\g)_{\nu_0}$ and the image of $\wt{\mu}$ with the
intersection of kernels of the screening operators. This will give us
an identification of $\on{Fun} \Op$ with ${\mb W}_{\nu_0}(\g)$.

\subsection{Screening operators}    \label{oper scr}

Since $\mu$ is an $N_{\F_{B,0}}$--bundle, we obtain that the image of
the homomorphism $\wt\mu$ is equal to the space of
$N_{\F_{B,0}}$--invariants of $\on{Fun} \MOp$, and hence to the space
of $\n_{\F_{B,0}}$--invariants of $\on{Fun} \MOp$. We fix a Cartan
subalgebra $\h$ in $\bb$, a coordinate $t$ on the disc $D$ and a
trivialization of $\F_B$. Then we identify $\on{Fun} \Op$ and
$\on{Fun} \MOp$ with rings of polynomials in infinitely many
variables. Each space has an action of $\DerO$ and $\wt\mu$ is a
$\DerO$--equivariant homomorphism between these rings. We will
describe the image of $\wt\mu$ as the intersection of the kernels of
certain linear operators which will turn out to be nothing but the
screening operators used in the definition of the classical
$\W$--algebra.

Using $t$ and our trivialization of $\F_B$, we write each oper in the
form
$$
\pa_t + p_{-1} + \sum_{i=1}^\el v_i(t) \cdot {\mb c}_i, \quad \quad
v_i(t) \in \C[[t]]
$$
(see \secref{canon}), where
$$
v_i(t) = \sum_{n<-d_i} v_{i,n} t^{-n-d_i-1}.
$$
Thus,
$$
\on{Fun} \Op = \C[v_{i,n_i}]_{i=1,\ldots,\ell;n_i<-d_i}.
$$
The action of $\DerO$ on this ring is determined from formulas
\eqref{cano}.

Likewise, we write each generic Miura oper as
\begin{equation}    \label{Miura oper}
\pa_t + p_{-1} + {\mathbf u}(t), \quad \quad {\mathbf u}(t) \in \h[[t]].
\end{equation}
Set $u_i(t) = \al_i({\mathbf u}(t)), i=1,\ldots,\el$, and
$$
u_i(t) = \sum_{n<0} u_{i,n} t^{-n-1}.
$$
Then
$$
\on{Fun} \MOp = \C[u_{i,n}]_{i=1,\ldots,\ell; n<0}.
$$
The action of $\DerO$ on this ring is determined as follows. If $s$ is
a new coordinate on $D$ such that $t = \varphi(s)$, then the same
Miura oper will appear as $\pa_s + p_{-1} + \wt{\mb u}(s)$, where
\begin{equation}    \label{express1}
\wt{\mb u}(s) = \varphi' \cdot {\mb u}(\varphi(s)) - \rho^\vee
\frac{\varphi''}{\varphi'}.
\end{equation}
This translates into the following formulas for the action of the
generators $L_n = -t^{n+1} \pa_t$ on the $u_{i,m}$'s:
\begin{align} \notag
L_n \cdot u_{i,m} &= - m u_{i,n+m}, \qquad -1\leq n<-m, \\
\label{Miura oper action of Ln}
L_n \cdot u_{i,-n} &= -n, \qquad n>0, \\ \notag
L_n \cdot u_{i,m} &= 0, \qquad n>-m.
\end{align}

\medskip

\noindent {\em Example.} We compute the Miura transformation $\mu$ in
the case when $\g=\sw_2$. In this case an oper has the form
$$
\pa_t + \left( \begin{matrix} 0 & v(t) \\ 1 & 0 \end{matrix} \right),
$$
and a generic Miura oper has the form
$$
\pa_t + \left( \begin{matrix} \frac{1}{2} u(t) & 0 \\ 1 &
- \frac{1}{2} u(t) \end{matrix} \right).
$$
To compute $\mu$, we need to find a element of $N[[t]]$ such that the
corresponding gauge transformation brings the Miura oper into the oper
form. We find that
$$
\left( \begin{matrix} 1 & - \frac{1}{2} u(t) \\ 0 & 1 \end{matrix}
\right) \left( \pa_t + \left( \begin{matrix} \frac{1}{2} u(t) & 0 \\ 1
& -\frac{1}{2} u(t) \end{matrix} \right) \right) \left( \begin{matrix}
  1 & \frac{1}{2} u(t) \\ 0 & 1 \end{matrix} \right)= \pa_t + \left(
\begin{matrix} 0 & \frac{1}{4} u(t)^2 + \frac{1}{2} \pa_t u(t) \\ 1 &
  0 \end{matrix} \right).
$$
Therefore we obtain
$$
\mu(u(t)) = v(t) = \frac{1}{4} u(t)^2 + \frac{1}{2} \pa_t u(t),
$$
which may also be written in the form
$$
\pa_t^2 - v(t) = \left( \pa_t + \frac{1}{2} u(t) \right) \left( \pa_t -
\frac{1}{2} u(t) \right).
$$
It is this transformation that was originally introduced by R. Miura
as the map intertwining the flows of the KdV hierarchy and the mKdV
hierarchy.

In the case of $\sw_2$ opers are projective connections and Miura
opers are affine connections. The Miura transformation is nothing but
the natural map from affine connections to projective connections.\qed

\medskip

Having chosen a trivialization of $\F_B$, we identify the twist
$\n_{\F_{B,0}}$ with $\n$. Now we choose the generators $e_i,
i=1,\ldots,\ell$, of $\n$ with respect to the action of $\h$ on $\n$
in such a way that for each $i$ the generator $e_i$ and the previously
chosen $f_i$ satisfy the standard relations of $\sw_2$. The
$N_{\F_{B,0}}$--action on $\MOp$ then gives rise to an infinitesimal
action of $e_i$ on $\MOp$. We will now compute the corresponding
derivation on $\on{Fun} \MOp$.

The action of $e_i$ is given by the infinitesimal gauge transformation
\begin{equation}    \label{delta}
\delta {\mathbf u}(t) = [x_i(t) \cdot e_i,\pa_t + p_{-1} + {\mathbf
u}(t)],
\end{equation}
where $x_i(t) \in \C[[t]]$ is such that $x_i(0) = 1$, and the right
hand side of formula \eqref{delta} belongs to $\h[[t]]$. It turns out
that these conditions determine $x_i(t)$ uniquely. Indeed, the right
hand side of \eqref{delta} reads:
$$
x_i(t) \cdot \al_i^\vee - u_i(t) x_i(t) \cdot e_i - \pa_t x_i(t) \cdot
e_i.
$$
Therefore it belongs to $\h[[t]]$ if and only if
\begin{equation}    \label{rec}
\pa_t x_i(t) = - u_i(t) x_i(t).
\end{equation}
If we write
$$
x_i(t) = \sum_{n\leq 0} x_{i,n} t^{-n},
$$
and substitute it into formula \eqref{rec}, we obtain that the
coefficients $x_{i,n}$ satisfy the following recurrence relation:
$$
n x_{i,n} = \sum_{k+m=n;k<0;m\leq 0} u_{i,k} x_{i,m}, \qquad n<0.
$$
We find from this formula that
\begin{equation}    \label{xin}
\sum_{n\leq 0} x_{i,n} t^{-n} = \exp \left( - \sum_{m>0}
\frac{u_{i,-m}}{m} t^m \right).
\end{equation}

Now we obtain that
$$
\delta u_j(t) = \al_j(\delta {\mathbf u}(t)) = a_{ij} x_i(t),
$$
where $(a_{ij})$ is the Cartan matrix of $\g$. In other words, the
operator $e_i$ acts on the algebra $\on{Fun} \MOp = \C[u_{i,n}]$ by
the derivation
\begin{equation}    \label{si}
\sum_{j=1}^\el a_{ij} \sum_{n\geq 0} x_{i,n} \frac{\pa}{\pa
  u_{i,-n-1}},
\end{equation}
where $x_{i,n}$ are given by formula \eqref{xin}. Thus, we
obtain the following characterization of $\on{Fun} \Op$ inside
$\on{Fun} \MOp$.

\begin{prop}
The image of $\on{Fun} \Op$ in $\on{Fun} \MOp$ under the Miura map
$\wt\mu$ is equal to the intersection of kernels of the operators
given by formula \eqref{si} for $i=1,\ldots,\ell$.
\end{prop}

\subsection{Back to the $\W$--algebra}
\label{ident with w and center}

Comparing formula \eqref{si} with formula \eqref{mbV} we find that if
we replace ${\mb b}'_{i,n}$ by $u_{i,n}$ in formula \eqref{mbV} for
the screening operator ${\mb V}_i[1]$, then we obtain formula
\eqref{si}. Therefore the intersection of the kernels of the operators
\eqref{si} is equal to the intersection of the kernels of the
operators ${\mb V}_i[1], i=1,\ldots,\ell$. But the latter is the
classical $\W$--algebra ${\mb W}_{\nu_0}(\g)$. Hence we obtain a
commutative diagram
\begin{equation}    \label{comm diag}
\begin{CD}
\pi_0(\g)_{\nu_0} @>{\sim}>> \on{Fun} \on{MOp}_G(\Db) \\
@AAA @AAA \\
{\mb W}_{\nu_0}(\g) @>{\sim}>> \on{Fun} \on{Op}_G(\Db)
\end{CD}
\end{equation}
where the top arrow is an isomorphism of algebras given on generators
by the assignment ${\mb b}'_{i,n} \mapsto u_{i,n}$.

In particular, we obtain from the above isomorphisms vertex Poisson
algebra structures on $\on{Fun} \on{Op}_G(\Db)$ and $\on{Fun}
\on{Op}_G(\Db)$ (the latter structure may alternatively be defined by
means of the Drinfeld-Sokolov reduction, as we show below). As before,
these structures depend on the choice of inner product $\nu_0$, which
from now on we will use as a subscript.

The above diagram gives us a geometric interpretation of the screening
operators ${\mb V}_i[1]$: they correspond to the action of the
generators $e_i$ of $\n$ on $\on{Fun} \on{MOp}_G(\Db)$.

As a byproduct, we obtain a character formula for ${\mb
W}_{\nu_0}(\g)$, as we promised in \secref{center and w}. Indeed, the
action of $N$ on $\on{MOp}_G(\Db)$ is free, so that $$H^i(\n,\on{Fun}
\on{MOp}_G(\Db)) = 0, \qquad i>0.$$ Therefore the character of
$\on{Fun} \on{Op}_G(\Db)$ is equal to the character of the Chavalley
complex computing the cohomology of $\n$ with coefficients in
$\on{Fun} \on{MOp}_G(\Db)$. The latter is equal to the character of
$\on{Fun} \on{MOp}_G(\Db)$, i.e., $\prod_{n>0} (1-q^n)^{-\ell}$, times
the character of $\bigwedge^\bullet \n^*$. But since $\deg e_i = -1$,
we obtain that that $\deg e_\al = -\al(\rho^\vee)$, and so
$$
\on{ch} \bigwedge{}^\bullet \n^* = \prod_{\al \in \De_+}
(1-q^{\al(\rho^\vee)}) = \prod_{i=1}^\ell \prod_{n_i=1}^{d_i}
(1-q^{n_i}).
$$
Thus, we obtain that the character of $\on{Fun} \on{Op}_G(\Db) \simeq
{\mb W}_{\nu_0}(\g)$ is given by the right hand side of formula
\eqref{char of z}.

Comparing formulas \eqref{Miura oper action of Ln} and \eqref{new
action of Ln}, we find that the top isomorphism in \eqref{comm diag}
is $\DerO$--equivariant. Since the screening operators commute with
the $\DerO$--action, we find that the bottom isomorphism is also
$\DerO$--equivariant. Thus we obtain the following

\begin{thm}    \label{w same as opers}
The classical $\W$--algebra ${\mb W}_{\nu_0}(\g)$ is
$\DerO$--equivariantly isomorphic to the algebra $\on{Fun}
\on{Op}_G(\Db)$ of functions on the space of opers on
$\Db = \on{Spec} \C[[t]]$.
\end{thm}

Combining \thmref{w same as opers} with \thmref{odin} and
\propref{intertw} we come to the following result, where by $^L G$ we
understand the group of inner automorphisms of $^L\g$.

\begin{thm}    \label{isom with opers}
There is an isomorphism $\zz(\ghat)_{\ka_0} \simeq \on{Fun}
\on{Op}_{^L G}(\Db)_{\ka_0^\vee}$ which preserves the vertex Poisson
structures and the $\DerO$--module structures on both sides. Moreover,
it fits into a commutative diagram of vertex Poisson algebras equipped
with $\DerO$--action:
\begin{equation}    \label{important}
\begin{CD}
\pi_0(\g)_{\ka_0} @>{\sim}>> \on{Fun} \on{MOp}_{^L
  G}(\Db)_{\ka^\vee_0} \\
@AAA @AAA \\
\zz(\ghat)_{\ka_0} @>{\sim}>> \on{Fun} \on{Op}_{^L G}(\Db)_{\ka_0^\vee}
\end{CD}
\end{equation}
where the upper arrow is given on the generators by the assignment
$b_{i,n} \mapsto - u_{-i,n}$.
\end{thm}

This theorem is consistent with the previously given descriptions of
the transformation properties of the $b_i(t)$'s and the
$u_i(t)$'s. Indeed, according to the computations of \secref{trans
form for fields}, the $b_i(t)$'s transform as components of a
connection on the $^L H$--bundle $\Omega^{-\rho}$ on the standard disc
$\Db = \on{Spec} \C[[t]]$. On the other hand, by formula \eqref{change
of var}, the $u_i(t)$'s transform as components of a connection on the
dual $^L H$--bundle $\Omega^{\rho}$ on $\Db$. The map $b_i(t) \mapsto
- u_i(t)$ sends the former to the latter.

Given any disc $D$, we may consider the twists of our algebras by the
$\AutO$--torsor $\Au$, defined as in \secref{coord-indep version}. We
will mark them by the subscript $D$. Then we obtain that
\begin{align*}
{\mb W}_{\nu_0}(\g)_D &\simeq \on{Fun} \Op_{\nu_0}, \\
\zz(\ghat)_{\ka_0,D} &\simeq \on{Fun} \on{Op}_{^L G}(D)_{\ka_0^\vee}.
\end{align*}

\subsection{The associated graded spaces}

In this section we describe the isomorphism $\zz(\ghat) \simeq
\on{Fun} \on{Op}_{^L G}(\Db)$ of \thmref{isom with opers} at the level
of associated graded spaces. Note that both algebras $\zz(\ghat)$ and
$\on{Fun} \on{Op}_{^L G}(\Db)$ have natural filtrations.

The filtration on $\zz(\ghat)$ is induced by the
Poincar\'e--Birkhoff--Witt filtration on the universal enveloping
algebra $U(\ghat_{\ka_c})$, see \secref{ass graded}. By \thmref{isom
  with w}, $\on{gr} \zz(\ghat) = (\on{gr}
V_{\ka_c}(\g))^{\g[[t]]}$. But
$$
\on{gr} V_{\ka_c}(\g) = \on{Sym} \g((t))/\g[[t]] \simeq \on{Fun}
\g^*[[t]] dt,
$$
independently of the choice of coordinate $t$ and inner product on
$\g$. In \propref{free pol} we gave a description of $(\on{gr}
V_{\ka_c}(\g))^{\g[[t]]}$. The coordinate-independent version of this
description is as follows. Let $C^\vee_\g = \on{Spec} \; (\on{Fun}
\g)^G$ and
$$
C^\vee_{\g,\Omega} = \Omega \underset{\C^\times}\times C^\vee_\g,
$$
where $\Omega = \C[[t]] dt$ is the topological module of differentials
on $\Db = \on{Spec} \C[[t]]$. Note that a choice of homogeneous
generators $P_i, i=1,\ldots,\ell$, of $(\on{Fun} \g)^G$ gives us an
identification $$C^\vee_{\g,\Omega} \simeq \bigoplus_{i=1}^\ell
\Omega^{\otimes(d_i+1)}.$$ By \propref{free pol} we have a canonical
isomorphism
\begin{equation}    \label{can isom}
\on{gr} \zz(\ghat) \simeq \on{Fun} C^\vee_{\g,\Omega}.
\end{equation}

Next we consider the map $$\al: \zz(\ghat) \to \pi_0(\g)$$ which is
equal to the restriction of the embedding $V_{\ka_c}(\g) \to
W_{0,\ka_c}$ to $\zz(\ghat)$. In the proof of \propref{injective} we
described a filtration on $W_{0,\ka_c}$ compatible with the PBW
filtration on $V_{\ka_c}(\g)$. This implies that the map $\zz(\ghat)
\to \pi_0(\g)$ is also compatible with filtrations. According to
\secref{trans form for fields},
$$
\pi_0(\g) = \on{Fun} \on{Conn}(\Omega^{-\rho})_{\Db}.
$$
The space $\on{Conn}(\Omega^{-\rho})_{\Db}$ of
connections on the $^L H$--bundle $\Omega^{-\rho}$ is an affine space
over the vector space $^L\h \otimes \Omega = \h^* \otimes \Omega$.
Therefore we find that
\begin{equation}    \label{gr of pi}
\on{gr} \pi_0(\g) = \on{Fun} \h^* \otimes \Omega.
\end{equation}

We have the Harish-Chandra isomorphism $(\on{Fun} \g^*)^G \simeq
(\on{Fun} \h^*)^W$, where $W$ is the Weyl group of $\g$, and hence an
embedding $(\on{Fun} \g^*)^G \to \on{Fun} \h^*$. This embedding gives
rise to an embedding
$$
\on{Fun} C^\vee_{\g,\Omega} \to \on{Fun} \h^* \otimes \Omega.
$$
It follows from the proof of \propref{injective} that this is
precisely the map $$\on{gr} \al: \on{gr} \zz(\ghat) \to \on{gr}
\pi_0(\g)$$ under the identifications \eqref{can isom} and \eqref{gr of
pi}.

Let us now look at the associated graded of the map
\begin{equation}    \label{from Miura opers to opers}
\on{Fun} \on{Op}_{^L G}(\Db) \to \on{Fun} \on{MOp}_{^L G}(\Db).
\end{equation}
Let $C_\g = (\on{Fun} \g)^G$ and
$$
C_{\g,\Omega} = \Omega \underset{\C^\times}\times C_\g.
$$
Using the canonical form of the opers described in \secref{canon},
we obtain an identification
\begin{equation}    \label{gr of opers}
\on{gr} \on{Fun} \on{Op}_{^L G}(\Db) \simeq \on{Fun}
C_{{}^L\g,\Omega}
\end{equation}
(see \cite{BD}, \S\S~3.1.12--3.1.14, for more details). On the other
hand, \propref{map beta} implies that the space $\on{MOp}_{^L G}(\Db)$
is isomorphic to the space of connections
$\on{Conn}(\Omega^{\rho^\vee})_{\Db}$, which is an affine space over
the vector space $^L\h \otimes \Omega$. Therefore
\begin{equation}    \label{gr of Miura opers}
\on{gr} \on{Fun} \on{MOp}_{^L G}(\Db) \simeq \on{Fun} {}^L\h \otimes
\Omega.
\end{equation}

We have the Harish-Chandra isomorphism $(\on{Fun} \g)^G \simeq
(\on{Fun} \h)^W$, where $W$ is the Weyl group of $\g$, and hence an
embedding $(\on{Fun} \g)^G \to \on{Fun} \h$. This embedding gives
rise to an embedding
$$
\on{Fun} C_{\g,\Omega} \to \on{Fun} \h \otimes \Omega.
$$
The explicit construction of the morphism
$\on{Conn}(\Omega^{\rho^\vee})_{\Db} \to \on{MOp}_{^L G}(\Db)$ given
in the proof of \propref{map beta} implies that the associated graded
of the map \eqref{from Miura opers to opers}, with respect to the
identifications \eqref{gr of opers} and \eqref{gr of Miura opers} is
nothing but the homomorphism
$$
\on{Fun} C_{{}^L\g,\Omega} \to \on{Fun} {}^L\h \otimes \Omega.
$$

Now recall that our isomorphism $\zz(\ghat) \simeq \on{Fun}
\on{Op}_{^L G}(\Db)$ fits into the commutative diagram
\eqref{important}. All maps in this diagram preserve the filtrations
and we have described above the maps of the associated graded induced
by the vertical arrows in \eqref{important}. Moreover, we have
$$
\on{gr} \pi_0(\g) \simeq \on{Fun} \h^* \otimes \Omega = \on{Fun} {}^L\h
\otimes \Omega \simeq \on{gr} \on{MOp}_{^L G}(\Db),
$$
and
$$
\on{gr} \zz(\ghat) \simeq \on{Fun} C^\vee_{\g,\Omega} = \on{Fun}
C_{{}^L\g,\Omega} \simeq \on{gr} \on{Fun} \on{Op}_{^L G}(\Db),
$$
where we used the canonical isomorphisms
$$
(\on{Fun} \g^*)^G = (\on{Fun} \h^*)^W = (\on{Fun} {}^L\h)^W =
(\on{Fun} {}^L \g)^{^L G}.
$$
According to \thmref{isom with opers}, under these identifications
the map
$$
\on{gr} \pi_0(\g) \to \on{gr} \on{MOp}_{^L G}(\Db)
$$
induced by the upper vertical arrow in the diagram \eqref{important}
is equal to the the operator $(-1)^{\on{deg}}$ which takes the value
$(-1)^n$ on elements of degree $n$. Therefore we obtain that the same
is true for the associated graded of our isomorphism
\begin{equation}    \label{isom again}
\zz(\ghat) \simeq \on{Fun} \on{Op}_{^L G}(\Db).
\end{equation}
Thus we obtain the following

\begin{thm}
The isomorphism \eqref{isom again} preserves filtrations. The
corresponding associated graded spaces are both isomorphic to
$\on{Fun} C^\vee_{\g,\Omega}$. The corresponding isomorphism of the
associated graded spaces is equal to $(-1)^{\on{deg}}$.
\end{thm}

\section{The center of the completed universal enveloping algebra of
  $\ghat_{\ka_c}$}

Let $U_{\ka}(\ghat)$ be the quotient of the universal enveloping
algebra $U(\ghat_{\ka})$ of $\ghat_{\ka}$ by the ideal generated
by $(K-1)$. Define its completion $\wt{U}_{\ka}(\ghat)$ as follows:
$$
\wt{U}_{\ka}(\ghat) = \underset{\longleftarrow}\lim \;
U_{\ka}(\ghat)/U_{\ka}(\ghat) \cdot (\g \otimes t^N\C[[t]]).
$$
It is clear that $\wt{U}_{\ka}(\ghat)$ is a topological algebra which
acts on all smooth representations of $\ghat_{\ka}$ (i.e., such that
any vector is annihilated by $\g \otimes t^N\C[[t]]$ for sufficiently
large $N$), on which $K$ acts as the identity. In this section we
describe the center $Z(\ghat)$ of $\wt{U}_{\ka_c}(\ghat)$ and show
that it is isomorphic to the topological algebra of functions on the
space $\on{Op}_{^L G}(\Db^\times)$ of $^L G$--opers on the punctured
standard disc $\Db$.

\subsection{Enveloping algebra of a vertex algebra}
\label{enveloping}

Recall from \cite{vertex}, \S~3.5, that to any vertex algebra $V$ we
may attach a topological Lie algebra
$$
U(V) = (V \otimes \C((t)))/\on{Im}(T \otimes 1 + 1 \otimes \pa_t).
$$
It is topologically spanned by elements $A_{[n]} = A \otimes t^n, A
\in V, n \in \Z$, with the commutation relations
\begin{equation}    \label{comm rel}
[A_{[m]},B_{[k]}] = \sum_{n \geq 0} \left( \begin{matrix} m \\ n
\end{matrix} \right) (A_{(n)} \cdot B)_{[m+k-n]}.
\end{equation}
Now let $V = V_{\ka}(\g)$. Then each $A \in V_{\ka}(\g)$ is a
linear combination of monomials in $J^a_n, n<0$, applied to the vacuum
vector. We construct a linear map $U(V_{\ka}(\g)) \to
\wt{U}_{\ka}(\ghat)$ sending $(J^{a_1}_{n_1} \ldots
J^{a_m}_{n_m} v_{\ka})_{[k]}$ to
\begin{multline*}
\on{Res}_{z=0} Y(J^{a_1}_{n_1} \ldots J^{a_m}_{n_m} v_{\ka},z) z^k dz
\\ = \frac{1}{(-n_1-1)!} \ldots \frac{1}{(-n_m-1)!} \on{Res}_{z=0}
\Wick \pa_z^{-n_1-1} J^{a_1}(z) \ldots \pa_z^{-n_m-1} J^{a_m}(z) \Wick
\; z^k dz,
\end{multline*}
where as before $J^a(z) = \sum_{n \in \Z} J^a_n z^{-n-1}$. The
commutation relations between Fourier coefficients of vertex operators
(see \cite{vertex}, \S~3.3.6) imply that this map is a homomorphism of
Lie algebras. This homomorphism is clearly injective. Moreover, the
universal enveloping algebra of $U(V_{\ka}(\g))$ also injects into
$\wt{U}_{\ka}(\ghat)$ and its completion is equal to
$\wt{U}_{\ka}(\ghat)$. This motivates the following general
construction of the enveloping algebra of a vertex algebra.

Let $V$ be a vertex algebra and $U(V)$ the corresponding Lie algebra
defined as above. Denote by $U(U(V))$ the universal enveloping algebra
of the Lie algebra $U(V)$. Define its completion
$$
\wt{U}(U(V)) = \underset{\longleftarrow}\lim \; U(U(V))/I_N,
$$
where $I_N$ is the left ideal generated by $A_{[n]}, A \in V,
n>N$. Then $\wt{U}(V)$ is by definition the quotient of $\wt{U}(U(V))$
by the two-sided ideal generated by the Fourier coefficients of the
series
$$
Y[A_{(-1)} B,z] - \Wick Y[A,z] Y[B,w] \Wick, \qquad A,B \in V,
$$
where $Y[A,z] = \sum_{n \in \Z} A_{[n]}
z^{-n-1}$, and the normal ordering is defined in the same way as for
the vertex operators $Y(A,z)$. In the case when $V = V_{\ka}(\g)$, we
have $\wt{U}(V_{\ka}(\g)) = \wt{U}_{\ka}(\ghat)$.

Clearly, the assignment $V \mapsto \wt{U}(V)$ gives rise to a functor
from the category of vertex algebras to the category of topological
associative algebras.

\begin{remark}
A. Beilinson and V. Drinfeld give a more conceptual construction of
the topological algebra $\wt{U}(V)$ in \cite{BD},
\S~3.7.1--3.7.3. Namely, they define $\wt{U}(V)$ as the space of
horizontal sections of the left $\D$--module $\V_{\Db^\times}$ on the
punctured disc $\Db^\times$ defined by a generalization of the
construction given in \secref{coord-indep version}.\qed
\end{remark}

\subsection{From $\zz(\ghat)$ to the center}

Let $B \in \zz(\ghat) \subset V_{\ka_c}(\g)$. Then $\g[[t]] \cdot B =
0$ and formula \eqref{comm rel} for the commutation relations implies
that all Fourier coefficients $B_{[k]}$ of the vertex operator
$Y(B,z)$ commute with the entire affine algebra $\ghat_{\ka_c}$. Thus
we obtain an injective map $U(\zz(\ghat)) \to Z(\ghat)$. Moreover,
applying the enveloping algebra functor $\wt{U}$ to $\zz(\ghat)$ we
obtain an injective homomorphism $\wt{U}(\zz(\ghat)) \to Z(\ghat)$.

The algebra $\wt{U}(\zz(\ghat))$ is a completion of a polynomial
algebra. Indeed, recall from \secref{center and w} that
$$
\zz(\ghat) = \C[S_{i}^{(n)}]_{i=1,\ldots,\ell;n\geq 0},
$$
where $S_{i}^{(n)} = T^n S_i$. Let $S_{i,[n]}, n \in \Z$, be the
Fourier coefficients of the vertex operator $Y(S_i,z)$, considered as
elements of $\wt{U}_{\ka_c}(\ghat)$. Then $\wt{U}(\zz(\ghat))$ is the
completion of the polynomial algebra $\C[S_{i,[n]}]_{i=1,\ldots,\ell;n
\in \Z}$ with respect to the topology in which the basis of open
neighborhoods of $0$ is formed by the ideals generated by
$S_{i,[n]}, i=1,\ldots,\ell; n>N$, for $N \in \Z$. Thus,  the
completed polynomial algebra $\wt{U}(\zz(\ghat))$ embeds into the
center $Z(\ghat)$ of $\wt{U}_{\ka_c}(\ghat)$.

\begin{prop}
$Z(\ghat)$ is equal to $\wt{U}(\zz(\ghat))$.
\end{prop}

The proof is given in \cite{BD}, Theorem 3.7.7. It is based on the
description of the associated graded of $Z(\ghat)$ with respect to the
PBW filtration, which is obtained by an argument similar to the one
used in the proof of \propref{free pol}. We obtain then that the
associated graded of $Z(\ghat)$ is equal to the completed polynomial
algebra in $P_{i,n}, i=1,\ldots,\ell;n \in \Z$, which are the Fourier
coefficients of the series $P_i(J^a(z))$ (see \secref{ass graded} for
the definition of the $P_i$'s). But by construction $P_{i,n}$ is the
symbol of $S_{i,[-n-1]}$. Therefore $Z(\ghat)$ cannot contain anything
but elements of the completed polynomial algebra in the $S_{i,[n]}$'s,
hence the result.

\subsection{Identification of the center with $\on{Fun} \on{Op}_{^L
    G}(\Db^\times)$}

Recall the $\AutO$--equiva\-riant isomorphism of vertex Poisson algebras
\begin{equation}    \label{isom of vpa}
\zz(\ghat)_{\ka_0} \simeq \on{Fun} \on{Op}_{^L G}(\Db)_{\ka_0^\vee}
\end{equation}
established in \thmref{isom with opers}. For a vertex Poisson algebra
$P$, the space $U(P)$ (also denoted by $\on{Lie}(P)$) has a natural
Lie algebra structure (see \cite{vertex}, \S~15.1.7), and therefore
the commutative algebra $\wt{U}(P)$ carries a Poisson algebra
structure. In particular, $\wt{U}(\zz(\ghat)_{\ka_0})$ is a Poisson
algebra. Under the identification $\wt{U}(\zz(\ghat)_{\ka_0}) \simeq
Z(\ghat)$ the corresponding Poisson structure on $Z(\ghat)_{\ka_0}$
may be described as follows. Consider $\wt{U}_{\ka}(\ghat)$ as the
one-parameter family $A_\ep$ of associative algebras depending on the
parameter $\ep = (\ka-\ka_c)/\ka_0$. Then $Z(\ghat)$ is the center of
$A_0$. Given $x, y \in Z(\ghat)$, let $\wt{x}, \wt{y}$ be their
liftings to $A_\ep$. Then the Poisson bracket $\{ x,y \}$ is defined
as the $\ep$--linear term in the commutator $[\wt{x},\wt{y}]$
considered as a function of $\ep$ (it is independent of the choice of
the liftings). The fact that this is indeed a Poisson structure was
observed by V. Drinfeld following the work \cite{Ha} of T. Hayashi. We
denote the center $Z(\ghat)$ equipped with this Poisson structure by
$Z(\ghat)_{\ka_0}$.

Likewise, the vertex Poisson algebra $\on{Fun} \on{Op}_{^L
G}(\Db)_{\ka_0^\vee}$ gives rise to a topological Poisson algebra
$\wt{U}(\on{Fun} \on{Op}_{^L G}(\Db)_{\ka_0^\vee})$. Then the
isomorphism \eqref{isom of vpa} gives rise to an isomorphism of
Poisson algebras
$$
Z(\ghat)_{\ka_0} \simeq \wt{U}(\on{Fun} \on{Op}_{^L
  G}(\Db)_{\ka_0^\vee}).
$$
Moreover, this isomorphism is $\AutO$--equivariant.

\begin{lem}    \label{on punctured disc}
The algebra $\wt{U}(\on{Fun} \on{Op}_{^L G}(\Db)_{\ka_0^\vee})$ is
canonically isomorphic to the topological algebra $\on{Fun}
\on{Op}_{^L G}(\Db^\times)$ of functions on the space of $^L G$--opers
on the punctured disc $\Db^\times = \on{Spec} \C((t))$.
\end{lem}

\begin{proof}
As explained in \secref{oper scr}, the algebra $\on{Fun} \on{Op}_{^L
G}(\Db)$ is isomorphic to the polynomial algebra
$\C[v_{i,n_i}]_{i=1,\ldots,\ell;n_i<-d_i}$, where $v_{i,n}$ are the
coefficients of the oper connection
\begin{equation}    \label{operr}
\nabla = \pa_t + p_{-1} + \sum_{i=1}^\ell v_i(t) {\mb c}_i,
\end{equation}
where $v_i(t) = \sum_{n<-d_i} v_{i,n} t^{-n-d_i-1}$. Then the above
construction of the functor $\wt{U}$ implies that the algebra
$\wt{U}(\on{Op}_{^L G}(\Db)_{\ka_0^\vee})$ is the completion of the
polynomial algebra in the variables $v_{i,n}, i=1,\ldots,\ell; n \in
\Z$, with respect to the topology in which the base of open
neighborhoods of $0$ is formed by the ideals generated by $v_{i,n},
n<N$. But this is precisely the algebra of functions on the space of
opers of the form \eqref{operr}, where $v_i(t) = \sum_{n \in \Z}
v_{i,n} t^{-n-d_i-1} \in \C((t))$.
\end{proof}

Since $\on{Fun} \on{Op}_{^L G}(\Db)_{\ka_0^\vee}$ is a vertex Poisson
algebra, we obtain that $\on{Fun} \on{Op}_{^L G}(\Db^\times)$ is a
Poisson algebra. We will use the subscript $\ka_0^\vee$ to indicate
the dependence of this Poisson structure on the inner product
$\ka^\vee_0$ on $^L\g$. We will give another definition of this
Poisson structure in the next section. Thus we obtain the following
result which was originally conjectured by V. Drinfeld.

\begin{thm}    \label{center isom to opers}
The center $Z(\ghat)_{\ka_0}$ is isomorphic, as a Poisson algebra, to
the Poisson algebra $\on{Fun} \on{Op}_{^L G}(\Db^\times)_{\ka_0^\vee}$.
Moreover, this isomorphism is $\AutO$--equivariant.
\end{thm}

\subsection{The Poisson structure on $\on{Fun}
  \on{Op}_{G}(\Db^\times)_{\nu_0}$}

The Poisson algebra of functions on the space
$\on{Op}_{G}(\Db^\times)$ of opers on the punctured disc $\Db^\times$
may be obtained by hamiltonian reduction called the Drinfeld-Sokolov
reduction \cite{DS}.

We start with the Poisson manifold $\on{Conn}_\g$ of connections on
the trivial $G$--bundle on $\Db^\times$, i.e., operators of the form
$\nabla = \pa_t + A(t)$, where $A(t) \in \g((t))$. The Poisson
structure on this manifold comes from its identification with a
hyperplane in the dual space to the affine Kac--Moody algebra
$\ghat_{\nu_0}$. Indeed, the topological dual space to $\ghat_{\nu_0}$
may be identified with the space of all $\la$--connections on the
trivial bundle on $\Db^\times$, see \cite{vertex}, \S~15.4. Namely, we
split $\ghat_{\nu_0} = \g((t)) \oplus \C K$ as a vector space. Then a
$\lambda$--connection $\pa_t + A(t)$ gives rise to a linear functional
on $\ghat_{\nu_0}$ which takes value $$\la b + \on{Res}
\nu_0(A(t),B(t)) dt$$ on $(B(t) + b K) \in \g((t)) \oplus \C K =
\ghat_{\nu_0}$. Note that under this identification the action of
$\AutO$ by changes of coordinates and the coadjoint (resp., gauge)
action of $G((t))$ on the space of $\la$--connections (resp., the dual
space to $\ghat_{\nu_0}$) agree. The space $\on{Conn}_\g$ is now
identified with the hyperplane in $\ghat_{\nu_0}^*$ which consists of
those functionals which take value $1$ on the central element $K$.

The dual space $\ghat_{\nu_0}^*$ carries a canonical Poisson structure
called the Kirillov--Kostant structure (see, e.g., \cite{vertex},
\S~15.4.1 for more details). Because $K$ is a central element, this
Poisson structure restricts to the hyperplane which we have identified
with $\on{Conn}_\g$. Therefore we obtain a Poisson structure on the
$\on{Conn}_\g$.

The group $N((t))$ acts on $\on{Conn}_\g$ by gauge transformations.
This action corresponds to the coadjoint action of the $N((t))$ on
$\ghat_{\nu_0}^*$ and is hamiltonian, the moment map being the
surjection $m: \on{Conn}_\g \to \n((t))^*$ dual to the embedding
$\n((t)) \to \ghat_{\nu_0}$. We pick a one--point coadjoint
$N((t))$--orbit in $\n((t))^*$ represented by the linear functional
$\psi$ which is equal to the composition $\n((t)) \to \n/[\n,\n]((t))
= \bigoplus_{i=1}^\ell \g_{\al_i}((t))$ and the functional
$$(x_i(t))_{i=1}^\ell \mapsto \sum_{i=1}^\ell \on{Res} x_i(t) dt.$$
One shows in the same way as in the proof of \lemref{free} that the
action of $N((t))$ on $m^{-1}(\psi)$ is free. Moreover, the quotient
$m^{-1}(\psi)/N((t))$, which is the Poisson reduced manifold, is
canonically identified with the space of $G$--opers on
$\Db^\times$. Therefore we obtain a Poisson structure on the
topological algebra of functions $\on{Fun}
\on{Op}_{G}(\Db^\times)$. On the other hand we have defined a Poisson
algebra structure on $\on{Fun} \on{Op}_{G}(\Db^\times)_{\nu_0}$ using
its identification with $\wt{U}(\on{Fun} \on{Op}_{G}(\Db)_{\nu_0})$,
the isomorphism $\on{Fun} \on{Op}_{G}(\Db)_{\nu_0} \simeq {\mb
W}_{\nu_0}(\g)_{\nu_0}$ obtained in \thmref{w same as opers}, and the
vertex Poisson algebra structure on ${\mb W}_{\nu_0}(\g)$.

\begin{lem}
The two Poisson structures coincide.
\end{lem}

\begin{proof}
In \cite{vertex}, \S\S~14.4 and 15.8, we defined a complex
$C^\bullet_\infty(\g)$ and showed that its $0$th cohomology is a
vertex Poisson algebra canonically isomorphic to ${\mb
W}_{\nu_0}(\g)$ (and all other cohomologies vanish). It is
clear from the construction that if we apply the functor $\wt{U}$ to
the complex $C^\bullet_\infty(\g)$ we obtain the BRST complex of the
Drinfeld-Sokolov reduction described above. Since the functor $\wt{U}$
is left exact, we obtain that the $0$th cohomology of
$\wt{U}(C^\bullet_\infty(\g))$, i.e., the Poisson algebra $\on{Fun}
\on{Op}_{G}(\Db^\times)_{\nu_0}$ is isomorphic to the Poisson algebra
$\wt{U}({\mb W}_{\nu_0}(\g))$ which is what we needed to
prove.
\end{proof}

\subsection{The Miura transformation as the Harish-Chandra
  homomorphism}

The Harish-Chandra homomorphism is the homomorphism from the center
$Z(\g)$ of $U(\g)$, where $\g$ is a simple Lie algebra to the algebra
$\on{Fun} \h^*$ of polynomials on $\h^*$. It identifies $Z(\g)$ with
the algebra $(\on{Fun} \h^*)^W_\rho$ of $W$--invariant polynomials on
$\h^*$ shifted by $\rho \in \h^*$. To construct this homomorphism,
one needs to assign a central character to each $\la \in \h^*$. This
central character is just the character of the center on the Verma
module $M_\la$.

In the affine case, we construct a similar homomorphism from the
center $Z(\ghat)$ of $\wt{U}_{\ka_c}(\ghat)$ to the topological
algebra $\on{Fun} \on{Conn}(\Omega^{-\rho})_{\Db^\times}$ of functions
on the space of connections on the $^L H$--bundle $\Omega^{-\rho}$ on
$\Db^\times$. According \propref{param of wak}, points of
$\on{Conn}(\Omega^{-\rho})_{\Db^\times}$ parameterize Wakimoto modules
of critical level. Thus, for each $\chi \in
\on{Conn}(\Omega^{-\rho})_{\Db^\times}$ we have the Wakimoto module
$W_\chi$ of critical level. The following theorem describes the affine
analogue of the Harish-Chandra homomorphism.

Note that $\on{Fun} \on{Conn}(\Omega^{-\rho})_{\Db^\times}$ is the
completion of the polynomial algebra in $b_{i,n}, i=1,\ldots,\ell; n
\in \Z$, with respect to the topology in which the base of open
neighborhoods of $0$ is formed by the ideals generated by $b_{i,n},
n<N$. In the same way as in the proof of \lemref{on punctured disc} we
show that it is isomorphic to $\wt{U}(\on{Fun}
\on{Conn}(\Omega^{-\rho})_{\Db})$. We also define the topological
algebra $\on{Fun} \on{MOp}_{^L G}(\Db^\times)$ as $\wt{U}(\on{Fun}
\on{MOp}_{^L G}(\Db^\times))$. It is the completion of the polynomial
algebra in $u_{i,n}, i=1,\ldots,\ell; n \in \Z$, with respect to the
topology in which the base of open neighborhoods of $0$ is formed by
the ideals generated by $u_{i,n}, n<N$. The isomorphism of
\propref{map beta} gives rise to an isomorphism of the topological
algebras
$$
\on{Fun} \on{Conn}(\Omega^{-\rho})_{\Db^\times} \to \on{Fun}
  \on{MOp}_{^L G}(\Db^\times),
$$
under which $b_{i,n} \mapsto - u_{i,n}$.

\begin{thm}    \label{miura and hc}
The center $Z(\ghat)$ acts by a central character on each Wakimoto
module $W_\chi$. The corresponding homomorphism of algebras $Z(\ghat)
\to \on{Fun} \on{Conn}(\Omega^{-\rho})_{\Db^\times}$ fits into a
commutative diagram
\begin{equation}    \label{important1}
\begin{CD}
\on{Fun} \on{Conn}(\Omega^{-\rho})_{\Db^\times} @>{\sim}>> \on{Fun}
  \on{MOp}_{^L G}(\Db^\times) \\
@AAA @AAA \\
Z(\ghat) @>{\sim}>> \on{Fun} \on{Op}_{^L G}(\Db^\times)
\end{CD}
\end{equation}
where the upper arrow is given on the generators by the assignment
$b_{i,n} \mapsto - u_{i,n}$.
\end{thm}

\begin{proof}
The action of $U(V_{\ka_c}(\g))$, and hence of $\wt{U}_{\ka_c}(\ghat)$,
on $W_\chi$ is obtained through the homomorphism of vertex algebra
$V_{\ka_c}(\g) \to M_\g \otimes \pi_0(\g)$. In particular, the action
of $Z(\ghat)$ on $W_\chi$ is obtained through the homomorphism of
commutative vertex algebras $\zz(\ghat) \to \pi_0(\g)$. But $\pi_0(\g)
= \on{Fun} \on{Conn}(\Omega^{-\rho})_{\Db}$. Therefore the statement
of the theorem follows from \thmref{isom with opers} by applying the
functor of enveloping algebras $V \mapsto \wt{U}(V)$ introduced in
\secref{enveloping}.
\end{proof}

Thus, we see that the affine analogue of the Harish-Chandra
homomorphism is nothing but the Miura transformation for the Langlands
dual group. In particular, its image (which in the finite-dimensional
case consists of all $W$--invariant polynomials) is equal to the
intersection of the kernels of the screening operators.

\subsection{Compatibility with the finite-dimensional Harish-Chandra
homomorphism}

We have a natural $\Z$--gradation on $Z(\ghat)$ by the operator $L_0 =
- t \pa_t \in \DerO$. Let us denote by the upper script $0$ the
degree $0$ part and by the upper script $<0$ be the span of all
elements of negative degrees in any $\Z$--graded object. Consider the
Wakimoto modules $W_{\chi}$, where $\chi = \pa_t + \la/t, \la \in
\h^*$. These Wakimoto modules are $\Z$--graded, with the degree
$0$ part $W_\chi^0$ being the subspace $\C[a^*_{\al,0}\vac]_{\al \in
  \Delta_+} \subset M_\g \subset W_\chi$. The Lie algebra $\g$ acts on
this space and it follows immediately from our construction of
Wakimoto modules that as a $\g$--module, $W_\chi^0$ is isomorphic to
the contragredient Verma module $M^*_\la$. The algebra
$$
Z' = Z(\ghat)^0/(Z(\ghat) \cdot Z(\ghat)^{<0})^0
$$
naturally acts on $W_\chi^0$ and commutes with $\g$. Therefore varying
$\la \in \h^*$ we obtain a homomorphism
\begin{equation}    \label{map from z prime}
Z' \to (\on{Fun} \h^*)^W_\rho,
\end{equation}
which factors through the Harish-Chandra homomorphism $Z(\g) \to
(\on{Fun} \h^*)^W_\rho$.

Let $\on{Op}_{^L G}(\Db^\times)_{\leq 1}$ be the subvariety of {\em
opers with regular singularity} in $\on{Op}_{^L G}(\Db^\times)$, whose
points are $B((t))$--gauge equivalence classes of operators of the
form
\begin{equation}    \label{oper with sing}
\pa_t + \dfrac{1}{t}\left(p_{-1} + {\mb v}(t) \right), \qquad {\mb
  v}(t) \in {}^L \bb[[t]]
\end{equation}
(see \cite{BD}, \S~3.8.8).

Likewise we introduce the subvariety $\on{MOp}_{^L
G}(\Db^\times)_{\leq 1}$ of {\em Miura opers with regular singularity}
in $\on{MOp}_{^L G}(\Db^\times)$, whose points are operators of the
form
\begin{equation}    \label{miura oper with sing}
\pa_t + p_{-1} + \frac{1}{t} {\mb u}(t), \qquad {\mb
  u}(t) \in {}^L \h[[t]].
\end{equation}
Thus, $\on{MOp}_{^L G}(\Db^\times)_{\leq 1}$ is identified with
the space $\on{Conn}(\Omega^\rho)_{\leq 1}$ of connections with
regular singularity on the $^L H$--bundle $\Omega^\rho$ on $\Db$.
The Miura transformation then restricts to a morphism
$$
\on{MOp}_{^L G}(\Db^\times)_{\leq 1} \to \on{Op}_{^L
  G}(\Db^\times)_{\leq 1}.
$$

We have a $\Z$--gradation on both $\on{Fun} \on{Op}_{^L
G}(\Db^\times)$ and $\on{Fun} \on{Conn}(\Omega^\rho)$ induced by
the action of the operator $L_0 = - t \pa_t \in \DerO$. This gradation
descends to $\on{Fun} \on{Op}_{^L G}(\Db^\times)_{\leq 1}$ and
$\on{Fun} \on{MOp}_{^L G}(\Db^\times)_{\leq 1}$ where it takes only
non-negative values.

Define residue maps
\begin{align*}
\on{Res}: \; & \on{Op}_{^L G}(\Db^\times)_{\leq 1} \to
\on{Spec} \; (\on{Fun} {}^L\g)^{^L G}, \\
\on{Res}: \; & \on{Conn}(\Omega^\rho)_{\leq 1}) \to
\on{Spec} {}^L\h,
\end{align*}
sending an oper \eqref{oper with sing} to $p_{-1}+{\mb v}(0)$ (resp.,
a Miura oper \eqref{miura oper with sing} to ${\mb u}(0)$). They give
rise to homomorphisms
\begin{align*}
(\on{Fun} {}^L\g)^{^L G} & \to (\on{Fun} \on{Op}_{^L
    G}(\Db^\times)_{\leq 1})^0, \\ \on{Fun} {}^L\h & \to (\on{Fun}
    \on{Conn}(\Omega^\rho)_{\leq 1})^0.
\end{align*}

\begin{lem}    \label{another diagram}
These homomorphisms are isomorphisms and they fit into the commutative
diagram
$$
\begin{CD}
(\on{Fun} {}^L\h) @>{\sim}>> (\on{Fun}
\on{Conn}(\Omega^\rho)_{\leq 1})^0 \\
@AAA @AAA \\
(\on{Fun} {}^L\g)^{^L G} @>{\sim}>> (\on{Fun} \on{Op}_{^L
  G}(\Db^\times)_{\leq 1})^0
\end{CD}
$$
where the left vertical arrow is the Harish-Chandra homomorphism
shifted by $\rho$.
\end{lem}

\begin{proof}
Note that 
$$
\pa_t + p_{-1} + \frac{1}{t}{\mb u}(t) = \rho(t) \, \left(
\pa_t + \frac{p_{-1}}{t} + \frac{1}{t}({\mb u}(t) - \rho) \right) \,
\rho(t)^{-1}.
$$
This implies that the above diagram is commutative. It is clear from
the definition that the upper horizontal arrow is an isomorphism. To
see that the lower horizontal arrow is an isomorphism we pass to the
associated graded spaces where it becomes obvious (see \secref{char of
  zz}).
\end{proof}

Since the ideal of $\on{Op}_{^L G}(\Db^\times)_{\leq 1}$ in $\on{Fun}
\on{Op}_{^L G}(\Db^\times)$ is equal to $\on{Fun} \on{Op}_{^L
G}(\Db^\times)^{<0}$ we obtain that the isomorphism of \thmref{center
isom to opers} gives rise to an isomorphism of algebras
\begin{equation}    \label{the above}
Z' \simeq (\on{Fun} \on{Op}_{^L G}(\Db^\times)_{\leq 1})^0.
\end{equation}
Then \lemref{another diagram} implies that the map \eqref{map from z
prime} is an isomorphism. Recall that $^L\h = \h^*$.
Now \thmref{miura and hc} implies the following

\begin{prop}
There is a commutative diagram of isomorphisms
$$
\begin{CD}
Z' @>{\sim}>> (\on{Fun} \on{Op}_{^L G}(\Db^\times)_{\leq 1})^0 \\
@AAA @AAA \\
(\on{Fun} \h^*)^W_\rho @>{\sim}>> (\on{Fun} \h^*)^W_\rho
\end{CD}
$$
where the lower horizontal arrow is given by $f \mapsto f^-$,
$f^-(\la+\rho) = f(-\la-\rho)$.
\end{prop}

This is the statement of Theorem 3.8.17 of \cite{BD}.

\end{document}